\documentclass[11pt]{article}
\usepackage{amsfonts}
\usepackage{amssymb, amsmath}
\usepackage{enumitem}
\usepackage[mathscr]{eucal}
\usepackage{mathrsfs}
\usepackage{graphics}
\usepackage{setspace}
\usepackage{tensor}
\usepackage{graphicx}
\usepackage{verbatim}
\usepackage{cases}
\usepackage{float}
\usepackage{url}
\usepackage{setspace} 
\usepackage[font=small,format=plain,labelfont=bf,up,textfont=it,up]{caption}

\DeclareFontFamily{OMX}{yhex}{}
\DeclareFontShape{OMX}{yhex}{m}{n}{<->yhcmex10}{}
\DeclareSymbolFont{yhlargesymbols}{OMX}{yhex}{m}{n}
\DeclareMathAccent{\wideparen}{\mathord}{yhlargesymbols}{"F3}

    \newcommand{\argmin}{\mathop{\rm argmin}}

    \newcommand{\vech}{\mathrm{vech}}
    \newcommand{\wh}{\widehat}

    \newcommand{\n}{\mathbf{N}}

    \newcommand{\wt}{\widetilde}

\newtheorem{theorem}{Theorem}[section]
\newtheorem{lemma}{Lemma}[section]

\newtheorem{corollary}{Corollary}[section]

\numberwithin{equation}{section}
\newtheorem{remark}{Remark}[section]

\newenvironment{proof}[1][\sc Proof.]{\begin{trivlist}
\item[\hskip \labelsep {\bfseries #1}]}{\end{trivlist}}

\newcommand{\qed}{\nobreak \ifvmode \relax \else
      \ifdim\lastskip<1.5em \hskip-\lastskip
      \hskip1.5em plus0em minus0.5em \fi \nobreak
      \vrule height0.75em width0.5em depth0.25em\fi}

\allowdisplaybreaks

\makeatletter
\renewenvironment{titlepage}
 {%
  \if@twocolumn
    \@restonecoltrue\onecolumn
  \else
    \@restonecolfalse\newpage
  \fi
  \thispagestyle{empty}%
 }
 {%
  \if@restonecol
    \twocolumn
  \else
    \newpage
  \fi
 }
\makeatother

\begin{document}
\begin{titlepage}
\title{Nonparametric Estimation of Surface Integrals on Level Sets}
\author{Wanli Qiao\thanks{This work is partially supported by NSF grant DMS 1821154 and the Thomas F. and Kate Miller Jeffress Memorial Trust Award.}\\
    Department of Statistics\\
    George Mason University\\
    4400 University Drive, MS 4A7\\
    Fairfax, VA 22030\\
    Email address: wqiao@gmu.edu}
\date{\today}
\maketitle
\begin{abstract} \noindent 
Surface integrals on density level sets often appear in asymptotic results in nonparametric level set estimation (such as for confidence regions and bandwidth selection). Also surface integrals can be used to describe the shape of level sets (using such as Willmore energy and Minkowski functionals), and link geometry (curvature) and topology (Euler characteristic) of level sets through the Gauss-Bonnet theorem. We consider three estimators of surface integrals on density level sets, one as a direct plug-in estimator, and the other two based on different neighborhoods of level sets. We allow the integrands of the surface integrals to be known or unknown. For both of these scenarios, we derive the rates of the convergence and asymptotic normality of the three estimators. 
\end{abstract}
%
%
{\small 
}

\end{titlepage}

\section{Introduction}
The $c$-level set of a density function $f$ on $\mathbb{R}^d$ is defined as
\begin{align*}
\mathcal{M}_c = f^{-1}(c)=\{x\in\mathbb{R}^d:\; f(x)=c\}.
\end{align*}

The set $\mathcal{L}_c=f^{-1}[c,\infty)= \{x\in\mathbb{R}^d:\; f(x)\geq c\}$ is called the super level set. Level set estimation finds its application in many areas such as clustering (Cadre et al. 2009), classification (Mammen and Tsybakov, 1999), and anomaly detection (Steinwart et al. 2005). It has received extensive study in the literature. See, e.g., Hartigan (1987), Polonik (1995), Tsybakov (1997), Walter (1997), Cadre (2006), Rigollet and Vert (2009), Singh et al. (2009), Rinaldo and Wasserman (2010), Steinwart (2015), and Chen et al. (2017). \\

We consider $d\geq2$ in this paper. For simplicity of notation, the subscript $c$ is often omitted for $\mathcal{M}_c$ and $\mathcal{L}_c$. When $f$ has no flat part at the level $c$, $\mathcal{M}$ is a $(d-1)$-dimensional submanifold of $\mathbb{R}^d$. Given an i.i.d. sample $X_1,\cdots,X_n$ from the density $f$, we investigate the estimation of the surface integral 
\begin{align*}
\lambda_c(f, {g}) = \lambda(f, {g}) =  \int_{\mathcal{M}} g(x) d\mathscr{H}(x), 
\end{align*}
for some integrable function $g$ (which can be known or unknown), where $\mathscr{H}$ is the $(d-1)$-dimensional normalized Hausdorff measure. Our estimators are based on the kernel density estimator of $f$ given by
\begin{align}\label{kde}
\widehat f(x) = \frac{1}{nh^d}\sum_{i=1}^n K\left(\frac{x-X_i}{h}\right),
\end{align}
where $h$ is the bandwidth and $K$ is a $d$-dimensional kernel function. The plug-in estimators of $\cal{M}$ and $\cal{L}$ are given by $\wh{\mathcal{M}} = \{x\in\mathbb{R}^d:\; \widehat f(x)=c\}$ and $\wh{\mathcal{L}} = \{x\in\mathbb{R}^d:\; \widehat f(x)\geq c\}$, respectively.\\

If $g$ is known (for example, when $g\equiv1$), we consider the following three estimators for $\lambda(f,g)$: 
\begin{align}
& \lambda(\wh f,g) = \int_{\wh f^{-1}(c)} g(x)d\mathscr{H}(x), \label{surfintegral} \\ 
& \lambda_{\epsilon_n}^*(\wh f,g)  := \frac{1}{2\epsilon_n} \int_{\wh f^{-1}[c-\epsilon_n,c+\epsilon_n]}  g(x) \|\nabla \wh f(x)\| dx,\label{star}\\
%
\text{and  } &\lambda_{\epsilon_n}^\dagger(\wh f,g) := \frac{1}{2\epsilon_n} \int_{\wh f^{-1}(c)\oplus \epsilon_n} g(x) dx,\label{dagger}
\end{align}
for some $\epsilon_n>0$, where $\nabla \wh f$ is the gradient of $\wh f$, and $\wh f^{-1}(c)\oplus \epsilon_n$ is the union of all balls with radius $\epsilon_n$ and centers at $\wh f^{-1}(c)$. Here $\lambda(f,g)$ is a direct plug-in estimator, and $\lambda_{\epsilon_n}^*(\wh f,g)$ and $\lambda_{\epsilon_n}^\dagger(\wh f,g)$ can be viewed as local averages of surface integrals close to $\lambda(\wh f,g)$ when $\epsilon_n$ is small. The integration in $\lambda_{\epsilon_n}^\dagger(\wh f,g)$ is over a band or a tube of constant width around $\wh f^{-1}(c)$, while the domain of integration in $\lambda_{\epsilon_n}^*(\wh f,g)$ has varying width. See subsection~\ref{neighborhood} for more comparison among the three estimators. If $g$ is unknown and suppose it has an estimator $\wh g$, then we consider $\lambda(\wh f,\wh g)$, $\lambda_{\epsilon_n}^*(\wh f,\wh g)$ and $\lambda_{\epsilon_n}^\dagger(\wh f,\wh g)$ as three estimators for $\lambda(f,g)$. The main results of this manuscript include the rates of convergence and asymptotic normality of these estimators for the cases when $g$ is either known or unknown.\\

The surface integral $\lambda(f,g)$ is an important quantity that is involved in asymptotic theory for level set estimation. For example, it appears in Cadre (2006) as the convergence limit of the set-theoretic measure of $\mathcal{L}\Delta\wh{\mathcal{L}}:= (\cal{L} \backslash \wh{\cal L}) \cup (\wh{\cal L} \backslash {\cal L})$, which is the symmetric difference between $\cal{L}$ and $\wh{\cal L}$. Using a similar measure as the risk criterion, Qiao (2018) shows that the optimal bandwidth for nonparametric level set estimation is determined by a ratio of two surface integrals in the form of $\lambda(f,g)$. Qiao and Polonik (2019) develop large sample confidence regions for $\mathcal{M}$ and $\mathcal{L}$, for which the surface area of $\mathcal{M}$ (a special form of $\lambda(f,g)$ when $g\equiv1$) is the only unknown quantity that needs to be estimated. The quantity $\lambda(f,g)$ is also a key component in the concept of vertical density representation (Troutt et al. 2004). Surface integrals on (regression) level sets appears in optimal tuning parameter selection for nearest neighbour classifiers (Hall and Kang 2005, Samworth 2012, and Cannings et al. 2017).\\

Some important concepts in differential geometry are in the form of surface integrals. For example, the Willmore energy of $(d-1)$-dimensional smooth submanifold $\mathcal{S}$ embedded in $\mathbb{R}^d$ is defined as
\begin{align}\label{WillmoreEn}
W(\mathcal{S}) = \int_{\mathcal{S}} |H(x)|^2 d\mathscr{H}(x),
\end{align}
where $H(x)$ is the mean curvature of $\mathcal{S}$ at $x$ (Willmore, 1965). $W(\mathcal{S})$ measures the total elastic bending of $\mathcal{S}$ from a sphere if $d=3$.  The Willmore energy is widely used in studying the shape of biological cell membranes (Seifert, 1997). When $\mathcal{S}$ is a level set, the Willmore energy finds its applications in image inpainting (Caselles et al., 2008) and segmentation of spinal Vertebrae (Lim et al., 2013). \\

Minkowski functionals are surface integrals of some functions of principal curvatures on manifolds. When the manifolds are level sets, Minkowski functionals are widely used as morphological descriptors (shape statistics) in studying cosmic microwave background. See Chapter 10.3 of Mart\'{i}nez and Saar (2001), Schmalzing and G\'{o}rski (1998), Kerscher (2000), Pranav et al. (2017) and the references therein. Minkowski functionals (and their ratios) can be used to characterize different shapes, for example, planarity, filamentarity, clusters etc. See Section~\ref{unknownintegrand} for more details.\\

Surface integrals also plays a very important role in linking geometry and topology of manifolds, reflected in the Gauss-Bonnet theorem (for $d=3$):
\begin{align}\label{Gauss-Bonnet}
\int_{\mathcal{S}} \kappa(s) d\mathscr{H}(x) = 2\pi \chi(\mathcal{S}),
\end{align}
where $\kappa$ is the Gaussian curvature, and $\chi(\mathcal{S})$ is the Euler characteristic of $\mathcal{S}$, which is an alternating sum of Betti numbers of $\mathcal{S}$ (the number of holes of different dimensions). The Euler characteristic and Betti number are topological invariants that play a critical role in topological data analysis, in particular in persistent homology which captures topological changes in the evolution of level sets $\mathcal{M}_c$ as the level $c$ varies (Turner et al., 2014).\\

Surface area has been used in the definition of ``contour index'', which is the ratio between perimeter and square root of the area, and has appeared in medical imaging and remote sensing (Canzonieri and Carbone, 1998, Garcia-Dorado et al., 1992, Salas et al., 2003). The estimation of surface area also has extensive applications in stereology (Baddeley et al, 1986, Baddeley and Jensen, 2005, Gokhale, 1990). Other examples of $\lambda(f,g)$ arise where $g$ can be observable temperature, humidity, or the density of some non-homogeneous material, and one is interested in the surface integrals of these quantities on an unknown manifold (Jim\'{e}nez and Yukich, 2011). \\

In the literature Cuevas et al. (2012) obtain a consistency result for the estimation of the surface area of $\mathcal{M}$, i.e. when $g\equiv 1$. There exists some recent work on the estimation of the surface area of the boundary of an unknown body $S\subset G$ where $G$ is a bounded set given a sample on $G$. The work is relevant to the study in this manuscript but in a setting different from what we consider here. There the surface area is defined as the Minkowski content, which coincides with the normalized Hausdorff measure in regular cases (Ambrosio et al. 2008). Armend\'{a}riz et al. (2009) and Trillos et al. (2017) obtain asymptotic normality results for the surface area (or perimeter) of $\partial S$ in a framework that assumes a uniform distribution on $G$ and the binary labels for $S$ and $G\backslash S$ are observed with the sample. The sampling scheme is called the ``inside-outside'' model in Cuevas and Pateiro-L\'{o}pez (2018), and the binary labels for $S$ and $G\backslash S$ contain important information of the location of $\partial S$. In contrast, we study the estimation of surface integral $\lambda(f,g)$, which is more general than surface area. More importantly, the location of $f^{-1}(c)$ is completely unknown and needs to be estimated. As a result even for $g\equiv1$ the approach we take is very different from the above work. In the setting of the ``inside-outside'' model, the rates of convergence of estimators for surface area of $\partial S$ are derived under different shape assumptions (Cuevas et al. 2007, Pateiro-L\'{o}pez and Rodr\'{i}guez-Casal 2008, 2009). The consistency of the surface integral estimation on $\partial S$ has been considered by Jim\'{e}nez and Yukich (2011), where the integrand function is assumed to be observable on the sample points. Assuming the i.i.d. sampling scheme only on $S$, Arias-Castro and Rodr\'{i}guez-Casal (2017) estimate the perimeter of $\partial S$ using the alpha-shape for $d=2$ and derive the rate of convergence. \\

In subsection~\ref{plugin} we relate $\lambda(\wh f,g) - \lambda(f,g)$ to $d_g(f,\wh f):=\int_{{\cal{L}} \Delta \wh{\cal L}} g(x) dx$ and show their heuristic connection and difference. In the literature $d_g(f,\wh f)$ has been studied in, e.g., Ba\'{i}llo et al. (2000), Ba\'{i}llo (2003), Cadre (2006), Cuevas et al. (2006), and Mason and Polonik (2009). The estimation of surface integral we study in this manuscript is also related to estimation of integral functionals of a density function and its derivatives, where the integral is over $\mathbb{R}^d$. See Levit (1978), Hall and Marron (1987), Bickel and Ritov (1988), Gin\'{e} and Nickl (2008), and Gin\'{e} and Mason (2008).\\

We organize the manuscript as follows. In the rest of this section we introduce some notation, geometric concepts and assumptions used for our results. Sections~\ref{knownintegrand} and \ref{unknownintegrand} are dedicated to the estimation of the surface integral $\lambda(f,g)$ when $g$ is assumed to be known and unknown, respectively. The main results of these two sections are Theorems~\ref{HausdorffIntegration} and \ref{finaltheorem}, respectively, which give the rates of convergence and asymptotic normality of our estimators for the two scenarios. In Section~\ref{eulercharact} we discuss the estimation of Euler characteristic of level sets as a special case of $\lambda(f,g)$ with an unknown $g$. The proofs of all the results in Sections \ref{knownintegrand} and \ref{unknownintegrand} are left to Section~\ref{proofsection}.\\ 

\subsection{Notation and geometric concepts}\label{geometriconcept}

%
For any $x\in\mathbb{R}^d$ and $A\subset\mathbb{R}^d$, let $d(x,A)=\inf_{y\in A}\|x-y\|$. The Hausdorff distance between any two sets $A, B\subset\mathbb{R}^d$ is $$d_H(A, B) = \max\left\{\sup_{x\in B} d(x,A),\; \sup_{x\in A}d(x,B)\right\}.$$ The ball with center $x$ and radius $\epsilon$ is denoted by $\mathcal{B}(x,\epsilon) = \{y\in\mathbb{R}^d:\; \|x-y\|\leq \epsilon\}.$ For any set $A\subset\mathbb{R}^d$ and $\epsilon>0$, we denote $A\oplus\epsilon = \bigcup_{x\in A} \mathcal{B}(x,\epsilon)$. Let the normal projection of $x$ onto $A$ be $\pi_A(x) =\{y\in A:\; \|x-y\|=d(x,A)\}$. Note that $\pi_A(x)$ may not be a single point. \\ 

\begin{figure}
\centering
\includegraphics[scale=0.36]{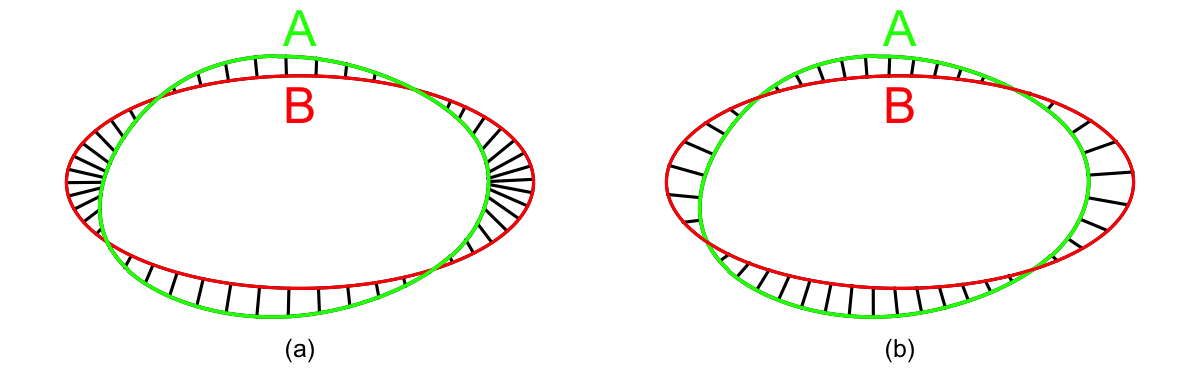}
\caption{This figure shows the normal compatibility between two curves $A$ (green) and $B$ (red). (a) shows the normal projection from $A$ to $B$. (b) shows the normal projection from $B$ to $A$.}
\label{normalcompat}
\end{figure}

We will also use the concept of {\em reach} of a manifold. For a set $\mathcal{S}\subset\mathbb{R}^d$, let Up$(\mathcal{S})$ be the set of points $x\in\mathbb{R}^d$ such that $\pi_{\mathcal{S}}(x)$ is unique. The reach of $\mathcal{S}$ is defined as
\begin{align*}
\rho(\mathcal{S}) = \sup\{\delta>0: \mathcal{S}\oplus\delta\subset \text{Up}(\mathcal{S}) \}.
\end{align*}
See Federer (1959). A positive reach is related to the concepts of ``$r$-convexity'' and ``rolling condition''. See Cuevas et al. (2012). These are common regularity conditions for the estimation of surface area. For two manifolds $A$ and $B$, if the normal projections $\pi_A: B \mapsto A$ and $\pi_B: A\mapsto B$ are homeomorphisms, then $A$ and $B$ are called {\em normal compatible} (see Chazal et al. (2007)). See Figure~\ref{normalcompat}.\\

Let $\|M\|_F$ and $\|M\|_{\max}$ be the Frobenius and max norms of a matrix $M$, respectively. For a twice differentiable function $g: \mathbb{R}^d\mapsto \mathbb{R}$, let $\nabla g$ and $\nabla^2 g$ be the gradient and Hessian matrix of $g$, respectively. For a vector or matrix $M$ and a positive integer $r$, let $M^{\otimes r}$ be the $r$-th Kronecker power of $M$. The vector of $r$-th derivatives of $g$ (if they exist) is defined as $\nabla^{\otimes r} f(x) = \frac{\partial^r f}{(\partial x)^{\otimes r}}(x)\in\mathbb{R}^{d^r}$, where we apply the $r$-th Kronecker power of the operator $\nabla$. For a bandwidth $h>0$, let $\gamma_{n,h}^{(k)} = \sqrt{\frac{\log{n}}{nh^{d+2k}}},$ which, under standard assumptions, is the uniform rate of convergence of the kernel estimator of the $k$-th density derivatives after being centered at their expectation. We denote $a\vee b=\max(a,b)$ for $a,b\in\mathbb{R}.$\\ 


\subsection{Assumptions and their discussion}
The kernel function $K$ is called a $\nu$th ($\nu\geq2$) order kernel if $\int_{\mathbb{R}^d} u^{\otimes\nu} K(u) du\neq0$ and finite while 
\begin{align*}
    \int_{\mathbb{R}^d} u^{\otimes l} K(u) du=
    \begin{cases}
      1, & \text{if } l=0, \\
      0, & \text{if } l=1,\cdots,\nu-1,\\
    \end{cases}
%
%
\end{align*}

where we use the convention $u^{\otimes0}=1$. For Kronecker power used in matrices and multivariable Taylor expansion as well as high-order kernels, we refer the reader to Chac\'{o}n and Duong (2018). \\

We introduce some assumptions that will be used in this paper. 
For $\delta>0$, denote $\mathcal{I}(\delta)= \mathcal{I}_f(\delta)= f^{-1}([c-\delta,c+\delta]) =  \{x\in\mathbb{R}^d:\; f(x)\in [c-\delta,c+\delta]\}$, which is a neighborhood of $\mathcal{M}$.\\

{\bf Assumptions}: \\
(K) $K$ is a two times continuously differentiable kernel function of $\nu$th order for $\nu\geq 2$, with bounded support. \\

(F1) The density function $f$ has continuous partial derivatives up to order $\nu$ on $\mathcal{I}(2\delta_0)$ for some $\delta_0>0$. There exists $\epsilon_0>0$ such that  $\|\nabla f(x)\| > \epsilon_0$ for all $x\in \mathcal{I}(2\delta_0)$.\\ 

(H) The bandwidth $h$ depends on $n$ such that $\gamma_{n,h}^{(2)}\rightarrow 0$ and $h\rightarrow0$ as $n\rightarrow\infty$.\\ 

{\bf Discussion of the assumptions}:\\
1. The assumption $\|\nabla f(x)\| > \epsilon_0$ for $x\in\mathcal{I}(2\delta_0)$ in (F1) implies that the Lebesgue measure of $\mathcal{M}$ on $\mathbb{R}^d$ is zero. The level $c$ is called a regular value of $f$ and it is implied that $\inf_{x\in\mathbb{R}} f(x) < c<\sup_{x\in\mathbb{R}} f(x)$. This is a typical assumption in the literature of level set estimation (see, e.g. Cadre (2006), Cuevas et al. (2006), Mason and Polonik (2009), Mammen and Polonik (2013)), and guarantees that $\mathcal{M}$ has no flat parts and is a compact $(d-1)$-dimensional manifold (see Theorem 2 in Walther, 1997). In particular, under assumption (F1) the following well-known margin condition (first introduced by Polonik, 1995) is satisfied: $P(|f(X)-c|<\epsilon)\leq C\epsilon$ for some positive constant $C$ and small $\epsilon$ (see Lemma 4 in Rinaldo and Wasserman, 2010). The requirement for the $\nu$th partial derivatives of $f$ is mainly to deal with the bias in the kernel density estimation. If our estimation target is $ \lambda(\mathbb{E}\wh f,g)$ instead of $\lambda(f,g)$ when $g$ is known, then we only need continuous second partial derivatives of $f$ on $\mathcal{I}(2\delta_0)$.\\ 

2. The estimator $\widehat f$ using a high-order kernel ($\nu>2$) can take negative values, and therefore loses its interpretability for practitioners (see, e.g., Silverman, 1986, page 69, Hall and Murison, 1993).  However, this is not a problem for our estimator, since we are only interested in the level set of $f$ at a positive level, and our estimators $\wh{\cal M}$ and $\wh{\cal L}$ do not directly use negative values of $\widehat f$. \\

3. Assumption (H) guarantees the strong uniform consistency of the second partial derivatives of $f$. This is used to establish the normal compatibility between $\mathcal{M}$ and $\wh{\mathcal{M}}$, so that we can explicitly define a homeomorphism between them. A similar assumption has been used in Chen et al. (2017) and Qiao and Polonik (2019) for the construction of confidence regions of level sets.

\section{Surface integral estimation with a known integrand}\label{knownintegrand}

We first consider the estimation of $\lambda(f,g)$ when $g$ is a known function. Three estimators $\lambda(\wh f,g)$, $\lambda_{\epsilon_n}^*(\wh f,g)$ and $\lambda_{\epsilon_n}^\dagger(\wh f,g)$ defined in (\ref{surfintegral}) - (\ref{dagger}) are investigated. In general we use $\wt\lambda_{\epsilon_n}(\wh f,g)$ to denote any one of them. In particular, when we use $\wt\lambda_{\epsilon_n}$ to denote $\lambda$, we take $\epsilon_n\equiv0$.\\ 

We give an asymptotic normality result for $\wt\lambda_{\epsilon_n}(\wh f,g) -\lambda(f,g)$ in the following theorem. To formulate the result we need to introduce more notation and geometric concepts (especially the curvatures on level sets). We denote tangent space of $\mathcal{M}$ at $x\in\mathcal{M}$ by $T_x(\mathcal{M})$. Let $\n(x) = \frac{\nabla f(x)}{\|\nabla f(x)\|}$ be the normalized gradient. For any $t\in\mathbb{R}$, let $T_{x,t}(\mathcal{M})=\{y+t\n(x),y\in T_x(\mathcal{M})\}$, which is an affine tangent plane. If $-\n$ is chosen as the surface normal of $\mathcal{M}$, then $\nabla \n(x)$ is called the shape operator (or Weingarten map). Let $G(x) = \mathbf{I}_d - \n(x)\n(x)^T$. Direct calculation shows that (see Kindlmann et al. 2003) 
 \begin{align}\label{Kindlmann}
 \nabla \n(x) = \|\nabla f(x)\|^{-1} G(x)\nabla^2 f(x).
 \end{align}
The following matrix is equivalent to $\nabla \n(x)$ as a linear map on $T_x(\mathcal{M})$ (see (\ref{shapeoppro})) and hence is also called the shape operator of $\mathcal{M}$ at $x$:
 \begin{align}\label{shapeoperator}
 S(x) = \nabla \n(x) G(x) = \|\nabla f(x)\|^{-1} G(x)\nabla^2 f(x) G(x).
 \end{align}
The shape operator $S(x)$ has $d-1$ non-zero eigenvalues, which are principal curvatures of $\mathcal{M}$ at $x$. The mean curvature, denoted by $H(x)$ is the average of the $(d-1)$ principal curvatures, i.e., $H(x) = \frac{1}{d-1}{\bf tr} [S(x)] =  \frac{1}{d-1}\|\nabla f(x)\|^{-1}{\bf tr} [G(x)\nabla^2 f(x)]$. 
%
Consider the plug-in estimators of the above quantities. Let $\wh \n(x)=\frac{\nabla\wh f(x)}{\|\nabla \wh f(x)\|}$, $\wh G(x) =\mathbf{I}_d - \wh \n(x)\wh \n(x)^T$, $\wh S(x) = \nabla \wh\n(x)\wh G(x)$, $\wh H(x)=\frac{1}{d-1}\mathbf{tr}[\wh S(x)]$. For a kernel function $K$ and $x\in\mathcal{M}$, let $R_K(x)=\int_{\mathbb{R}} \left(\int_{T_{x,t}(\mathcal{M})} K\left(u\right)d\mathscr{H}(u)\right)^2dt$. If $K$ is spherically symmetric, then $R_K(x)$ is a constant for all $x\in\mathcal{M}$ and we write $R(K)=R_K(x)$. 

\begin{theorem}\label{HausdorffIntegration}
Let $g:\mathcal{I}(2\delta_0) \mapsto \mathbb{R}$ be a function with bounded continuous second partial derivatives. Let $\wt\lambda_{\epsilon_n}$ be either $\lambda$, $\lambda_{\epsilon_n}^*$ or $\lambda_{\epsilon_n}^\dagger$. Let $\alpha_{n,h}:=\frac{1}{nh^{d+2}}+ h^\nu + \frac{1}{\sqrt{nh}} + \epsilon_n^2=o(1)$. Under assumptions (K), (F1) and (H), we have the following.\\
(i) 
\begin{align}\label{assertion1}
\wt\lambda_{\epsilon_n}(\wh f,g) -\lambda(f,g) = O_p\left(\alpha_{n,h} \right).
\end{align}
(ii) If in addition $nh^{2d+3}\rightarrow \infty$, $nh^{1+2\nu}\rightarrow\gamma\geq0$ and $\epsilon_n^2\sqrt{nh}\rightarrow0$, we have
\begin{align}\label{assertion2}
\sqrt{nh} [\wt\lambda_{\epsilon_n}(\wh f,g) -\lambda(f,g)] \rightarrow_D \mathscr{N}(\sqrt{\gamma}\mu,\sigma^2),
\end{align}
where $\mu$ is given in (\ref{asymmean}) and $\sigma^2=c\lambda(f,w_g^2R_K)$ with $w_g(x) = \|\nabla f(x)\|^{-1}[ \n(x)^T \nabla g(x) + (d-1)H(x)g(x)].$\\
(iii) Suppose we choose to use a spherically symmetric kernel function $K$, and then $\sigma^2 = c R(K)\lambda(f,w_g^2)$ in (ii). For a nonnegative sequence $\tau_n=O(h)$, let $\wh \sigma_{\tau_n}^2 = cR(K)\wt\lambda_{\tau_n}(\wh f,\wh w_g^2)$, where $\wh w_g(x) = \|\nabla \wh f(x)\|^{-1}[ \wh \n(x)^T \nabla g(x) +  (d-1)\wh H(x)g(x)]$. Under all the assumptions in (ii) with $\gamma=0$ and suppose $\sigma^2\neq0$ and $K$ is four times continuously differentiable, we have 
\begin{align}\label{assertionnorm}
\sqrt{nh} \; \wh \sigma_{\tau_n} ^{-1}[\wt\lambda_{\epsilon_n}(\wh f,g) -\lambda(f,g)] \rightarrow_D \mathscr{N}(0,1).
\end{align}
\end{theorem}

\begin{remark}\label{discuss} $\;$

{\em
%
a) The above result can be used to construct a confidence interval for $\lambda(f,g)$. For $0<\alpha<1$, Let $z_{\alpha/2}$ be the $(1-\alpha/2)$ quantile of $\mathscr{N}(0,1)$. With rates of $h$, $\tau_n$ and $\epsilon_n$ chosen complying with the above theorem, a $(1-\alpha)$ asymptotic confidence interval for $\lambda(f,g)$ is given by $$\left[\wt\lambda_{\epsilon_n}(\wh f,g) - \frac{1}{\sqrt{nh}}\wh \sigma_{\tau_n}  z_{\alpha/2},\; \wt\lambda_{\epsilon_n}(\wh f,g) + \frac{1}{\sqrt{nh}}\wh \sigma_{\tau_n} z_{\alpha/2} \right].$$
Note that the choice of $\wt\lambda$ used in $\wh\sigma_{\tau_n}$ does not have to the same one as in $\wt\lambda_{\epsilon_n}(\wh f,g)$.

%
b) Below we give interpretation of the rates in $\alpha_{n,h}$. Here $O_p(\epsilon_n^2)$ accounts for the difference between $\wt\lambda_{\epsilon_n}(\wh f, g)$ and $\lambda(\wh f, g)$, i.e., the error in approximating the surface integral by an integral over a neighborhood of width $\epsilon_n$ around the level set. We then focus on $\lambda(\wh f, g)-\lambda(f, g)$ and it turns out that we have the approximation $\lambda(\wh f, g) - \lambda(f,g)\approx  T_{n,1} + T_{n,2},$ where $T_{n,1}=\lambda(f, w_g\times(f-\wh f))$ and $T_{n,2}=\lambda(f, p(\nabla\wh f-\nabla f))$ with $p(\nabla\wh f-\nabla f)$ a quadratic form of $\nabla\wh f- \nabla f$. It is known that $\wh f(x)-\mathbb{E}\wh f(x)$ has a standard rate of $1/\sqrt{nh^{d}}$. Due to the integrals on a $(d-1)$-dimensional manifold $\mathcal{M}$, we gain $(d-1)$ powers of $h$, which results in the rate $1/\sqrt{nh}$ for the stochastic part of $T_{n,1}$. The fact that $(nh^{d+2})^{-1/2}$ is the rate of convergence of $\nabla\wh f(x)-\mathbb{E}\nabla\wh f(x)$ explains $(nh^{d+2})^{-1}$ as a rate for the bias part of $T_{n,2}$. The rate $h^\nu$ accounts for the remaining bias in the overall estimation. To have the asymptotic normality in (\ref{assertion2}), we use the conditions $nh^{2d+3}\rightarrow \infty$ and $nh^{1+2\nu}\rightarrow 0$ to make the bias asymptotically negligible. Under assumption (H), these conditions require $\nu>d+1$. 
%

c) The asymptotic normality in (\ref{assertionnorm}) holds because $\wh \sigma_{\tau_n}^2$ is a consistent estimator of $\sigma^2$. Estimating $\sigma^2$ by $\wh \sigma_{\tau_n}^2$ is a special case of estimating surface integrals when the integrand is unknown, which is further investigated in Section~\ref{unknownintegrand}.
%
%

d) We discuss some special cases of $g$. If $g\equiv 1$, then $w_g(x) = (d-1)\|\nabla f(x)\|^{-1}H(x)$ and the surface integral $\lambda(f,w_g^2)$ in $\sigma^2$ is a weighted Willmore energy, whose original definition is given in (\ref{WillmoreEn}). In some extreme cases the choice of $g$ can make $w_g\equiv0$, and then the asymptotic normal distribution in part (ii) is degenerate. We have excluded such a degenerate scenario in part (iii).

}
\end{remark}

In order to prove the rates of convergence and asymptotic normality results in Theorem~\ref{HausdorffIntegration}, we need derive a few results that are also interesting in their own right. We first consider the direct plug-in estimator $\lambda(\wh f,g)$, and then $\lambda_{\epsilon_n}^*(\wh f,g)$ and $\lambda_{\epsilon_n}^\dagger(\wh f,g)$, which are based on neighborhoods of $\wh{\mathcal{M}}$.

\subsection{Direct plug-in estimation}\label{plugin}

The proof of Theorem~\ref{HausdorffIntegration} when $\wt\lambda_{\epsilon_n}=\lambda$ is built upon the results below in this section. One of the main challenges is that the domains of integrals in $\lambda(f,g)$ and $\lambda(\wh f,g)$ are not the same, which makes the comparison difficult. Briefly speaking, our strategy is to establish a diffeomorphism between the two domains and utilize the area formula in differential geometry to convert the surface integrals into those with the same integral domains. We provide some heuristic for this procedure, which also explains the role of the curvature of level sets in $\sigma^2$. \\ 

\begin{figure}
\centering
\includegraphics[scale=0.4]{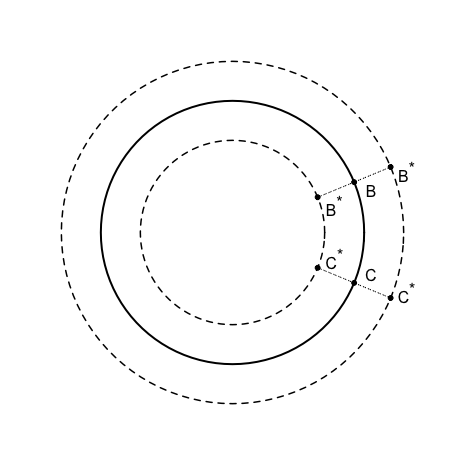}
\caption{This figure provides some heuristic to understand the connection and difference between the estimation of Hausdorff measure of $\mathcal{M}$ and the Lebesgue measure of $\wh{\mathcal{L}}\Delta \mathcal{L}$.}
\label{heuristic}
\end{figure}

In Figure~\ref{heuristic}, we consider a circle with radius $r$ as our level set $\mathcal{M}$. An ideal estimator $\wh{\mathcal{M}}$ is a circle with the same center and radius $r+\Delta r$. Here we allow $\Delta r$ to be either positive or negative, and the two dotted circles on the graph are two possible versions of $\wh{\mathcal{M}}$. In this heuristic we only consider the perimeter of the circles, corresponding to $g\equiv1$ in the surface integral on level sets. From elementary geometry it is known that the difference between the perimeters of $\wh{\mathcal{M}}$ and $\mathcal{M}$ is $2\pi\Delta r$.  We focus on the local geometry to better understand the behavior of this difference. We consider short arcs $\wideparen{BC}$ on $\mathcal{M}$ and $\wideparen{B^*C^*}$ on $\wh{\mathcal{M}}$, where both $BB^*$ and $CC^*$ can be extended to go through the center of the circles. The difference between the lengths of the two arcs is $|\wideparen{B^*C^*}| - |\wideparen{BC}| = \Delta r r^{-1} |\wideparen{BC}|,$ where $r^{-1}$ can be understood as the curvature of $\mathcal{M}$, and we denote it by $H$ (the same notation for mean curvature). Now imagine that both $\mathcal{M}$ and $\wh{\mathcal{M}}$ are slightly deformed from circles, and we allow both $\Delta r$ and $H$ to depend on location $x$, then the difference of arc lengths can be approximated by a surface integral $\int_{\wideparen{BC}} \Delta r(x) H (x) d\mathscr{H}(x)$. It is straightforward to extend this approximation to the entire surface, i.e., $\mathscr{H}({\wh{\mathcal{M}}})  - \mathscr{H}({\mathcal{M}}) \approx \int_{\mathcal{M}} \Delta r(x) H (x) d\mathscr{H}(x).$ Another key aspect is $\Delta r(x)\approx \frac{f(x) - \wh f(x)}{\|\nabla f(x)\|}, $ where $\|\nabla f(x)\|$ reflects a ratio between vertical and horizontal variations in level set estimation. See Lemma~\ref{pointlevelset} below and its remark for details. Note that this approximation has included the sign of $\Delta r$ into consideration. So overall 
%
\begin{align}\label{approximation1}
\mathscr{H}({\wh{\mathcal{M}}})  - \mathscr{H}({\mathcal{M}}) \approx \int_{\mathcal{M}} \frac{f(x) - \wh f(x)}{\|\nabla f(x)\|} H (x) d\mathscr{H}(x).
\end{align}

The above approximation explains the role of curvature in the estimation of surface integrals on level sets. In Remark~\ref{discuss}b we have given the approximation $\lambda(\wh f, g) - \lambda(f,g)\approx T_{n,1}+T_{n,2}$, which coincides with (\ref{approximation1}) when $g\equiv 1$ if $T_{n,2}$ is ignored.\\ 

Using Figure~\ref{heuristic} again, we heuristically explain the difference and connection between $\mathscr{H}({\wh{\mathcal{M}}})  - \mathscr{H}({\mathcal{M}})$ and the volume of the symmetric difference $\wh{\mathcal{L}}\Delta \mathcal{L}$, which corresponds to the band between $\mathcal{M}$ and $\wh{\mathcal{M}}$. From elementary geometry, it is known that the volume of this band is approximately $\pi r|\Delta r|$, if $\Delta r$ is small.  Using the arc notation, we have that $\text{Volume} (\wideparen{BCC^*B^*}) \approx |\Delta r| |\wideparen{BC}|$. 
So overall
\begin{align}\label{approximation2}
\text{Volume} (\wh{\mathcal{L}}\Delta \mathcal{L}) \approx  \int_{\mathcal{M}} \frac{|f(x) - \wh f(x)|}{\|\nabla f(x)\|}  d\mathscr{H}(x),
\end{align} 
which is an $L_1$ type of integral (see Qiao, 2018). The approximations in (\ref{approximation1}) and (\ref{approximation2}) provides some insight into the connection and difference between estimating the surface integral on $\mathcal{M}$ and the set-theoretic measure of $\wh{\mathcal{L}}\Delta \mathcal{L}$: both can be approximated by a surface integral on $\mathcal{M}$ and the integrands are related to $f-\wh f$, but the latter is integrating an absolute value and only the former is (asymptotically) impacted by the curvature of level sets.\\

Below we formulate the results that have been outlined in the heuristic above. Essentially we need to establish a homeomorphism between $\mathcal{M}$ and $\wh{\mathcal{M}}$, which is given in Lemma~\ref{compatible} below. We show the approximation $\Delta r(x)\approx \frac{f(x) - \wh f(x)}{\|\nabla f(x)\|} $ in Lemma~\ref{pointlevelset}. The asymptotic normality of $\lambda(f,w_g\times(\wh f - \mathbb{E}\wh f))$ is given in Theorem~\ref{asympnorm}. \\

The following lemma is regarding the reach of the true and estimated density level sets as well as their normal compatibility (see subsection~\ref{geometriconcept} for these geometric concepts). It is similar to Lemma 1 in Chen et al. (2017), but our result is under slightly different assumptions and holds uniformly for a collection of density level sets. Also see Theorems 1 and 2 in Walther (1997) for relevant results. 

\begin{lemma}\label{compatible}
Under assumptions (F1), (K), and (H), we have the following results.\\
(i) There exists $r_0>0$, such that the reach $\rho(\mathcal{M}_\tau)>r_0$ for all $\tau\in[c-\delta_0,c+\delta_0]$. \\
(ii)  With probability one, there exists $r_0>0$, such that the reach $\rho(\wh{\mathcal{M}}_\tau)>r_0$ for all $\tau\in[c-\delta_0,c+\delta_0]$ and $n$ large enough.\\
(iii) With probability one, $\wh{\mathcal{M}}_\tau$ and $\mathcal{M}_\tau$ are normal compatible for all $\tau\in[c-\delta_0,c+\delta_0]$ and $n$ large enough.

\end{lemma}
%
%

With this result we can explicitly define homeomorphisms between level sets and their plug-in estimators. For $x\in\mathcal{M}_\tau$ with $|\tau-c|\leq\delta_0$, and $t\in\mathbb{R}$, let
\begin{align}\label{projection}
\zeta_x(t) = x +  t\n(x).
\end{align}
Note that $\{\zeta_x(t):\; t\in\mathbb{R}\}$ is orthogonal to the tangent space of $\mathcal{M}$ at $x$. Furthermore, let $t_n(x) = \argmin_{t}\left\{|t|:\; \zeta_x( t ) \in \wh{\mathcal{M}}_\tau \right\},$ which is the first time point when $\zeta_x(t)$ hits $\wh{\mathcal{M}}_\tau$ starting from a point $x$ on $\mathcal{M}_\tau$, and 
\begin{align}\label{map}
P_n(x) = \zeta_x( t_n(x)).
\end{align}

Under the assumptions in Lemma \ref{compatible}, $\wh{\mathcal{M}}_\tau$ and $\mathcal{M}_\tau$ are normal compatible for large $n$. Hence $t_n(x)$ is uniquely defined, and $P_n$ is a homeomorphism between $\mathcal{M}_\tau$ and $\wh{\mathcal{M}}_\tau$. 
%
%
%
The following lemma gives asymptotic approximations for $t_n(x)$ and its gradient. 

\begin{lemma}\label{pointlevelset}
Suppose that assumptions (F1), (K), and (H) hold. For any point $x\in \mathcal{I}(\delta_0)$, we have
\begin{align}
t_n(x) &= \frac{f(x) - \widehat f(x)}{\|\nabla f(x)\|} + \delta_{n,1}(x),\label{delta1nexpr}\\
\nabla t_n(x) &= \frac{\nabla f(x) - \nabla \widehat f(x)}{\|\nabla f(x)\|} - \frac{  \nabla \|\nabla f(x)\|}{\|\nabla f(x)\|^2 } t_n(x) + \delta_{n,2}(x),\label{delta23nexpr}
\end{align}
where 
\begin{align}
\sup_{x\in \mathcal{I}(\delta_0)} |\delta_{n,1}(x)| & = O_p\left(\gamma_{n,h}^{(0)}\;\gamma_{n,h}^{(1)}\right) +o_p\left(h^{\nu-1}(\gamma_{n,h}^{(0)} + h^\nu)\right),\label{delta1nrate}\\ 
\sup_{x\in \mathcal{I}(\delta_0)} \|\delta_{n,2}(x)\| & = O_p\left((\gamma_{n,h}^{(1)})^2\right) +o_p\left(h^{\nu-2}(\gamma_{n,h}^{(0)} + h^\nu)\right).\label{delta2nrate} 
%
%
\end{align}
\end{lemma}

\begin{remark}\label{verticalhorz} $\;${\em
%
Note that $|t_n(x)| = \|P_n(x)-x\|$. Then the above results link the local horizontal variation $t_n(x)$ with the local vertical variation $f(x) - \widehat f(x)$ as well as their derivatives. Here $\|\nabla f(x)\|$ on the right-hand side of (\ref{delta1nexpr}) can be understood as a directional derivative of $f$ in the gradient direction, i.e., $\|\nabla f(x)\| = \langle \nabla f(x), \n(x)\rangle$, and it reflects the asymptotic rate of change between the local vertical and horizontal variations, as $P_n(x)-x$ is parallel to the direction of $\n(x)$.} 
\end{remark}

As heuristically explained in (\ref{approximation1}), $\lambda(\wh f,g)-\lambda(f,g)$ can be approximated by $\lambda(f,w_g\times (f-\wh f))$. The following theorem gives the asymptotic normality of $\sqrt{nh}\lambda(f,w\times (\wh f-\mathbb{E}\wh f))$ for a weight function $w$.

\begin{theorem}\label{asympnorm}
Suppose that assumptions (F1) and (K) hold. 
%
let $w: \mathcal{M}\mapsto\mathbb{R}$ be a bounded continuous function. If $nh\rightarrow\infty$ and $h\rightarrow0$, then
\begin{align}\label{asymptocnormdens}
\sqrt{nh} \int_{\mathcal{M}} w(x)[\widehat f(x) -\mathbb{E} \widehat f(x) ] dx \rightarrow_D \mathscr{N}(0,\sigma^2),
\end{align}
where $\sigma^2=c\lambda(f,w^2R_K)$.
\end{theorem}

\subsection{Estimation using integrals over neighborhoods of level sets}\label{neighborhood}

We consider $ \lambda_{\epsilon_n}^*(\wh f,g) $ and $ \lambda_{\epsilon_n}^\dagger(\wh f,g)$ as two alternative estimators of $\lambda(f,g)$, where $ \lambda_{\epsilon_n}^*$ and $\lambda_{\epsilon_n}^\dagger$ are given in (\ref{star}) and (\ref{dagger}). It can be seen that (by Proposition A.1 of Cadre (2006) and Section 3.4 of Gray, 2004)
\begin{align}
& \lambda_{\epsilon_n}^*(\wh f,g) =\frac{1}{2\epsilon_n} \int_{-\epsilon_n}^{\epsilon_n}  \lambda_{c+\epsilon}(\wh f,g)d\epsilon, \label{altstar}\\
& \lambda_{\epsilon_n}^\dagger(\wh f,g)  = \frac{1}{2\epsilon_n} \int_{0}^{\epsilon_n} \int_{\partial (\wh f^{-1}(c) \oplus \epsilon)} g(x) d\mathscr{H}(x)d\epsilon\label{altdagger},
\end{align}
where $\partial (\wh f^{-1}(c) \oplus \epsilon) = \{x: d(x,\wh f^{-1}(c)) = \epsilon\}$. So $ \lambda_{\epsilon_n}^*(\wh f,g) $ and $ \lambda_{\epsilon_n}^\dagger(\wh f,g) $ can be viewed as two local averages of surface integrals in a neighborhood of $\wh f^{-1}(c)$. They may have some computational advantage over $\lambda(\wh f,g)$, because numerical approximation for integrals over subsets of $\mathbb{R}^d$ in (\ref{star}) and (\ref{dagger}) is less challenging than numerical approximation to a surface integral in (\ref{surfintegral}) for $d\geq 3$.\\

Even though both of $ \lambda_{\epsilon_n}^*(\wh f,g) $ and $ \lambda_{\epsilon_n}^\dagger(\wh f,g) $ are based on integration over some neighborhoods of $\wh f^{-1}(c)$, their integration domains have different shapes, which can potentially impact their performance. The region $\wh f^{-1}[c-\epsilon_n,c+\epsilon_n]$ is an ``implicit'' tube defined through variation of the vertical levels of $\wh f$, while the tube $\wh f^{-1}(c)\oplus\epsilon_n$ has a constant radius, grown through horizontal variation from $\wh f^{-1}(c)$. Note that the roles of $\epsilon_n$ in $\lambda_{\epsilon_n}^*$ and $\lambda_{\epsilon_n}^\dagger$ are the magnitude of vertical and horizontal variations, respectively, and the relevant geometric interpretation for why $\|\nabla \wh f\|$ appears in the integrand in (\ref{star}) can be found in Remark~\ref{verticalhorz}. When $\epsilon_n$ are appropriately chosen, the two different types of tubes have been used as asymptotic confidence regions for $\wh f^{-1}(c)$ in Chen et al. (2017) and Qiao and Polonik (2019). Finite sample performance of these two types of confidence regions are compared in Qiao and Polonik (2019) and it has been observed that the two types confidence regions behave differently, especially when $c$ is close to 0 or or a critical value of $f$. We expect these scenarios also impact the integration in $\lambda_{\epsilon_n}^*(\wh f,g)$ and $\lambda_{\epsilon_n}^\dagger(\wh f,g)$. \\

Let $\wt\lambda_{\epsilon_n}$ be one of $\lambda_{\epsilon_n}^*$ and $\lambda_{\epsilon_n}^\dagger$. Then it can be seen that $\wt\lambda_{\epsilon_n}(\wh f, g) - \lambda(f,g) = [\wt\lambda_{\epsilon_n}(\wh f, g)-\lambda(\wh f, g)] + [\lambda(\wh f, g)-\lambda(f, g)].$ The results developed in subsection~\ref{plugin} are useful to study $\lambda(\wh f, g)-\lambda(f, g)$. The difference $\wt\lambda_{\epsilon_n}(\wh f, g)-\lambda(\wh f, g)$ accounts for the extra error in the estimation caused by replacing a surface integral by an integration over a tubular neighborhood of $\wh f^{-1}(c)$ with (vertical or horizontal) radius $\epsilon_n$. The following theorem states that this difference is of the order of $\epsilon_n^2$. 

\begin{theorem}\label{lambdastardaggertheorem}
Let $g:\mathcal{I}(2\delta_0) \mapsto \mathbb{R}$ be a function with bounded continuous second partial derivatives. Let $\wt\lambda_{\epsilon_n}$ be either $\lambda_{\epsilon_n}^*$ or $\lambda_{\epsilon_n}^\dagger$. Under assumptions (F1), (K), and (H), as $\epsilon_n\rightarrow0$ we have 
\begin{align}\label{assertion1star}
|\wt\lambda_{\epsilon_n}(\wh f, g)-\lambda(\wh f, g)|=O_p(\epsilon_n^2)
\end{align}
\end{theorem}
%

In Theorems~\ref{HausdorffIntegration} and \ref{lambdastardaggertheorem}, the integrand $g$ is assumed to have bounded continuous second partial derivatives. We can extend the above result when we replace $g$ in $\wt\lambda_{\epsilon_n}(\wh f, g) - \lambda(f, g)$ by a sequence of random functions $b_n$. This is useful when the case of surface integral estimation with unknown integrand $g$ is studied, where we can take $b_n = \wh g - g$, if $\wh g$ is an estimator of $g$. The following corollary give the rate of convergence of $\wt\lambda_{\epsilon_n}(\wh f, b_n) - \lambda(f, b_n)$. To get this result, essentially we need to closely examine how $\wt\lambda_{\epsilon_n}(\wh f, g) - \lambda(f, g)$ depends on $g$ and its first and second derivatives in the proof of Theorem~\ref{lambdastardaggertheorem}.

\begin{corollary}\label{Hausdorff}
Suppose that assumptions (F1), (K), and (H) hold. For a sequence of functions $b_n:\mathcal{I}(2\delta_0) \mapsto \mathbb{R}$, suppose that there exist sequences $\alpha_n$, $\beta_n$ and $\eta_n$ such that $\sup_{x\in \mathcal{I}(2\delta_0)} |b_n(x)| = O_p(\alpha_n)$ and $\sup_{x\in \mathcal{I}(2\delta_0)} \|\nabla b_n(x)\| = O_p(\beta_n)$ and $\sup_{x\in \mathcal{I}(2\delta_0)} \|\nabla^2 b_n(x)\|_F = O_p(\eta_n)$. 
%
Let $\wt\lambda_{\epsilon_n}$ be $\lambda$ or $\lambda_{\epsilon_n}^*$ or $\lambda_{\epsilon_n}^\dagger$ with $\epsilon_n\rightarrow0$. Then 
%
\begin{align}\label{Dgn0rate2}
\wt\lambda_{\epsilon_n}(\wh f, b_n) - \lambda(f, b_n) =  O_p\left\{\left(\alpha_n+\beta_n\right) \left(\gamma_{n,h}^{(0)}+ h^\nu\right) + (\alpha_n+\beta_n+\eta_n) \epsilon_n^2 \right\}.
\end{align}
\end{corollary}
\section{Surface integral estimation with an unknown integrand}\label{unknownintegrand}

In this section we consider the estimation of $\lambda(f,g)$ when $g$ is unknown and has to be estimated. Let $\wh g$ be an estimator of $g$. The estimators for $\lambda(f,g)$ we consider are $\wt\lambda_{\epsilon_n}(\wh f, \wh g)$, where $\wt\lambda_{\epsilon_n}$ is one of $\lambda$, $\lambda_{\epsilon_n}^*$ and $\lambda_{\epsilon_n}^\dagger$. Before we show the convergence rates and asymptotic normality for these estimators, recall some application examples of surface integral estimation with an unknown integrand: (i) the Willmore energy $W(\mathcal{M})$ defined in (\ref{WillmoreEn}), 
(ii) the Euler characteristic of level sets through the Gauss Bonnet theorem in (\ref{Gauss-Bonnet}), (iii) the variance in the asymptotic normality result in Theorem~\ref{HausdorffIntegration}, (iv) Minkowski functionals of level sets, and (v) plug-in bandwidth selection for nonparametric estimation of density level sets. The surface integrals in examples (i)-(iii) have been defined. Below we will give some description of surface integrals in examples (iv) and (v). \\

{\bf Minkowski functionals. } The Minkowski functionals in the cosmology literature are mainly used as integrals on the sets of two different models: (1) union of balls with varying radii $r$ centered at a point set, as a realization of a Poisson process (see Mecke et al, 1994), and (2) super level set ${\cal L}_c$ with varying levels $c$ of a density function $f$ (see Schmalzing and G\'{o}rski, 1998). We give the definition of Minkowski functionals in the framework of (2). Let $\omega_j = \frac{2\pi^{(j+1)/2}}{\Gamma[(j+1)/2]}$, where $\Gamma$ is the Gamma function. The $d+1$ Minkowski functionals are
\begin{align*}
V_j(\cal L) =  \begin{cases}
\mathscr{H}(\cal{L}), & j=0\\
\left[\omega_j {d \choose j} \right]^{-1} \lambda(f,F_j(\kappa_1,\cdots,\kappa_{d-1})), & j=1,\cdots,d,
\end{cases}
\end{align*}
where $F_j$ is defined via the polynomial expansion 
\begin{align}\label{minkfunction}
\prod_{i=1}^{d-1} (x+\kappa_i) = \sum_{j=1}^d x^{d-j} F_j(\kappa_1,\cdots,\kappa_{d-1}).
\end{align}

Note that $F_1=1$, $F_2=\kappa_1 + \cdots + \kappa_{d-1}$, and $F_d=\kappa_1\cdots\kappa_{d-1}$. Here $\kappa_1,\cdots,\kappa_{d-1}$ together with zero are the roots of the characteristic polynomial $\text{det}(t\mathbf{I}_d - S)$, where $S$ is the shape operator given in (\ref{shapeoperator}). In other words, for $j=1,\cdots,d$, $(-1)^{j-1}F_j(\kappa_1(x),\cdots,\kappa_{d-1}(x))$ is the coefficient before $t^{d-j+1}$ in $\text{det}(t\mathbf{I}_d - S)$, and so $F_j$'s only depend on the first and second derivatives of the density. Some of Minkowski functionals has very intuitive geometric or topological interpretation:  $V_0$ is the $d$-dimensional volume of $\cal L$, $V_1$ is the $(d-1)$-dimensional volume of the boundary $\cal M$, $V_2$ is the total mean curvature of $\cal M$, and $V_d$ is the Euler characteristic of $\cal M$ (by the generalized Gauss-Bonnet theorem) up to some constants.\\

{\bf Estimation of surface integral in bandwidth selection for level set estimation.}
Bandwidth selection is important for kernel-type estimators, such as $\wh{\mathcal{M}}$ and $\wh{\mathcal{L}}$. The risk $\mathbb{E} \int_{{\cal{L}} \Delta \wh{\cal L}} g(x) dx$ for some integrable function $g$ is widely used as a measure of dissimilarity between $\cal{L}$ and $\wh{\cal L}$ (or between $\cal{M}$ and $\wh{\cal M}$). 
In particular, Qiao (2018) shows that under regularity conditions for an integer $p\geq 1$, $m_p(h):=\mathbb{E} \int_{{\cal{L}} \Delta \wh{\cal L}} |f(x)-c|^{p-1} dx = \frac{1}{p}m_p(h)\{1+o(1)\},$ where
$$ \wt m_p(h) =\int_{\cal{M}} \frac{\mathbb{E}  |\widehat f(x) - f(x)|^p}{\|\nabla f(x)\|}d\mathscr{H}(x). $$
Choosing $p=2$ makes the above risk an $L_2$ type of risk. Notice that 
\begin{align*}
\wt m_2(h) = \int_{\mathcal{M}} \frac{|\mathbb{E} \wh f(x) -f(x)|^2 + \text{Var} (\wh f(x)) }{\|\nabla f(x)\| } d\mathscr{H}(x), 
%
\end{align*}
where if we use a second order kernel ($\nu=2$),
\begin{align*}
& \mathbb{E}\wh f(x) -f(x) = \frac{1}{2}h^2 \left[ \int_{\mathbb{R}^d} y^{\otimes 2} K(y) dy\right]^T \nabla^{\otimes 2} f(x)\{1+o(1)\},\\
& \text{Var}(\wh f(x)) = \frac{1}{nh^d} c \int_{\mathbb{R}^d} K^2(x) dx\{1+o(1)\}.
\end{align*}

Clearly, the asymptotic minimizer $h_{\text{opt}}$ of $m_2(h)$ (or the exact minimizer of $\wt m_2(h)$), is in the form of $n^{-1/(d+4)}\lambda(f,g_1)/\lambda(f,g_2)$, where $g_1$ and $g_2$ are functions of the first and second derivatives of $f$. \\ 

Now we turn to the question of estimating $\lambda(f,g)$ when $g$ is unknown. Suppose $g$ has an estimator $\wh g$, we consider $\wt\lambda_{\epsilon_n}(\wh f,\wh g)$ as an estimator of $\lambda(f,g)$, where $\wt\lambda_{\epsilon_n}$ can be any one of $\lambda$, $\lambda_{\epsilon_n}^*$ and $\lambda_{\epsilon_n}^\dagger$. The specification of the function $g$ and its estimator $\wh g$ will of course impact the results we will get. Below we give a form of $g$ that is satisfied by all the examples (i)-(v). Note that all the principal curvatures $\kappa_1,\cdots,\kappa_{d-1}$ of $\cal M$ are the eigenvalues of the shape operator $S$ given in (\ref{shapeoperator}), which only depends on the first and second partial derivatives of $f$. This is also the requirement we have for $g$ that we consider below.\\

We need to introduce more notation. Let $\gamma=(\gamma_1,\cdots,\gamma_d)$ be a $d$-dimensional index, where $\gamma_i$'s are non-negative integers. Let $|\gamma| = \gamma_1+\cdots+\gamma_d$. For a smooth function $m:\mathbb{R}^d\mapsto\mathbb{R}$, define $$m^{(\gamma)}(x) = \frac{\partial^{|\gamma|} }{\partial^{\gamma_1}x_1\cdots \partial^{\gamma_d}x_d} m(x),\; x\in\mathbb{R}^d.$$
A closer examination reveals that in the examples (i)-(iv), the integrands $g$ in $\lambda(f,g)$ are all in the form of 
\begin{align}\label{goriginal}
g(x) = \wt\phi(\nabla f(x)) \prod_{j=1}^Q f^{(\beta_j)}(x),
\end{align}
for a function $\wt\phi:\mathbb{R}^d\mapsto \mathbb{R},$ a finite nonnegative interger $Q$ and $d$-dimensional indices $|\beta_j|=2,$ $j=1,\cdots,Q$, or a sum of finitely many such terms (see Remark~\ref{formg} below for more discussion). 
%
Let $d_{f}=(d_{f,1}^T,d_{f,2}^T)^T$ with $d_{f,1}=\nabla f$ and $d_{f,2} = \vech \nabla^2 f$, where $\vech$ is the half-vectorization operator, that is, it vectorizes the lower triangular portion of symmetric matrices. In other words, $d_{f,1}$ and $d_{f,2}$ include all the first and second partial derivatives of $f$, respectively. Let $a_d=d+d(d+1)/2$, which is the dimension of $d_f$. For $g$ in (\ref{goriginal}) we define a function $\phi:\mathbb{R}^{a_d}\mapsto\mathbb{R}$ such that we can write
\begin{align} \label{altgexpress}
g(x)=\phi(d_f(x)) = \wt\phi(d_{f,1}(x)) \prod_{j=1}^Q f^{(\beta_j)}(x).
\end{align}
Note that if $Q=0$, then $g$ only depends on $d_{f,1}$ but we still express $\phi$ as a function of $d_f$ for convenience. 
%
%
Let $d_{\wh f}(x)$ be the plug-in estimator of $d_f(x)$ and 
\begin{align}\label{altgexpressest}
\wh g(x) = \phi(d_{\wh f}(x)) = \wt\phi(d_{\wh f,1}(x)) \prod_{j=1}^Q \wh f^{(\beta_j)}(x).  
\end{align}
We need the following assumptions.\\

(P) For $g$ in (\ref{altgexpress}), $Q$ is a finite nonnegative interger and if $Q\geq 1$, we assume $|\beta_j|=2,$ for $j=1,\cdots,Q$. There exist constants $\delta_1>0$ and $C>0$ such that $\wt\phi(y)$ is $\nu\vee(Q+1)$ times differentiable for $y\in \mathcal{G}\oplus\delta_1$, where $\mathcal{G}=\{d_{f,1}(x):\;x\in\mathcal{M}\}$, and 
\begin{align}\label{nablajbound}
\sup_{y\in \mathcal{G}\oplus\delta_1} \|\nabla^{\otimes j} \wt\phi(y)\|_{\max} \leq C, \; j=1,2,\cdots,\nu\vee(Q+1).
\end{align}

(F2) The density function $f$ has bounded and continuous partial derivatives up to order $\nu+l_Q$ on $\mathcal{I}(2\delta_0)$, where $l_Q$ is defined through
\begin{align*}
l_Q = \begin{cases}
1 & \text{ if } Q \leq 0,\\
2 & \text{ if } Q \geq 1.
\end{cases}    
\end{align*}

\begin{remark}\label{formg}
%
{\em The assumption (P) is satisfied by all the examples (i)-(v) introduced at the beginning of this section. Notice that the integrands in these examples are finite degree polynomials of $1/\|\nabla f\|$, $\nabla f$, $\nabla^2 f$, and $F_j(\kappa_1,\cdots,\kappa_{d-1})$ for $j=1,\cdots,d$, where $F_j$'s are given in (\ref{minkfunction}). Also $F_j(\kappa_1,\cdots,\kappa_{d-1})$ is a finite degree polynomial of $1/\|\nabla f\|$, $\nabla f$, $\nabla^2 f$ because it is the sum of all the $j\times j$ principal minors of $S(x)$ (cf. Theorem 7.1.2 in Mirksy, 1995).} 
\end{remark}

Let $\nabla_1 \phi$ and $\nabla_2 \phi$ be the vectors of partial derivatives of $\phi$ with respect to the first $d$ variables, and the last $(a_d-d)$ variables, respectively. For $l=1,2$, let \begin{align}%
\sigma_l^2 = c\lambda(f,m_{K,l}) \text{ and } \wh\sigma_l^2 = c\lambda(\wh f,\wh m_{K,l}),\label{asyptoticvarexp}
\end{align}
where
\begin{align*}
&  m_{K,l}(x)=\int_{\mathbb{R}} \left(\int_{T_{x,t}(\mathcal{M})} \nabla_l\phi(d_f(x))^T d_{K,l}(u) d\mathscr{H}(u) \right)^2dt, \;   x\in\mathcal{M},\\
&  \wh m_{K,l}(x)=\int_{\mathbb{R}} \left(\int_{T_{x,t}(\wh{\mathcal{M}})} \nabla_l\phi( d_{\wh f}(x))^T d_{K,l}(u) d\mathscr{H}(u) \right)^2dt, \; x\in\wh{\mathcal{M}}.
\end{align*}


The following is the main result of this section, which gives the rates of convergence and asymptotic normality of $\wt\lambda_{\epsilon_n}(\wh f,\wh g) - \lambda(f,g)$.

\begin{theorem}\label{finaltheorem}
Let $g$ and $\wh g$ be given in (\ref{altgexpress}) and (\ref{altgexpressest}), respectively. Let $\wt\lambda_{\epsilon_n}$ be either $\lambda$, $\lambda_{\epsilon_n}^*$, $\lambda_{\epsilon_n}^\dagger$. Let $\alpha_{n,h,Q}:=\frac{1}{\sqrt{nh^{1+2l_Q}}} + \frac{1}{nh^{d+2l_{Q-1}}} + h^\nu + \epsilon_n^2=o(1)$. Under assumptions (F1), (F2), (K), (H), and (P) and if $K$ is four times continuously differentiable, we have
\begin{align}\label{asymptotunrat}
\wt\lambda_{\epsilon_n}(\wh f,\wh g) - \lambda(f,g)= O_p\left(\alpha_{n,h,Q}\right).    
\end{align}
If further $nh^{2d+4l_{Q-1}-2l_Q-1}\rightarrow\infty$, $nh^{2\nu+2l_{Q}+1}\rightarrow0$ and $\sqrt{nh^{1+2l_Q}}\epsilon_n^2\rightarrow0$, then 
\begin{align}\label{asymptotnormunkn}
\sqrt{nh^{1+2l_Q}}[\wt\lambda_{\epsilon_n}(\wh f,\wh g) - \lambda(f,g)]\rightarrow_D \mathscr{N}(0,\sigma_{l_Q}^2),   
\end{align}
where $\sigma_{l_Q}^2$ given in (\ref{asyptoticvarexp}), and if $\sigma_{l_Q}^2\neq0$,
\begin{align}\label{asymptotnormunknpl}
\sqrt{nh^{1+2l_Q}}\wh\sigma_{l_Q}^{-1}[\wt\lambda_{\epsilon_n}(\wh f,\wh g) - \lambda(f,g)]\rightarrow_D \mathscr{N}(0,1).    
\end{align}
%
\end{theorem}

\begin{remark}\label{discuss2}$\;$
%

{\em
a) More explicitly we can write  
\begin{align*}
 \alpha_{n,h,Q} = \begin{cases}
 \frac{1}{\sqrt{nh^{3}}} + \frac{1}{nh^{d+2}} + h^\nu + \epsilon_n^2 & \text{if } Q=0,\\
 \frac{1}{\sqrt{nh^{5}}} + \frac{1}{nh^{d+2}} + h^\nu + \epsilon_n^2 & \text{if } Q=1,\\
 \frac{1}{\sqrt{nh^{5}}} + \frac{1}{nh^{d+4}} + h^\nu + \epsilon_n^2 & \text{if } Q\geq 2.
 \end{cases}  
\end{align*}
It turns out that $\alpha_{n,h,Q}$ is also the rate of convergence of $\wt\lambda_{\epsilon_n}(f,\wh g) - \lambda(f,g)$ and rate of the difference $\wt\lambda_{\epsilon_n}(\wh f,\wh g) - \wt\lambda_{\epsilon_n}(f,\wh g)$ can be absorbed into $\alpha_{n,h,Q}$. The rates in $a_{n,h,Q}$ can be heuristically interpreted in a way similar to Remark~\ref{discuss}b for Theorem~\ref{HausdorffIntegration}. In particular, to understand the rate $1/\sqrt{nh^{1+2l_Q}}$, notice that  $\wh g(x) -\mathbb{E}\wh g(x)\approx \nabla\phi(d_f(x))^T[d_{\wh f}(x) - \mathbb{E} d_{\wh f}(x)]$ and thus has a rate of $1/\sqrt{nh^{d+2l_Q}}$. We then gain $(d-1)$ powers of $h$ through an integration on $\mathcal{M}$. The rate $1/nh^{d+2l_{Q-1}}$ is due to a quadratic form of $d_{\wh f,l_{Q-1}}(x) - \mathbb{E} d_{\wh f,l_{Q-1}}(x)$, i.e., the second term in the Taylor expansion of $\wh g(x) -\mathbb{E}\wh g(x)$. To get the asymptotic normality in (\ref{asymptotnormunkn}), we make the bias asymptotically negligible by assuming $nh^{2d+4l_{Q-1}-2l_Q-1}\rightarrow\infty$ and $nh^{2\nu+2l_{Q}+1}\rightarrow0$, which require $\nu>d+2l_{Q-1}-l_Q-1$.\\
b) The result in this theorem can be generalized if $\prod_{j=1}^Q f^{(\beta_j)}(x)$ in (\ref{altgexpress}) is replaced by a polynomial of second partial derivatives of $f$ and the plug-in estimator in (\ref{altgexpressest}) is changed accordingly. If so then the results in the theorem remain the same with $Q$ determined by the order of the polynomial. 
%
}
\end{remark}

Notice that 
\begin{align}\label{lambdadiffdecomp}
\wt\lambda_{\epsilon_n}(\wh f,\wh g) - \lambda(f,g)
=&\wt\lambda_{\epsilon_n}(\wh f,\wh g-g) + [ \wt\lambda_{\epsilon_n}(\wh f,g) -  \lambda(f,g)]\nonumber\\
= &\underbrace{[ \wt\lambda_{\epsilon_n}(\wh f,g) -  \lambda(f,g)]}_{\text{I}_n} +  \underbrace{[\wt\lambda_{\epsilon_n}(\wh f,\wh g-g)-\lambda( f,\wh g-g)]}_{\text{II}_n} + \underbrace{\lambda( f,\wh g-g)}_{\text{III}_n}.
\end{align}

This means the study of the difference $\wt\lambda_{\epsilon_n}(\wh f,\wh g) - \lambda(f,g)$ can be performed through the investigation of the three differences above. The difference I$_n$ is the main object we have studied in Section~\ref{knownintegrand}. The difference $\text{II}_n$ is also dealt with if we apply Corollary~\ref{Hausdorff} with $b_n = \wh g - g$. Below we study the difference III$_n$, which also corresponds to surface integral estimation when $\mathcal{M}$ is known but the integrand is unknown. 

\subsection{Surface integral estimation when $\mathcal{M}$ is known}\label{mknowncase}
The following theorem gives the rate of convergence and asymptotic normality of $\text{III}_n$ in (\ref{lambdadiffdecomp}). 
\begin{theorem}\label{asympototicnormalitymknown}
Suppose the assumptions (F1), (F2), (K), (H), and (P) hold. Then $\lambda(f,\mathbb{E}\wh g-g) = O(h^\nu+1/(nd^{d+2l_{Q-1}}))$, and
\begin{align}\label{asymp1stderiv}
\sqrt{nh^{1+2l_Q}} [\lambda(f,\wh g-\mathbb{E}\wh g)] \rightarrow_D \mathscr{N}(0,\sigma_{l_Q}^2),
\end{align}
where $\sigma_{l_Q}^2$ is given in (\ref{asyptoticvarexp}).
\end{theorem}

When $Q\geq2$, the bias $\lambda(f,\mathbb{E}\wh g-g)$ is of the order of $O(h^\nu + 1/nh^{d+4})$. It turns out that we can develop a U-statistic type estimator for $\lambda(f,g)$ with faster rate of convergence of the bias in some scenarios, while maintaining the asymptotic normality of the centered estimator. Our estimator for such $g(x)$ is given by
\begin{align}\label{gexpress2}
\wh g(x) =   \frac{(n-Q)!}{n! h^{(d+2)Q} } \sum_{(l_1,\cdots,l_Q)\in\pi_n^Q} \left[\wt\phi(d_{\wh f,1}^*(x)) \prod_{j=1}^Q K^{(\beta_j)}\left(\frac{x-X_{l_j}}{h}\right)  \right], 
\end{align}
where $\pi_n^Q$ is the set of all $Q$-element permutations of $\{1,\cdots,n\}$, and we use $d_{\wh f,1}^*(x)$ to denote the plug-in estimator of $d_{f,1}(x)$ using the data excluding $\{X_{l_1},\cdots,X_{l_Q}\}$. 
%

\begin{corollary}\label{asympototicnormalityspec}
Let $g$ and $\wh g$ be given in (\ref{altgexpress}) and (\ref{gexpress2}), respectively, with $Q\geq2$. Under the assumptions (F1), (F2), (K), (H), and (P), we have $\lambda(f,\mathbb{E}\wh g-g) = O(h^\nu + 1/(nd^{d+2}))$ and
\begin{align}\label{asympnormustat}
\sqrt{nh^5} [\lambda(f,\wh g-\mathbb{E}\wh g)] \rightarrow_D \mathscr{N}(0,\sigma^2),   
\end{align}
where $\sigma^2$ is given in (\ref{sigma2express}).
\end{corollary}

\begin{remark}\label{asympototicnormalityspecre}$\;$
%

{\em a) It is known that using U-statistic estimators can reduce the bias in estimating the integrated squared density derivatives on $\mathbb{R}^d$, compared with a direct plug-in estimator (see Hall and Marron, 1987 and Gin\'{e} and Nickl, 2008). We use a similar idea here but our estimator $\lambda(f,\wh g)$ based on $\wh g$ in (\ref{gexpress2}) is only partially like a U-statistic due to the possibly non-polynomial term $\phi(d_{f,1})$ in $g$. Using this estimator, the rate of convergence of the bias $\lambda(f,\mathbb{E}\wh g-g)$ is reduced from an order of $O(h^\nu + 1/(nd^{d+4}))$ in Theorem~\ref{asympototicnormalitymknown} to an order of $O(h^\nu + 1/(nd^{d+2}))$. If $g$ is completely a polynomial of the first and second partial derivatives of $f$, say, $g(x) = \prod_{i=1}^P f^{(\alpha_i)}(x) \prod_{j=1}^Q f^{(\beta_j)}(x)$ with $P+Q\geq 1$, and write
\begin{align*}
\wh g(x) = \frac{(n-P-Q)!}{n! h^{(d+1)P + (d+2)Q} } \sum_{(k_1,\cdots,k_P,l_1,\cdots,l_Q)\in \pi_n^{P+Q}}  \prod_{i=1}^P K^{(\alpha_i)}\left(\frac{x-X_{k_i}}{h}\right)  \prod_{j=1}^Q K^{(\beta_j)}\left(\frac{x-X_{l_j}}{h}\right) ,
\end{align*}
then the bias $\lambda(f,\mathbb{E}\wh g-g)$ is further reduced to $O(h^\nu)$ and an asymptotic normality result similar to (\ref{asympnormustat}) can be derived.\\
b) By (\ref{lambdadiffdecomp}) and using $\wh g$ in (\ref{gexpress2}), a result similar to Theorem~\ref{finaltheorem} with $1/(nd^{d+4}))$ replaced by $1/(nd^{d+2}))$ in $a_{n,h}$ for $Q\geq2$ is expected to hold. Here we need the uniform rates of convergence of $\wh g-g$ and its first and second derivatives when applying Corollary~\ref{Hausdorff}. Although not derived in this manuscript, its appears these uniform rates of convergence may be obtained using the approach in Dony and Mason (2008). 
}
\end{remark}

\subsection{Estimation of Euler characteristic of level sets}\label{eulercharact}

In this section we consider the estimation of Euler characteristic of level sets $\chi(\mathcal{M})$ when $d$ is odd ($\chi(\mathcal{M})=0$ when $d$ is even), which is a special case of $\lambda(f,g)$ studied in this section. Note that $\chi(\mathcal{M})$ is also one of the Minkowski functionals (see the beginning of Section~\ref{unknownintegrand}). Let $s_d= \frac{1\cdot3\cdots(d-2)}{(2\pi)^{(d-1)/2}} .$ By the generalized Gauss-Bonnet theorem (see Gray, 2004, Chapter 5), we have 
$$\chi(\mathcal{M}) =s_d \lambda(f,\kappa),$$
where $\kappa(x)$ is the Gaussian curvature of $\mathcal{M}$ at $x$, defined as the product of all the principal curvatures. For a level set $\mathcal{M}$, it is known that
\begin{align}\label{gaussiancurvature}
\kappa(x) = (-1)^{d-1} \frac{\n(x)^T[\nabla^2 f(x)]^\star\n(x)}{\|\nabla f\|^{d-1}},
\end{align}
where $[\nabla^2 f(x)]^\star$ denotes the adjugate matrix, i.e., the transpose of the cofactor matrix, of the Hessian. See Goldman (2005).\\

Let $\wh \kappa$ be the plug-in kernel estimator of the Gaussian curvature $\kappa$ given in (\ref{gaussiancurvature}), that is,
\begin{align*}
\wh \kappa(x)  = (-1)^{d-1}  \frac{\wh\n(x)^T[\nabla^2 \wh f(x)]^\star\wh \n(x)}{\|\nabla \wh f\|^{d-1}}. 
\end{align*}

The asymptotic performance of the three estimators $s_d\lambda(\wh f, \wh \kappa)$, $s_d\lambda_{\epsilon_n}^*(\wh f, \wh \kappa)$, and $s_d\lambda_{\epsilon_n}^\dagger(\wh f, \wh \kappa)$ for $\chi(\mathcal{M})$ has been considered in Theorem~\ref{finaltheorem}, using $\wh\kappa$ as a special case of $\wh g$. \\

Notice that $s_d\lambda(\wh f, \wh \kappa) = \chi(\wh{\mathcal{M}})$ by the generalized Gauss-Bonnet theorem, and using (\ref{altstar}) we have
\begin{align}\label{eulerestimator1}
s_d\lambda_{\epsilon_n}^*(\wh f, \wh \kappa) = \frac{1}{2\epsilon_n} \int_{-\epsilon_n}^{\epsilon_n} \chi(\wh{\mathcal{M}}_{c+\epsilon}) d\epsilon.
\end{align}

The function $\upsilon(\epsilon):=  \chi(\wh{\mathcal{M}}_{c+\epsilon})$ defines an Euler characteristic curve (see Turner et al., 2014) over $[-\epsilon_n, \epsilon_n]$. The estimator $s_d\lambda_{\epsilon_n}^*(\wh f, \wh \kappa)$ is understood as the average of $\chi(\wh{\mathcal{M}}_{c+\epsilon})$ for $\epsilon\in[-\epsilon_n, \epsilon_n]$. 
%
Inspired by this, we consider another estimator for $\chi(\cal M)$: $\frac{1}{2\epsilon_n} \int_{-\epsilon_n}^{\epsilon_n} \chi(\wh{\mathcal{M}}_{c}\uplus\epsilon) d\epsilon$, where ${\wh {\cal M}}_c\uplus\epsilon :=\{x+\epsilon\wh\n(x):\;x\in{\wh {\cal M}}_c\}$ is a parallel surface to ${\wh {\cal M}}_c$. Note here that the function $\wt\upsilon(\epsilon):=  \chi(\wh{\mathcal{M}}_{c}\uplus\epsilon)$ is also an Euler characteristic curve over $[-\epsilon_n, \epsilon_n]$. Using the generalized Gauss-Bonnet theorem again, we have
\begin{align}\label{eulerestimator2}
s_d\lambda_{\epsilon_n}^\dagger(\wh f, \wh \kappa^\sharp) = \frac{1}{2\epsilon_n} \int_{-\epsilon_n}^{\epsilon_n} \chi(\wh{\mathcal{M}}_{c}\uplus\epsilon) d\epsilon,
\end{align}
where $\wh \kappa^\sharp$ is the Gaussian curvature of ${\wh {\cal M}}_c\uplus\epsilon$. A similar idea of estimating some topological invariant of super level set $\mathcal{L}$ based on the topology of $\wh{\mathcal{M}}_{c+\epsilon}$ for a range of $\epsilon$ near zero is also used in Bobrowski et al. (2017). The following theorem describes a large sample behaviour of our estimators due to the property that $\chi(\mathcal{M})$ is integer-valued.\\ 

%

\begin{theorem}\label{eulercharactheorem}
Let $\wh {\chi ({\cal M})}$ be one of $s_d\lambda(\wh f, \wh \kappa)$, $s_d\lambda_{\epsilon_n}^*(\wh f, \wh \kappa)$, and $s_d\lambda_{\epsilon_n}^\dagger(\wh f, \wh \kappa^\sharp)$. Under assumptions (F1), (K), and (H), with probability one we have $\wh {\chi ({\cal M})} = \chi ({\cal M})$ for $n$ large enough.

\end{theorem}

\begin{remark}{\em
The computations of $s_d\lambda(\wh f, \wh \kappa)$, $s_d\lambda_{\epsilon_n}^*(\wh f, \wh \kappa)$, and $s_d\lambda_{\epsilon_n}^\dagger(\wh f, \wh \kappa^\sharp)$ are based on integration. Below we provide the expression of $\wh\kappa^\sharp$. For any $x\in{\wh {\cal M}}_c\oplus \epsilon_n$, let $\epsilon(x) = \text{sign}(\wh f(x)- c) d(x,{\wh {\cal M}}_c)$ so that $x\in {\wh {\cal M}}_c\uplus\epsilon(x)$. In other words, $\epsilon(x)$ is the signed distance from $x$ to ${\wh {\cal M}}_c$. Recall that $\pi_{{\wh {\cal M}}_c}(x)$ is the normal projection of $x$ onto ${\wh {\cal M}}_c$. For simplicity we write $\pi(x)=\pi_{{\wh {\cal M}}_c}(x)$. 
 It is known (c.f. Page 212, do Carmo, 1976) that the principal curvatures of ${\wh {\cal M}}_c \uplus \epsilon(x)$ at $x$ are 
 $$\frac{\wh\kappa_1(\pi(x))}{1+\epsilon(x)\wh\kappa_1(\pi(x))},\cdots,\frac{\wh\kappa_{d-1}(\pi(x))}{1+\epsilon(x)\wh\kappa_{d-1}(\pi(x))}.$$ 
Therefore $$\wh \kappa^\sharp(x) = \frac{\wh\kappa_1(\pi(x))\cdots \wh\kappa_{d-1}(\pi(x))}{\prod_{j=1}^{d-1}[1+\epsilon(x)\wh\kappa_j(\pi(x))]} = \frac{\wh \kappa(\pi(x))}{\prod_{j=1}^{d-1}[1+\epsilon(x)\wh\kappa_j(\pi(x))]} .$$}
\end{remark}

\section{Discussion}
In this manuscript we consider nonparametric estimation of $\lambda(f,g)$, which is a surface integral on density level set $f^{-1}(c)$, where the integrand $g$ can be known or unknown. We study three types of estimators: $\lambda(\wh f, g)$, $\lambda_{\epsilon_n}^*(\wh f, g)$ and $\lambda_{\epsilon_n}^\dagger(\wh f, g)$ if $g$ is known. If $g$ is unknown (but depends on the derivatives of $f$) and has an estimator $\wh g$, we replace $g$ by $\wh g$ in the above estimators. Among the three types of estimators, $\lambda(\wh f, g)$ (or $\lambda(\wh f, \wh g)$) is a direct plug-in estimator, and the estimators using $\lambda_{\epsilon_n}^*$ and $\lambda_{\epsilon_n}^*$ are based on neighborhoods of $\wh f^{-1}(c)$. Our main results are the rates of convergence and asymptotic normality for these estimators. We introduce a few examples that can use our results, such as the estimation of Willmore energy, Minkowski functionals and Euler characteristic of level sets, and bandwidth selection for level set estimation. In particular, for the estimation of Euler characteristic of level sets we provide the fourth estimator. \\

Our methods can be extended to level set estimation in regression and classification problems. For instance, Hall and Kang (2005) and Samworth (2012) show the optimal nearest neighbor classifiers depend on surface integrals on level sets of a regression function. We expect that our approach will work in the estimation of these surface integrals.\\

To make the bias have asymptotically negligible effect, the bandwidth used in the estimation of surface integrals on level sets need to be appropriately selected. This is especially important if one would like to construct confidence intervals for the population surface integrals using the asymptotic normality results. An alternative approach is to use debiased estimators to explicitly reduce the bias. See Chen (2017) and Qiao and Polonik (2019). However, this only has an effect on $h^\nu$ in the rates of convergence in our main theorems, while the rates due to the quadratic form of the derivatives (see Remarks~\ref{discuss}a and \ref{discuss2}a) are not impacted. It appears that one can remediate the latter by using U-statistic type estimators. See Corollary~\ref{asympototicnormalityspec} and Remark~\ref{asympototicnormalityspecre}.\\    

\section{Proofs}\label{proofsection}

{\bf Proof of Lemma~\ref{compatible}}
\begin{proof}
Using Lemma 3 in Arias-Castro et al. (2016), it follows by assumptions (K), (F1) and (H) that
\begin{align}
&\sup_{x\in\mathcal{I}(2\delta_0)} \|\widehat f -  \mathbb{E}\widehat f\| = O\left(\gamma_{n,h}^{(0)}\right), \; a.s. \label{supvar0} \\
&\sup_{x\in\mathcal{I}(2\delta_0)}\| \nabla\widehat f -  \mathbb{E}\nabla\widehat f\| = O\left(\gamma_{n,h}^{(1)}\right), \; a.s.  \label{supvar1}\\
&\sup_{x\in\mathcal{I}(2\delta_0)}\| \nabla^2\widehat f -  \mathbb{E}\nabla^2\widehat f\|_F = O\left(\gamma_{n,h}^{(2)}\right), \; a.s.  \label{supvar2}
\end{align} 
Also we follow a regular derivation procedure for kernel densities by applying change of variables and Taylor expansion, and obtain from assumptions (K), (F1) and (H) that
\begin{align}
&\sup_{x\in\mathcal{I}(2\delta_0)}| \mathbb{E}\widehat f(x) - f(x)| = O(h^\nu),\label{supbias0}\\
&\sup_{x\in\mathcal{I}(2\delta_0)}\| \mathbb{E}\nabla\widehat f(x) - \nabla f(x)\| = o(h^{\nu-1}),\label{supbias1}\\
&\sup_{x\in\mathcal{I}(2\delta_0)}\| \mathbb{E}\nabla^2\widehat f(x) - \nabla^2 f(x)\|_F = o(h^{\nu-2}).\label{supbias2}
\end{align}
This means we have strong uniform consistence of the first two derivatives of KDE under our assumptions. Assertions (i)-(iii) follow from similar arguments as given in the proof of Lemma 1 in Chen et al. (2017), which we briefly sketch here. Essentially bounded $\sup_{x\in\mathcal{I}(2\delta_0)}\|\nabla^2 f(x)\|_F$ and positive $\sup_{x\in\mathcal{I}(2\delta_0)}\|\nabla f(x)\|$ lead to positive $\inf_{\tau\in[c-\delta_0,c+\delta_0]}\rho(\mathcal{M}_\tau)$. With probability one $\sup_{x\in\mathcal{I}(2\delta_0)}\|\nabla^2 \wh f(x)\|_F$ is bounded and  $\sup_{x\in\mathcal{I}(2\delta_0)}\|\nabla \wh f(x)\|$ is positive for large sample once we use (\ref{supvar0}) - (\ref{supbias2}), which then implies $\inf_{\tau\in[c-\delta_0,c+\delta_0]}\rho(\wh{\mathcal{M}}_\tau)$ is positive. By Theorem 2 in Cuevas et al. (2006), we have 
\begin{align}\label{Hausdorffrate}
\sup_{\tau\in[c-\delta_0,c+\delta_0]} d_H(\wh{\mathcal{M}}_\tau, \mathcal{M}_\tau) = O\left( \sup_{x\in\mathcal{I}(2\delta_0)}\|\hat f - f\| \right)=O\left(\gamma_{n,h}^{(0)}\right) + O(h^\nu), \;a.s.
\end{align}
Therefore with probability one for $n$ large enough we have 
\begin{align*}
\sup_{\tau\in[c-\delta_0,c+\delta_0]} d_H(\wh{\mathcal{M}}_\tau, \mathcal{M}_\tau) \leq (2-\sqrt{2}) \inf_{\tau\in[c-\delta_0,c+\delta_0]} \min(\rho(\wh{\mathcal{M}}_\tau),\rho(\mathcal{M}_\tau)),
\end{align*}
which by Theorem 1 in Chazal et al. (2007) implies that $\wh{\mathcal{M}}_\tau$ and $\mathcal{M}_\tau$ are normal compatible for all $\tau\in[c-\delta_0,c+\delta_0]$. \hfill$\square$
\end{proof}

{\bf Proof of Lemma \ref{pointlevelset}}
\begin{proof}
%
%

We first focus on $t_n(x)$. The fact that $|t_n(x)| = \|P_n(x) - x\|$ and (\ref{Hausdorffrate}) imply
\begin{align}\label{tnx}
\sup_{x\in \mathcal{I}(\delta_0)} |t_n(x)| \leq \sup_{\tau\in[c-\delta_0,c+\delta_0]} d_H(\wh{\mathcal{M}_\tau}, \mathcal{M}_\tau) =o_p(1).
\end{align}
Since $\widehat f(P_n(x)) - f(x) =0$, using Taylor expansion for $\widehat f(P_n(x))$ we obtain
\begin{align}\label{tnxeq}
0 = \widehat f\left(x + t_n(x)\n(x) \right) - f(x) = \widehat f(x) -f(x) + t_n(x)  \n(x)^T\nabla \widehat f(x)+ e_n(x),
\end{align}
where $e_n(x) = \frac{1}{2} t_n(x)^2 \n(x)^T\nabla^2 \widehat f\left(x + s_1t_n(x)\n(x)\right) \n(x),$ for some $0<s_1<1$. Plugging $\nabla \widehat f(x) = \nabla f(x) - [\nabla f(x) - \nabla \widehat f(x)]$ into (\ref{tnxeq}), we have
\begin{align}\label{tnexpresss}
t_n(x) = \frac{1}{\|\nabla f(x)\| }[f(x) - \widehat f(x)] + \delta_{n,1}(x),
%
\end{align}
where
\begin{align*}
\delta_{n,1}(x) = \frac{1}{\|\nabla f(x)\| }\n(x)^T [\nabla f(x) - \nabla \widehat f(x)] t_n(x) - \frac{1}{\|\nabla f(x)\|}e_n(x).
\end{align*}
Under assumption (F1), we then have
\begin{align}\label{deltanineq}
&\sup_{x\in \mathcal{I}(\delta_0)} |\delta_{n,1}(x)|\nonumber\\
 \leq & \frac{1}{\epsilon_0} \sup_{x\in \mathcal{I}(\delta_0)} \|\nabla f(x) - \nabla \widehat f(x)\| \sup_{x\in \mathcal{I}(\delta_0)} |t_n(x)| + \frac{1}{2\epsilon_0} \sup_{x\in \mathcal{I}(\delta_0)} |t_n(x)|^2 \sup_{x\in \mathcal{I}(\delta_0)} \|\nabla^2\widehat f(x)\|_F.
\end{align}
Then by (\ref{supvar0}) - (\ref{supbias2}), (\ref{deltanineq}) implies that $\sup_{x\in \mathcal{I}(\delta_0)} |\delta_{n,1}(x)| = o_p\left( \sup_{x\in \mathcal{I}(\delta_0)} |t_n(x)| \right).$
%
Furthermore by (\ref{tnexpresss}) we have
\begin{align*}
\sup_{x\in \mathcal{I}(\delta_0)} |t_n(x)|  \leq \frac{1}{\epsilon_0} \sup_{x\in \mathcal{I}(\delta_0)} |\widehat f(x) - f(x)| + \sup_{x\in \mathcal{I}(\delta_0)} |\delta_{n,1}(x)| = O_p\left(\gamma_{n,h}^{(0)} + h^\nu\right).
 \end{align*}
Using  (\ref{deltanineq}) again, we obtain (\ref{delta1nrate}).\\

Next we study $\nabla t_n(x)$. Another Taylor expansion for $\widehat f(P_n(x))$ with fewer terms than (\ref{tnxeq}) gives
\begin{align}\label{tnxeq2}
0 = \widehat f\left(x + \n(x) t_n(x)\right) - f(x) = \widehat f(x) -f(x) + t_n(x) \n(x)^T r_n(x), 
\end{align}
where 
\begin{align}\label{rnexpress}
r_n(x) = \int_0^1 \nabla \widehat f\left(x + s\n(x)t_n(x)\right) ds.
\end{align}
Taking gradient on both sides of (\ref{tnxeq2}), we obtain
\begin{align*}
0 = & \nabla \widehat f(x) - \nabla f(x) + \nabla \left(\n(x)^T r_n(x)\right) t_n(x) +  \n(x)^T r_n(x)\nabla t_n(x)\\
= & \nabla \widehat f(x) - \nabla f(x) +  \|\nabla f(x)\|\nabla t_n(x) + s_n(x),
\end{align*}
where
\begin{align}\label{nbtn}
s_n(x) = \nabla \left(\n(x)^T r_n(x)\right) t_n(x) +  \n(x)^T [r_n(x)-\nabla f(x)]\nabla t_n(x).
\end{align}
Therefore
\begin{align}\label{nablatnequation}
\nabla t_n(x) = \frac{\nabla f(x) - \nabla \widehat f(x)}{\|\nabla f(x)\|} -\frac{s_n(x)}{\|\nabla f(x)\|}.
\end{align}

Let $u_{n}(x) = \nabla \left(\n(x)^T r_n(x) \right) - \nabla \|\nabla f(x)\|\text{ and }  v_{n}(x) = \n(x)^T [r_n(x) - \nabla f(x)].$ From (\ref{nbtn}) we have
\begin{align*}
-\frac{s_n(x)}{\|\nabla f(x)\|}  = \frac{ - (\nabla \|\nabla f(x)\|+u_{n}(x))t_n(x) + v_{n}(x)\nabla t_n(x) }{\|\nabla f(x)\| } .
\end{align*}

Plugging this into (\ref{nablatnequation}), we have
\begin{align*}
\nabla t_n(x) = \frac{\nabla f(x) - \nabla \widehat f(x)}{\|\nabla f(x)\|} - \frac{  \nabla \|\nabla f(x)\|}{\|\nabla f(x)\| } t_n(x) + \delta_{n,2}(x),
\end{align*}
where 
\begin{align}\label{delta2nexactexpress}
\delta_{n,2}(x) = \frac{-u_n(x)t_n(x) + v_{n}(x)\nabla t_n(x) }{\|\nabla f(x)\|} .
\end{align}

We then have $r_n(x) -\nabla f(x) = \nabla \widehat f(x) - \nabla f(x)  + \frac{1}{2} \nabla^2 \widehat f\left(x + s_1\n(x)t_n(x)\right) \n(x) t_n(x)$ by comparing (\ref{tnxeq}) and (\ref{tnxeq2}). 
Therefore
\begin{align}\label{vnrate}
\sup_{x\in \mathcal{I}(\delta_0)}|v_{n}(x)| = O_p\left(\gamma_{n,h}^{(1)}\right) + o_p(h^{\nu-1}).
\end{align}

Next we consider $u_{n}(x)$. Note that
\begin{align*}
u_{n}(x) & = \frac{\nabla^2 f(x) [r_n(x) - \nabla f(x)] - \nabla f(x)^T[r_n(x) - \nabla f(x)] \nabla \|\nabla f(x)\|}{\|\nabla f(x)\|^2} \\
& \hspace{1cm} + \frac{[\nabla r_n(x) -\nabla^2f(x)]\n(x)}{\|\nabla f(x)\|}.
\end{align*}
Using (\ref{rnexpress}) we have $\nabla r_n(x) = \int_0^1 \nabla^2 \widehat f\left(x + s\n(x)t_n(x)\right)\left( \mathbf{I}_d + s\nabla \left( \n(x)t_n(x)\right)\right) ds,$ which leads to $\sup_{x\in \mathcal{I}(\delta_0)}\|\nabla r_n(x) - \nabla^2 f(x) \|_F = O_p\left(\gamma_{n,h}^{(2)} + h^\nu\right).$ This implies that
\begin{align}\label{unrate}
\sup_{x\in \mathcal{I}(\delta_0)} \|u_{n}(x)\| = O_p\left(\gamma_{n,h}^{(2)}\right) + o_p(h^{\nu-2}).
\end{align}

Therefore by (\ref{delta2nexactexpress}), (\ref{vnrate}) and (\ref{unrate}) we have (\ref{delta2nrate}) and our proof is completed. \hfill$\square$
%
%
\end{proof}
%

{\bf Proof of Theorem~\ref{asympnorm}} 
\begin{proof}
Note that
\begin{align}\label{sumxi}
\lambda(f,w\times(\wh f - \mathbb{E}\wh f)) = \frac{1}{n}\sum_{i=1}^n [\xi_i -\mathbb{E}\xi_i],
\end{align}
where $\xi_i = \frac{1}{h^d}\int_{\mathcal{M}} w(x) K \left( \frac{x - X_i}{h}\right) d\mathscr{H}(x).$ Next we will calculate the variance of $\xi_i$. First we have
\begin{align*}
\mathbb{E}( \xi_i^2) 
& =  \frac{1}{h^{2d}} \int_{\mathbb{R}^d}\int_{\mathcal{M}}  w(x)  K \left(\frac{x-z}{h}\right)d \mathscr{H}(x)  \int_{\mathcal{M}} w(y)  K \left(\frac{y-z}{h}\right)d \mathscr{H}(y) f(z) dz\nonumber\\
& =\frac{1}{h^{d}} \int_{\mathcal{M}} \int_{\mathcal{M}} \int_{\mathbb{R}^d} w(x)  K\left(u\right)  w(y)  K\left(\frac{y-x}{h} + u \right) f(x- h u) du d\mathscr{H}(y)  d\mathscr{H}(x),
\end{align*}
where we have used the variable transformation $u=h^{-1}(x-z)$. Since $K$ is assumed to have bounded support, there exists a finite constant $C>0$ such that the kernel function $K$ has support contained in $\mathcal{B}(0,C/2)$. Then we can write $ \mathbb{E}( \xi_i^2) =\frac{1}{h^{d}} \int_{\mathcal{M}} U_G(x)\mathscr{H}(x) $, where
\begin{align*}
U_G(x) = & \int_{\mathcal{M}\cap \mathcal{B}(x,Ch)} \int_{\mathbb{R}^d} w(x)  K\left(u\right)  w(y)  K\left(\frac{y-x}{h} + u \right) f(x- h u) du d\mathscr{H}(y),\; x\in\mathcal{M}.
%
%
%
%
\end{align*}
%

Next we will study $U_G(x)$. Let $\mathcal{M}\ominus x = \{y-x: y\in\mathcal{M}\}$, i.e., the manifold obtained by shifting $\mathcal{M}$ so that $x$ becomes the origin. Then another variable transformation leads to
\begin{align*}
%
U_G(x) = & \int_{(\mathcal{M}\ominus x)\cap \mathcal{B}(0,Ch)} \int_{\mathbb{R}^d} w(x)  K\left(u\right)  w(x+y)  K\left(\frac{y}{h} + u \right) f(x- h u) du d\mathscr{H}(y).
\end{align*}

Without loss of generality, we may assume that $\mathbf{N}(x)=(0,\cdots,0,1)^\prime$ and $T_x(\mathcal{M})=\mathbb{R}^{d-1}\times\{0\}$. When $h$ is small enough, $(\mathcal{M}\ominus x)\cap \mathcal{B}(0,Ch)$ coincides with the graph $y=(y_1,\cdots,y_{d-1},0)\in T_x(\mathcal{M}) \mapsto \phi(y):=(y_1,\cdots,y_{d-1}, p(y))\in\mathbb{R}^d$ with a smooth function $\phi(y)$, where for $y\in[T_x(\mathcal{M})\cap \mathcal{B}(0,Ch)]$, $p$ has a quadratic approximation (see page 141, Lee 1997)
\begin{align}\label{quadratic}
p(y) =\frac{1}{2}\sum_{i=1}^{d-1} \kappa_i(x)\langle y,p_i(x)\rangle^2 + O(h^3).
\end{align}
where $\kappa_i(x)$, $i=1,\cdots,d-1$ are the principal curvatures of $\mathcal{M}$ at $x$, and $p_i(x)$, $i=1,\cdots,d-1$ are the corresponding principal directions. Both $\kappa_i(x)$ and $p_i(x)$ can be obtained as the eigenvalues and eigenvectors of the shape operator $S(x)$ in (\ref{shapeoperator}). Note that the rate in the big O term in (\ref{quadratic}) is uniform over $\mathcal{M}$. Immediately from (\ref{quadratic}) we have 
\begin{align}\label{tangentapp}
\sup_{y\in[T_x(\mathcal{M})\cap \mathcal{B}(0,Ch)]}\|\phi(y)-y\|=O(h^2).
\end{align}

Then $\phi$ defines a diffeomorphism between $T_x(\mathcal{M})\cap \mathcal{B}(0,Ch)$ and its image under $\phi$. The Jacobian determinant of $\phi$ is 
\begin{align}\label{jacobapp}
J_\phi(y) = \sqrt{1+\left(\frac{\partial }{\partial y_1}p(y)\right)^2+\cdots + \left(\frac{\partial }{\partial y_{d-1}}p(y)\right)^2} = 1+O(h),
\end{align}
uniformly in $\phi(T_x(\mathcal{M})\cap \mathcal{B}(0,Ch))$, where we have used (\ref{quadratic}). \\

Notice that $(\mathcal{M}\ominus x)\cap \mathcal{B}(0,Ch) \subseteq \phi(T_x(\mathcal{M})\cap \mathcal{B}(0,Ch))$, and the $U_G(x)$ remains unchanged if its domain of integration is changed from $(\mathcal{M}\ominus x)\cap \mathcal{B}(0,Ch)$ to $\phi(T_x(\mathcal{M})\cap \mathcal{B}(0,Ch))$. In view of this and by using (\ref{tangentapp}) and (\ref{jacobapp}), we get
\begin{align*}
U_G(x) =& \int_{T_x(\mathcal{M})\cap \mathcal{B}(0,Ch)} \int_{\mathbb{R}^d} w(x)  K\left(u\right)  w(x+\phi(y))  K\left(\frac{\phi(y)}{h} + u \right) f(x- h u) du J_\phi(y)d\mathscr{H}(y)\\
=& \int_{T_x(\mathcal{M}) \cap \mathcal{B}(0,Ch)} \int_{\mathbb{R}^d} w(x)  K\left(u\right)  w(x+y)  K\left(\frac{y}{h} + u \right) f(x- h u) du d\mathscr{H}(y) (1+o(1)),
\end{align*}
where the little o is uniform in $x\in\mathcal{M}$. 
%
Hence
\begin{align}\label{secondmoment}
& \mathbb{E}( \xi_i^2)  \nonumber\\
=&\frac{1}{h^{d}}  \int_{\mathcal{M}}  \int_{T_x(\mathcal{M}) \cap \mathcal{B}(0,Ch)} \int_{\mathbb{R}^d} w(x)  K\left(u\right)  w(x+y)  K\left(\frac{y}{h} + u \right) f(x- h u) du d\mathscr{H}(y)  d\mathscr{H}(x) (1+o(1)) \nonumber\\
=&\frac{1}{h}  \int_{\mathcal{M}}  \int_{T_x(\mathcal{M})} \int_{\mathbb{R}^d} w(x)  K\left(u\right)  w(x+hv)  K\left(v+ u \right) f(x) du d\mathscr{H}(v)  d\mathscr{H}(x) (1+o(1)) \nonumber\\
%
%
=& \frac{1}{h} c \rho(w,K,f)   (1+o(1)),
\end{align}
where $\rho(w,K,f) = \int_{\mathcal{M}} w(x)^2 \int_{T_x(\mathcal{M})} \int_{\mathbb{R}^d}  K\left(u\right)  K\left(v+ u \right) du d\mathscr{H}(v)\mathscr{H}(x).$ Notice that $\mathbb{R}^d=\cup_{t\in\mathbb{R}} T_{x,t}(\mathcal{M})$ for any $x\in\mathcal{M}$. So we can further write 
\begin{align*}
\rho(w,K,f) = & \int_{\mathcal{M}} w(x)^2 \int_{T_x(\mathcal{M})} \int_{\mathbb{R}}\int_{T_{x,t}(\mathcal{M})} K\left(u\right)  K\left(v+ u \right) d\mathscr{H}(u)dtd\mathscr{H}(v)\mathscr{H}(x)\\
=& \int_{\mathcal{M}} w(x)^2  \int_{\mathbb{R}}\int_{T_{x,t}(\mathcal{M})} K\left(u\right) \int_{T_x(\mathcal{M})}   K\left(v+ u \right) d\mathscr{H}(v) d\mathscr{H}(u)dt\mathscr{H}(x)\\
=& \int_{\mathcal{M}} w(x)^2  \int_{\mathbb{R}} \left(\int_{T_{x,t}(\mathcal{M})} K\left(u\right)    d\mathscr{H}(u)\right)^2dt\mathscr{H}(x),
\end{align*}
where the last step holds because we use change of variables $z=u+v$ and $z\in T_{x,t}(\mathcal{M})$ when $u\in T_{x,t}(\mathcal{M})$ and $v\in T_x(\mathcal{M})$. 
%
Similarly, using change of variables and Taylor expansion we get
\begin{align}\label{expectationsq}
\mathbb{E} ( \xi_i ) &= \frac{1}{h^{d}} \int_{\mathcal{M}} \int_{\mathbb{R}^d} w(x)  K \left(\frac{x-z}{h}\right) f(z) dz d\mathscr{H}(x)\nonumber\\
& =  \int_{\mathcal{M}} \int_{\mathbb{R}^d} w(x)  K \left(u\right) f(x-hu) du\mathscr{H}(x)\nonumber\\
& =   \int_{\mathcal{M}} \int_{\mathbb{R}^d} w(x)  K \left(u\right) f(x) dud\mathscr{H}(x)(1+o(1))\nonumber\\
&=  c \int_{\mathcal{M}} w(x) \mathscr{H}(x)(1+o(1)).
\end{align}

From (\ref{secondmoment}) and (\ref{expectationsq}) we have $\text{Var}(\xi_i) = h^{-1} c \rho(w,K,f) \{1+o(1)\}$ and therefore by (\ref{sumxi}),
\begin{align*}
\text{Var}\left(\int_{\mathcal{M}} w(x) \left[  \widehat f(x) - \mathbb{E}\widehat f(x) \right]dx\right) = \frac{\text{Var}(\xi_i)}{n} = \frac{1}{nh} c \rho(w,K,f) (1+o(1)).
\end{align*}
It only remains to show the asymptotic normality of $ \frac{1}{n} \sum_{i=1}^n (\xi_i - \mathbb{E} \xi_i)$. Note that following a similar procedure of obtaining (\ref{secondmoment}), we can show that $\mathbb{E}(|\xi_i|^3) = O(\frac{1}{h^{3d}} h^{3d-2}) = O(\frac{1}{h^2})$. Then for the third absolute central moment of $\xi_i$, we have
\begin{align*}
\mathbb{E} [|\xi_i - \mathbb{E}(\xi_i)|^3] \leq 8  \mathbb{E}(|\xi_i|^3) = O\left(\frac{1}{h^2}\right).
\end{align*}
Therefore
\begin{align*}
\frac{\{\sum_{i=1}^n \mathbb{E} [|\xi_i - \mathbb{E}(\xi_i)|^3]\}^{1/3}}{[\sum_{i=1}^n \text{Var}(\xi_i)]^{1/2}} = O\left(\frac{(n/h^2)^{1/3}}{(n/h)^{1/2}}\right) = O\left(\frac{1}{n^{1/6}h^{1/6}}\right) = o(1).
\end{align*}
Hence the Liapunov condition is satisfied and the asymptotic normality in (\ref{asymptocnormdens}) is verified. \\

\hfill$\square$
\end{proof}
%

{\bf Proof of Theorem~\ref{HausdorffIntegration}} 
\begin{proof}
We only show the proof for the case of $\wt\lambda_{\epsilon_n}=\lambda$. Once this is finished, the results for the cases of $\wt\lambda_{\epsilon_n}=\lambda_{\epsilon_n}^*$ or $\lambda_{\epsilon_n}^\dagger$ are direct consequences of Theorem~\ref{lambdastardaggertheorem} by using $\wt\lambda_{\epsilon_n}(\wh f, g) - \lambda(f,g) = [\wt\lambda_{\epsilon_n}(\wh f, g)-\lambda(\wh f, g)] + [\lambda(\wh f, g)-\lambda(f, g)].$\\

It is known that $\mathcal{M}$ is a compact $(d-1)$-dimensional submanifold embedded in $\mathbb{R}^d$ with assumption (F1) (see Theorem 2 in Walther, 1997). It admits an atlas $\{(U_\alpha,\psi_\alpha):\alpha\in\mathscr{A}\}$ indexed by a finite set $\mathscr{A}$, where $\{U_\alpha:\alpha\in\mathscr{A}\}$ is an open cover of $\mathcal{M}$, and for an open set $\Omega_\alpha\subset \mathbb{R}^{d-1}$, $\psi_\alpha: \Omega_\alpha \mapsto U_\alpha$ is a diffeomorphism. We denote $\psi_\alpha(\Omega_\alpha)=U_\alpha$. Let $B_\alpha$ be the Jacobian matrix of $\psi_\alpha$ and  $J_\alpha = B_\alpha^TB_\alpha$ be the Gram matrix. Immediately, we have
\begin{align}\label{boundedJac}
0<\inf_{\alpha\in\mathscr{A}}\inf_{\theta\in\Omega_\alpha} \det [J_\alpha(\theta)] \leq \sup_{\alpha\in\mathscr{A}}\sup_{\theta\in\Omega_\alpha} \det [J_\alpha(\theta)] <\infty.
\end{align}

By Lemma~\ref{pointlevelset}, $P_n(x) = x + t_n(x) \n(x)$ in (\ref{map}) is a diffeomorphism between $\mathcal{M}_\tau$ and $\wh{\mathcal{M}_\tau}$ for $|\tau-c|\leq\delta_0/2$. Let $\wh U_\alpha =\{P_n(x):x\in U_\alpha\}$. Then $\{(\wh U_\alpha,P_n\circ \psi_\alpha):\; \alpha\in\mathscr{A}\}$ is a finite atlas for $\wh{\mathcal{M}}$.\\

For simplicity, we will suppress the subscript $\alpha$ for $U_\alpha$, $\wh U_\alpha$, $\Omega_\alpha$, $\psi_\alpha$, $B_\alpha$ and $J_\alpha$.
By using the area formula on manifolds (cf. page 117, Evans and Gariepy, 1992), we have that 
\begin{align}\label{area1}
\int_{U} g(x) d\mathscr{H}(x) =   \int_{\Omega}g(\psi(\theta)) \Large\{\det [ J(\theta)] \Large\}^{1/2}d\theta.
\end{align}

Let $A_n$ be the Jacobian matrix of $P_n$. Then we can write $A_n(x) = \mathbf{I}_d +  R_n(x)$, $x\in \mathcal{I}(\delta_0)$, where $R_n(x) = \nabla\n(x) t_n(x) + \n(x)  \nabla t_n(x)^T.$ We further write $R_n(x)=R_{n,1}(x) + R_{n,2}(x)$, where 
\begin{align}\label{Rn1def}
R_{n,1}(x)=\|\nabla f(x)\|^{-1}\n(x)(\nabla f(x)-\nabla \wh f(x))^T,  
\end{align}
and $R_{n,2}(x)$ is the remainder term.
%
%
By Lemma \ref{pointlevelset} and using (\ref{supvar1}) and (\ref{supbias1}), we have
\begin{align}\label{tildeAmax}
\sup_{x\in \mathcal{I}(\delta_0)} \|R_{n,1}(x)\|_F = O_p\left(\gamma_{n,h}^{(1)}\right) + o_p\left(h^{\nu-1}\right),
\text{ and } \sup_{x\in \mathcal{I}(\delta_0)} \|R_{n,2}(x)\|_F = O_p\left(\gamma_{n,h}^{(0)}\right) + O_p\left(h^{\nu}\right).
\end{align}

By the chain rule, the Jabobian matrix of $P_n \circ \psi: \Omega\mapsto \wh{U}$ is given by $A_n(\psi(\theta)) B(\theta)$. Since $P_n \circ \psi$ is a parameterization function of $\wh{U}$ from $\Omega$, another application of the area formula leads to
\begin{align}\label{area2}
\int_{\wh U} g(x) d\mathscr{H}(x) =   \int_{\Omega}g(P_n(\psi(\theta))) \left\{\det \left[J_n(\theta)\right] \right\}^{1/2}d\theta,
\end{align}
where $J_n(\theta) = B(\theta)^TA_n(\psi(\theta))^T A_n(\psi(\theta)) B(\theta) $. 
%
%
%
Then it follows from (\ref{area1}) and (\ref{area2}) that
\begin{align}\label{hausdiff}
%
&\int_{\wh U} g(x) d\mathscr{H}(x) - \int_U g(x) d\mathscr{H}(x) \nonumber\\
= &\int_{\Omega}\left[ g(P_n(\psi(\theta))) \left\{\det \left[J_n(\theta)\right] \right\}^{1/2} - g(\psi(\theta)) \left\{\det [J(\theta) ] \right\}^{1/2}\right]d\theta\nonumber\\
=& \text{I}_n + \text{II}_n + \text{III}_n,
\end{align}
where 
\begin{align}
&\text{I}_n =  \int_{\Omega} [g(P_n(\psi(\theta))) - g(\psi(\theta))] \left\{\det [J(\theta) ] \right\}^{1/2}d\theta, \label{HausInexpr}\\
&\text{II}_n = \int_{\Omega} g(\psi(\theta)) \left\{ \left\{\det \left[J_n(\theta)\right] \right\}^{1/2} -  \left\{\det [J(\theta) ] \right\}^{1/2} \right\}d\theta,\label{HausIInexpr}\\
&\text{III}_n =  \int_{\Omega} [g(P_n(\psi(\theta))) - g(\psi(\theta))] \left\{ \left\{\det \left[J_n(\theta)\right] \right\}^{1/2} -  \left\{\det [J(\theta) ] \right\}^{1/2} \right\} d\theta. \label{HausIIInexpr}
\end{align}

We first study $\text{I}_n$. Notice that $\text{I}_n = \int_{U} [g(P_n(x)) - g(x)] d\mathscr{H}(x)$. Since $P_n(x) - x= t_n(x)\n(x) $ by definition, using Taylor expansion and Lemma \ref{pointlevelset} we have that
\begin{align}\label{gdiff}
\text{I}_n = \int_U \left[ \frac{f(x) - \widehat f(x)}{\| \nabla f(x) \| }\n(x)^T\nabla g\left(x \right) \right] d\mathscr{H}(x) + L_n,
%
\end{align}
where 
\begin{align*}
L_n = \int_U \left[  \delta_{n,1}(x)\n(x)^T\nabla g\left(x \right) + \frac{1}{2} t_n^2(x)\n(x)^T\nabla^2 g\left(x + s \n(x)t_n(x)\right) \n(x) \right] d\mathscr{H}(x),
\end{align*}
%
%
%
%
for some $0<s<1$ and $\delta_{n,1}$ is given in (\ref{delta1nexpr}). We apply Theorem~\ref{asympnorm} to the integrand of $L_n$ and then obtain
\begin{align}\label{gdiffrate}
L_n = O_p\left(\gamma_{n,h}^{(0)}\;\gamma_{n,h}^{(1)}\right) +o_p\left(h^{\nu-1}(\gamma_{n,h}^{(0)} + h^\nu)\right).
\end{align}

Next we focus on $\text{II}_n$. Notice that by using $A_n(x) = \mathbf{I}_d + R_n(x)$, we have
\begin{align}\label{determinant0}
 \det \left[J_n(\theta) \right] = & \det \left\{B(\theta)^T[ \mathbf{I}_d + R_n(\psi(\theta))]^T [\mathbf{I}_d+ R_n(\psi(\theta))] B(\theta) \right\}\nonumber\\
%
%
= & \det \left\{ J(\theta) + B(\theta)^T [R_n(\psi(\theta))^T + R_n(\psi(\theta)) +  R_n(\psi(\theta))^T R_n(\psi(\theta)) ] B(\theta) \right\} \nonumber\\
= & \det \left[J(\theta)\right] \det\left[\mathbf{I}_{d-1} - T_n(\theta)\right],
\end{align}
where $T_n(\theta) = - [J(\theta)]^{-1}B(\theta)^T\left\{  [\nabla\n(\psi(\theta)) + \nabla\n(\psi(\theta))^T ] t_n(\psi(\theta)) +  R_n(\psi(\theta))^T R_n(\psi(\theta))\right\}B(\theta).$ Notice that we have used the fact that  $\n(\psi(\theta))^TB(\theta)=0$, which can be seen by taking derivatives of the equation $f(\psi(\theta))=c$, $\theta\in\Omega$. 
%
Therefore 
\begin{align}\label{determinediff}
\{ \det \left[J_n(\theta) \right] \}^{1/2} - \{ \det \left[J(\theta) \right] \}^{1/2} =  \{ \det \left[J(\theta) \right] \}^{1/2}  \left\{  \{\det\left[\mathbf{I}_{d-1} - T_n(\theta)\right]\}^{1/2} - 1\right\}.
\end{align}

Let 
\begin{align}\label{qndef}
Q_n(x) = T_n(\psi^{-1}(x)) = Q_{n,1}(x)+Q_{n,2}(x), 
\end{align}
where with $R_{n,1}$ given in (\ref{Rn1def}), $Q_{n,1}(x)=- [J(\psi^{-1}(x))]^{-1}B(\psi^{-1}(x))^TR_{n,1}(x)^TR_{n,1}(x)B(\psi^{-1}(x))$ and $Q_{n,2}(x)$ is the remainder term. Apply the area formula with (\ref{determinediff}) we get $$\text{II}_n = \int_U  \left\{  \{\det\left[\mathbf{I}_{d-1} - Q_n(x)\right]\}^{1/2} - 1\right\} g(x) d\mathscr{H}(x).$$ Let $D(\theta) = B(\theta)[J(\theta)]^{-1}B(\theta)^T$ for $\theta\in\Omega$.  Continuing from (\ref{determinant0}), we have 
\begin{align}\label{determinant}
&\det\left[\mathbf{I}_{d-1} - Q_n(x)\right] \nonumber\\
=&1 - {\bf tr} \left[J(\psi^{-1}(x))^{-1}B(\psi^{-1}(x))^T  [\nabla\n(x) + \nabla\n(x)^T ] B(\psi^{-1}(x))\right] t_n(x) + \eta_n(x),
%
%
%
%
\end{align}
where $\eta_n(x) = \eta_{n,1}(x)+ \eta_{n,2}(x)$ with $\eta_{n,1}(x) = - {\bf tr}[D(\psi^{-1}(x)) R_n(x)^T R_n(x)]$ and $\eta_{n,2}(x) = \det(\mathbf{I}_{d-1}-Q_n(x) ) - [1 -  {\bf tr}(Q_n(x) )]$. 
Next we study $\eta_{n,1}$ and $\eta_{n,2}$. We write $\eta_{n,1}=\eta_{n,1}^{(1)}+\eta_{n,1}^{(2)}$ where $\eta_{n,1}^{(1)}(x)=- {\bf tr}[D(\psi^{-1}(x)) R_{n,1}(x)^T R_{n,1}(x)]$ and $\eta_{n,1}^{(2)}(x)$ is the remainder term. It follows from (\ref{tildeAmax}) that 
\begin{align}
&\sup_{x\in U} |\eta_{n,1}^{(1)}(x)| =  O_p\left((\gamma_{n,h}^{(1)})^2\right) +o_p\left(h^{\nu-1}(\gamma_{n,h}^{(0)} + h^{\nu-1})\right),\label{Delta1ratepart1}\\
&\sup_{x\in U} |\eta_{n,1}^{(2)}(x)| =  O_p\left(\gamma_{n,h}^{(0)}\;\gamma_{n,h}^{(1)}\right) +o_p\left(h^{\nu-1}(\gamma_{n,h}^{(0)} + h^\nu)\right).\label{Delta1ratepart2}
\end{align}
Then notice that $\eta_{n,2}(\theta)=0$ for $d=2$ and for $d\geq3$, $\eta_{n,2}(\theta) =\sum_{j=2}^{d-1} (-1)^j\mathbf{tr}(\Lambda^j Q_n(x))$, where $\mathbf{tr}(\Lambda^j Q_n(x))$ is the trace of the $j$th exterior power of $Q_n(x) $, which can be computed as the sum of all the $j\times j$ principal minors of $Q_n(x)$ ( Theorem 7.1.2, Mirsky, 1955). We write $\eta_{n,2}=\eta_{n,2}^{(1)}+\eta_{n,2}^{(2)}$, where $\eta_{n,2}^{(1)}(x)=\mathbf{tr}(\Lambda^2 Q_{n,1}(x))$ and $\eta_{n,2}^{(2)}(x)$ is the remainder term. Using (\ref{tildeAmax}) and the perturbation bound theory for matrix determinants (Theorem 3.3, Ipsen and Rehman, 2008), we have that 
\begin{align}
&\sup_{x\in U} |\eta_{n,2}^{(1)}(x)| =  O_p\left((\gamma_{n,h}^{(1)})^2\right) +o_p\left(h^{\nu-1}(\gamma_{n,h}^{(0)} + h^{\nu-1})\right),\label{Delta2ratepart1}\\
&\sup_{x\in U} |\eta_{n,2}^{(2)}(x)| =  O_p\left(\gamma_{n,h}^{(0)}\;\gamma_{n,h}^{(1)}\right) +o_p\left(h^{\nu-1}(\gamma_{n,h}^{(0)} + h^\nu)\right).\label{Delta2ratepart2}
\end{align}

Let $\eta_{n}^{(1)}(x)=\eta_{n,1}^{(1)}(x)+\eta_{n,2}^{(1)}(x)$ and $\eta_{n}^{(2)}(x)=\eta_{n,1}^{(2)}(x)+\eta_{n,2}^{(2)}(x)$. Then $\eta_{n}^{(1)}(x)$ is a quadratic form of elements of $\nabla f(x)-\nabla\wh f(x)$, similar to $T_{n,2}(x)$ defined in (\ref{Tnjdef}). Following similar arguments as in the proof of Theorem~\ref{asympototicnormalitymknown}, we get 
\begin{align}\label{etaiint}
\int_U \eta_{n}^{(1)}(x) d\mathscr{H}(x) = O_p\left(\frac{1}{nh^{d+2}}\right) + o_p\left(h^{2(\nu-1)}\right).    
\end{align}

Using Taylor expansion and (\ref{determinant}), we have
\begin{align}\label{determinant2}
& \{\det\left[\mathbf{I}_{d-1} - Q_n(x)\right]\}^{1/2} - 1\nonumber\\
 = & \left[1 + {\bf tr} \left[D(\psi^{-1}(x))  [\nabla\n(x) + \nabla\n(x)^T ] \right] t_n(x)  + \eta_n(x) \right]^{1/2}-1\nonumber\\
 = & \frac{1}{2}  {\bf tr} \left[D(\psi^{-1}(x))  [\nabla\n(x) + \nabla\n(x)^T ] \right] t_n(x) + \frac{1}{2} \eta_{n}^{(1)}(x) + \Delta_n(x), 
 %
%
%
 \end{align}
where $\Delta_n(x)$ denotes the remainder term in the Taylor expansion. Specifically,
\begin{align*}
\Delta_n(x) = &\frac{1}{2}  \eta_{n}^{(2)}(x)  + \sum_{i=2}^\infty (-1)^i\frac{(2i-3)!!}{2^ii!}  \left[ {\bf tr} \left[D(\psi^{-1}(x))  [\nabla\n(x) + \nabla\n(x)^T ] \right] t_n(x)  + \eta_n(x) \right]^i.
\end{align*}
It follows from (\ref{Delta1ratepart1}) - (\ref{Delta2ratepart2}) that
 \begin{align}\label{Delta3rate}
\sup_{x \in U} |\Delta_n(x)| =  O_p\left(\gamma_{n,h}^{(0)}\;\gamma_{n,h}^{(1)}\right) +o_p\left(h^{\nu-1}(\gamma_{n,h}^{(0)} + h^\nu)\right).
\end{align}

Using (\ref{etaiint}), (\ref{determinant2}), (\ref{Delta3rate}) and Theorem~\ref{asympnorm}, we have 
\begin{align}\label{IInexpress1}
\text{II}_n 
= & \frac{1}{2}\int_U {\bf tr} \left[D(\psi^{-1}(x))  [\nabla\n(x) + \nabla\n(x)^T ]\right] \frac{f(x) - \wh f(x)}{\|\nabla f(x)\|} g(x)  d\mathscr{H}(x) \nonumber \\
&\hspace{2cm}+ O_p\left(\frac{1}{nh^{d+2}} \right) +o_p\left(h^{\nu-1}(\gamma_{n,h}^{(0)} + h^{\nu-1})\right).
\end{align}

 In the integrand above, 
 \begin{align*}
& {\bf tr} \left[D(\psi^{-1}(x)) [\nabla\n(x) + \nabla\n(x)^T ] \right] \\
= &{\bf tr} \left[J(\psi^{-1}(x))^{-1} B(\psi^{-1}(x))^T  [\nabla\n(x) + \nabla\n(x)^T ] B(\psi^{-1}(x)) \right],     
 \end{align*}
which in fact does not depend on the parameterization $\psi$, as we show below. Recall that $G(x) = \mathbf{I}_d - \n(x)\n(x)^T$ and (\ref{Kindlmann}). Note that for any $u,v\in T_x(\mathcal{M})$, we have that 
\begin{align}\label{shapeoppro}
u^T S(x) v =  u^T \nabla \n(x) v = u^T \nabla \n(x)^T v .
\end{align}
The above equality can be derived by observing that the matrix $G(x)$ is symmetric and has the property that for any $u\in T_x(\mathcal{M})$, $G(x)u=u$ and $u^TG(x)=u^T$. 
Notice that all the columns of the Jacobian matrix $B(\psi^{-1}(x))$ are vectors in $T_x(\mathcal{M})$. Therefore using (\ref{shapeoppro}) we have
\begin{align*}
&\frac{1}{2}\mathbf{tr} [J( \psi^{-1}(x))^{-1} B(\psi^{-1}(x))^T [\nabla\n(x) + \nabla\n(x)^T ] B( \psi^{-1}(x))] \\
=& \mathbf{tr} [J( \psi^{-1}(x))^{-1} B(\psi^{-1}(x))^TS(x)B( \psi^{-1}(x))].
\end{align*}
%
Since $S(x)$ is symmetric and $S(x)\n(x)=0$, there exists an orthonormal basis $\{p_1(x),\cdots,p_{d-1}(x)\}$ for $T_x(\mathcal{M})$ such that with 
$P(x)=[p_1(x),\cdots,p_{d-1}(x),\n(x)]$, 
\begin{align*}
S(x) = P(x) K(x) P(x)^T,
\end{align*}
where $K(x)=\text{diag}(\kappa_{1}(x),\cdots,\kappa_{d-1}(x),0)$,
with $\kappa_1(x),\cdots,\kappa_{d-1}(x)$ being the principal curvatures of $\mathcal{M}$ at $x$. Let $\wt K(x)=\text{diag}(\kappa_{1}(x),\cdots,\kappa_{d-1}(x))$, $\wt B(x) = B( \psi^{-1}(x)) [J( \psi^{-1}(x))]^{-1/2}$, and $\check B(x) = [\wt B(x),\n(x)]$. 
%
Notice that $\n(x)^T\wt B(x)=0$ and therefore $\check B(x)^T \check B(x) = \mathbf{I}_d$. Then
\begin{align*}
& \mathbf{tr} [J( \psi^{-1}(x))^{-1} B(\psi^{-1}(x))^TS(x)B( \psi^{-1}(x))]\\
%
= & \mathbf{tr} [ \wt B(x)^TP(x) K(x) P(x)^T  \wt B(x)]\\
= &  \mathbf{tr} \left[\begin{pmatrix}
\wt B(x)^T\wt P(x) & 0
\end{pmatrix} 
\begin{pmatrix}
\wt K(x) & 0\\
0 & 0
\end{pmatrix} 
\begin{pmatrix}
\wt P(x)^T\wt B(x)\\
0
\end{pmatrix} \right]\\
= & \mathbf{tr} \left[\begin{pmatrix}
\wt B(x)^T\wt P(x) & 0\\
0 & 1
\end{pmatrix} 
\begin{pmatrix}
\wt K(x) & 0\\
0 & 0
\end{pmatrix} 
\begin{pmatrix}
\wt P(x)^T\wt B(x) & 0\\
0 & 1
\end{pmatrix} \right]\\
 = &  \mathbf{tr} [ \check B(x)^TP(x) K(x) P(x)^T  \check B(x)]\\
 = & \mathbf{tr} [ K(x)]\\
 = & (d-1) H(x),
 \end{align*}
where $H(x) =\frac{1}{d-1} \sum_{i=1}^{d-1}\kappa_{i}(x)$ is the mean curvature of $\mathcal{M}$ at $x$. 
%
%
Therefore from (\ref{IInexpress1}) we get
\begin{align}\label{IInasymptotic}
\text{II}_n = (d-1)  \int_U  H(x) \frac{f(x) - \wh f(x)}{\|\nabla f(x)\|} g(x)  d\mathscr{H}(x) +  O_p\left(\frac{1}{nh^{d+2}} \right) +o_p\left(h^{\nu-1}(\gamma_{n,h}^{(0)} + h^{\nu-1})\right).
\end{align}
Next we study $\text{III}_n$. Notice that $\text{III}_n = \int_U [g(P_n(x)) -g(x)]  \left\{  \{\det\left[\mathbf{I}_d - Q_n(x)\right]\}^{1/2} - 1\right\}  d\mathscr{H}(x)$.  

Lemma~\ref{pointlevelset} leads to 
\begin{align}\label{gdiffrate}
\sup_{x\in U} |g(P_n(x)) - g(x)| \leq \sup_{x\in \mathcal{I}(\delta_0)}\|\nabla g(x)\| \sup_{x\in\mathcal{M}} |t_n(x)| = O_p\left(\gamma_{n,h}^{(0)} + h^\nu\right). 
\end{align}

It follows from (\ref{determinant2}), (\ref{Delta3rate}) and Lemma~\ref{pointlevelset} that 
\begin{align}\label{detdiffrate}
\sup_{x\in U} \left|  \{\det\left[\mathbf{I}_{d-1} - Q_n(x)\right]\}^{1/2} - 1 \right| = O_p\left(\gamma_{n,h}^{(0)} + h^\nu  \right) . 
\end{align}

By (\ref{gdiffrate}) and (\ref{detdiffrate}) we have
\begin{align}\label{gdiffdetrate}
 \left| \text{III}_n\right| \leq & \mathscr{H}(U) \sup_{x\in U} |g(P_n(x)) - g(x)| \sup_{x\in U} \left|  \{\det\left[\mathbf{I}_{d-1} - Q_n(x)\right]\}^{1/2} - 1 \right| 
%
%
 =  O_p\left((\gamma_{n,h}^{(0)} + h^\nu)^2\right) .
\end{align}

Now collecting the above results for $\text{I}_n$ in (\ref{gdiff}), $\text{II}_n$ in (\ref{IInasymptotic}) and $\text{III}_n$ in (\ref{gdiffdetrate}) we can show
\begin{align*}
\int_{\wh U} g(x) \mathscr{H}(x) - \int_{U} g(x) \mathscr{H}(x) 
&= \int_U w_g(x) [f(x) - \wh f(x)]  d\mathscr{H}(x)  \\
&\hspace{1cm}+ O_p\left(\frac{1}{nh^{d+2}} \right) +o_p\left(h^{\nu-1}(\gamma_{n,h}^{(0)} + h^{\nu-1})\right).
\end{align*}
Recall that $\{U_\alpha:\alpha\in\mathscr{A}\}$ is a finite cover of $\mathcal{M}$. Immediately, we have 
%
%
\begin{align}\label{lambdadiffasympt}
%
\lambda(\wh f, g) - \lambda(f,g) =&  \int_\mathcal{M} w_g(x) [f(x) - \wh f(x)]  d\mathscr{H}(x)  + O_p\left(\frac{1}{nh^{d+2}} \right) +o_p\left(h^{\nu-1}(\gamma_{n,h}^{(0)} + h^{\nu-1})\right).
\end{align}

Following standard calculation for the bias of kernel density estimation (see page 83, Chac\'{o}n and Duong, 2018), similar to (\ref{expectationsq}) we can show that $\int_\mathcal{M} w(x) [f(x) - \mathbb{E} \wh f(x)] d\mathscr{H}(x) = h^\nu \mu (1+o(1))$, where 
\begin{align}\label{asymmean}
\mu = -(\nu!)^{-1} \int_\mathcal{M} \left[\int_{\mathbb{R}^d} y^{\otimes\nu} K(y)dy\right]^T \nabla^{\otimes\nu} f(x) w(x) \mathscr{H}(x).
\end{align}
Also by Theorem~\ref{asympnorm} we have that $\int_\mathcal{M} w_g(x) [\mathbb{E}\wh f(x) - \wh f(x)]  d\mathscr{H}(x)\rightarrow_D N(0,\sigma^2)$, where $\sigma^2=c\lambda(f,w_g^2R_K).$\\

The assertions (i) and (ii) follow from (\ref{lambdadiffasympt}) and by noticing that $1/(nh^{d+2}) = o(1/\sqrt{nh^d})$ when $nh^{2d+3}\rightarrow\infty$. Assertion (iii) is a consequence of $|\wh\sigma_{\tau_n}^2 - \sigma^2|=o_p(1)$. Below we briefly describe how to show this. Similar to (\ref{lambdadiffdecomp}), we can write $\wh\sigma_{\tau_n}^2 - \sigma^2=D_{n,1}+D_{n,2}+D_{n,3}$, where $D_{n,1}=cR(K)[\wt\lambda_{\tau_n}(\wh f,w_g) -  \lambda(f,w_g)]$, $D_{n,2}=cR(K)[\wt\lambda_{\tau_n}(\wh f,\wh w_g-w_g)-\lambda( f,\wh w_g-w_g)]$, and $D_{n,3}=cR(K)\lambda( f,\wh w_g-w_g)$. Immediately we have $D_{n,1}=o_p(1)$ using assertion (i), and $D_{n,3}=o_p(1)$ by notice that $\sup_{x\in\mathcal{I}(2\delta_0)}|\wh w_g(x)-w_g(x)|=O_p(\gamma_{n,h}^{(2)})+o_p(h^{\nu-2})$. We also have $D_{n,2}=o_p(1)$ under assumption (H) and $\tau_n=O(h)$ by applying Corollary~\ref{Hausdorff}. 

%
\hfill$\square$ 
\end{proof}

{\bf Proof of Theorem~\ref{lambdastardaggertheorem}}
\begin{proof}
Below we discuss the cases of $\wt\lambda_{\epsilon_n}=\lambda_{\epsilon_n}^*$ and $\wt\lambda_{\epsilon_n}=\lambda_{\epsilon_n}^\dagger$, respectively.\\

{\em Case of $\wt\lambda_{\epsilon_n}=\lambda_{\epsilon_n}^*$.} Many arguments are similar to the proofs of Lemma~\ref{pointlevelset} and Theorem~\ref{HausdorffIntegration} so we only give main steps. We would like to find the Taylor expansion of $\lambda_{c+{\epsilon}}(\wh f,g)=\int_{\wh f^{-1}(c+\epsilon)}g(x)d\mathscr{H}(x)$ as a function of $\epsilon$ when $\epsilon\rightarrow0$. Similar to (\ref{projection}), define $\wh\zeta_x(t) = x + t\wh\n(x)$. For any $x\in \wh{\mathcal{M}}_{c}$, let 
$$\wh t_\epsilon(x) = \argmin_t \{ |t|:\; \wh\zeta_x(t)  \in \wh{\mathcal{M}}_{c+\epsilon} \}.$$
Define $P_\epsilon(x)=\wh\zeta_x(\wh t_\epsilon(x))=x+ \wh t_\epsilon(x)\wh\n(x)$. Using Taylor expansion we have
\begin{align*}
c+\epsilon = \wh f(P_{\epsilon}(x)) = c + \wh t_{\epsilon}(x)\wh\n(x)^T \nabla \wh f(x) + \delta_\epsilon(x),
\end{align*}
where $\delta_\epsilon(x)=\frac{1}{2}\wh t_\epsilon(x)^2 \wh\n(x)^T\nabla^2 \wh f(\wh\zeta_x(s_1)) \wh\n(x)$ for some $s_1$ between 0 and $\wh t_{\epsilon}(x)$. Therefore $\wh t_{\epsilon}(x)=\|\wh f(x)\|^{-1}(\epsilon-\delta_\epsilon(x)),$ and 
\begin{align}\label{tepsilonrate}
\sup_{x\in \wh f^{-1}(c)} |\wh t_{\epsilon}(x)| = O_p(\epsilon).
\end{align}
%
%

We also can write 
%
$\epsilon = \wh f(P_{\epsilon}(x)) - \wh f(x) = \wh t_{\epsilon}(x)\wh\n(x)^T \wh r_{\epsilon}(x),$ 
where $\wh r_{\epsilon}(x) = \int_0^1 \nabla \wh f(\wh\zeta_x(s\wh t_{\epsilon}(x))) ds.$ Therefore $0 = \nabla  \wh t_{\epsilon}(x) [\wh\n(x)^T \wh r_{\epsilon}(x)] +  \wh t_{\epsilon}(x)\nabla [\wh\n(x)^T \wh r_{\epsilon}(x)].$ Here $\nabla [\wh\n(x)^T r_{\epsilon}(x)] = \nabla \wh\n(x)^T r_{\epsilon}(x) + \nabla r_{\epsilon}(x)^T \wh\n(x), $ where
\begin{align*}
 \nabla \wh r_{\epsilon}(x) = \int_0^1 \nabla^2 \wh f(\wh\zeta_x(s\wh t_{\epsilon}(x)))[\mathbf{I}_d+ s\nabla(\wh t_{\epsilon}(x)\wh\n(x))]ds.
\end{align*}
Note that $\wh\n(x)^T \wh r_{\epsilon}(x)-\|\nabla \wh f(x)\| =o_p(1)$. Using (\ref{tepsilonrate}) we have 
\begin{align}\label{tepsilonrate2}
\sup_{x\in \wh f^{-1}(c)} \|\nabla \wh t_{\epsilon}(x)\| = \sup_{x\in \wh f^{-1}(c)} \left\|\frac{-\wh t_\epsilon(x)\nabla [\wh\n(x)^T \wh r_{\epsilon}(x)]}{[\wh\n(x)^T \wh r_{\epsilon}(x)] }\right\| =O_p(\epsilon).
\end{align}
%
%
%
%

Since with probability one the reach of $\wh{\mathcal{M}}_c$ is positive for large $n$ by Lemma~\ref{compatible}, the map $\wh P_{\epsilon}$ is a diffeomorphism between $\wh{\mathcal{M}}_c$ and $\wh{\mathcal{M}}_{c+\epsilon}$ when $\epsilon$ is small. The Jacobian of $P_{\epsilon}$ is given by $A_{\epsilon}(x)=\mathbf{I}_d + R_{\epsilon}(x)$, where $R_{\epsilon}(x) = \nabla \wh\n(x)\wh t_{\epsilon}(x) + \wh\n(x)\nabla \wh t_{\epsilon}(x)^T$. By (\ref{tepsilonrate}) and (\ref{tepsilonrate2}) we have $\sup_{x\in \wh f^{-1}(c)} \|R_{\epsilon}(x)\| =O_p(\epsilon)$.\\

Suppose $\{(U_\alpha,\psi_\alpha):\; \alpha\in \mathscr{A}\}$ is a finite atlas for $\wh f^{-1}(c)$. Let $\mathscr{U}=\{U_\alpha:\;  \alpha\in \mathscr{A}\}$. Let $\{U,\psi\}$ be a chart on $\wh{\mathcal{M}}$ (see the beginning of the proof of Theorem~\ref{HausdorffIntegration}) such that $U=\psi(\Omega)$ for a subset $\Omega\subset \mathbb{R}^{d-1}$. Let the Jacobian of $\psi$ be $B$. Let $J=B^TB$, $D=BJ^{-1}B^T$, and
$$Q_\epsilon(x) = - [J(\psi^{-1}(x))]^{-1}B(\psi^{-1}(x))^T\left\{  2\nabla\wh\n(x) t_\epsilon(x) +  R_\epsilon(x)^T R_\epsilon(x)\right\}B(\psi^{-1}(x)).$$
Then following a similar argument as in the proof of Theorem~\ref{HausdorffIntegration}, we have
\begin{align*}
 \int_{P_\epsilon(U)}  g(x) d\mathscr{H}(x) -  \int_{U}  g(x) d\mathscr{H}(x) = \text{I}_\epsilon +  \text{II}_\epsilon +  \text{III}_\epsilon,
\end{align*}
where 
\begin{align*}
\text{I}_\epsilon & =  \int_{U}  [g(P_\epsilon(x))-g(x)] d\mathscr{H}(x), \\
\text{II}_\epsilon & = \int_{U} g(x)\left\{  \left\{ \det[\mathbf{I}_{d-1} - Q_\epsilon(x)] \right\}^{1/2} -  1 \right\}d\mathscr{H}(x),\\
%
\text{III}_\epsilon  & = \int_{U}  [g(P_\epsilon(x))-g(x)]  \left\{  \left\{ \det[\mathbf{I}_{d-1} - Q_\epsilon(x)] \right\}^{1/2} -  1 \right\}d\mathscr{H}(x).
\end{align*}

Notice that $\text{I}_\epsilon = \epsilon\int_{U} \|\nabla \wh f(x)\|^{-1} \wh\n(x)^T\nabla g(x)\mathscr{H}(x) + r_1(U,\epsilon),$ where 
\begin{align*}
r_1(U,\epsilon) = & -\int_U  \frac{\delta_{\epsilon}(x)}{\|\nabla\wh f(x)\|}  \wh\n(x)^T\nabla g\left(x \right) d\mathscr{H}(x)\\ 
& + \int_U \frac{1}{2}\wh t_\epsilon(x)^2\wh\n(x)^T\nabla^2 g\left(x + s_2 \wh\n(x)\wh t_\epsilon(x)\right) \wh\n(x) d\mathscr{H}(x),
\end{align*}
for some $0<s_2<1$. Hence $\sup_{U\in \mathscr{U}} |r_1(U,\epsilon)| = O_p(\epsilon^2)$ using (\ref{tepsilonrate}).\\ 

Next we study $\text{II}_\epsilon$. Let $Q_\epsilon(x) = Q_{\epsilon,1}(x) + Q_{\epsilon,2}(x)$ with 
\begin{align*}
&Q_{\epsilon,1}(x) = 2[J(\psi^{-1}(x))]^{-1}B(\psi^{-1}(x))^T  \nabla\wh \n(x) B(\psi^{-1}(x)) t_\epsilon(x),\\
&Q_{\epsilon,2}(x) =   [J(\psi^{-1}(x))]^{-1}B(\psi^{-1}(x))^T R_\epsilon(x)^T R_\epsilon(x) B(\psi^{-1}(x)).
\end{align*}
Let $\eta_\epsilon(x)=\eta_{\epsilon,1}(x)+\eta_{\epsilon,2}(x)$, where $\eta_{\epsilon,1}=\mathbf{tr}[D(\psi^{-1}(x)) R_\epsilon(x)^T R_\epsilon(x)]$ and $\eta_{\epsilon,2}(x) = \det(\mathbf{I}_{d-1}-Q_\epsilon(x)) - [1-\mathbf{tr}(Q_\epsilon(x))]$. Also let
\begin{align*}
\Delta_\epsilon(x) = &\frac{1}{2}  \eta_\epsilon(x)  + \sum_{i=2}^\infty (-1)^i\frac{(2i-3)!!}{2^ii!}  \left[ 2{\bf tr} \left[D(\psi^{-1}(x))  \nabla\wh\n(x) \right] t_\epsilon(x)  + \eta_\epsilon(x) \right]^i.
%
\end{align*}
Then $\text{II}_\epsilon =  \epsilon\int_{U} \|\nabla \wh f(x)\|^{-1} \mathbf{tr}[\wh S(x)] g(x)\mathscr{H}(x) + r_2(U,\epsilon),$ where $r_2(U,\epsilon) = \int_U \Delta_\epsilon(x) g(x)d\mathscr{H}(x).$ We can show that $\sup_{U\in\mathscr{U}} |r_2(U,\epsilon)| = O_p(\epsilon^2).$ 
%
%
%
%
Furthermore we can show that $\text{III}_\epsilon  =  r_3(U,\epsilon)$, where $\sup_{U\in\mathscr{U}} |r_3(U,\epsilon)|=O_p(\epsilon^2).$\\ 

Overall with $\wh w_g(x) = \|\nabla \wh f(x)\|^{-1}[ \wh \n(x)^T \nabla g(x) +  (d-1)\wh H(x)g(x)]$, we have
\begin{align*}
 \int_{P_\epsilon(U)}  g(x) d\mathscr{H}(x) -  \int_{U}  g(x) d\mathscr{H}(x) = \epsilon \int_U \wh w_g(x) d\mathscr{H}(x) + [r_1(U,\epsilon)+r_2(U,\epsilon)+r_3(U,\epsilon)].
\end{align*}

Therefore $\int_{\wh{\mathcal{M}}_{c+\epsilon}}  g(x) d\mathscr{H}(x) -  \int_{\wh{\mathcal{M}}_c}  g(x) d\mathscr{H}(x) = \epsilon \int_{\wh{\mathcal{M}}_c} \wh w_g(x) d\mathscr{H}(x) + r(\epsilon),$ where $|r(\epsilon)| =O_p(\epsilon^2).$ 
%
It then follows from (\ref{altstar}) that
\begin{align*}
|\lambda_{\epsilon_n}^*(\wh f,g) - \lambda(\wh f,g)|=&\left| \frac{1}{2\epsilon_n} \int_{-\epsilon_n}^{\epsilon_n}  \int_{\wh{\mathcal{M}}_{c+\epsilon}}  g(x) d\mathscr{H}(x) d\epsilon -  \int_{\wh{\mathcal{M}}_c}  g(x) d\mathscr{H}(x) \right| \\
 = & \left|  \frac{1}{2\epsilon_n}  \int_{-\epsilon_n}^{\epsilon_n}  \epsilon d\epsilon \int_{\wh{\mathcal{M}}_c}w(x) d\mathscr{H}(x)  +  \frac{1}{2\epsilon_n} \int_{-\epsilon_n}^{\epsilon_n}  r(\epsilon) d\epsilon \right|\\
 =&\left| \frac{1}{2\epsilon_n} \int_{-\epsilon_n}^{\epsilon_n}  r(\epsilon) d\epsilon\right| = O_p(\epsilon_n^2).
 %
\end{align*}
%

{\em Case of $\wh\lambda=\lambda_{\epsilon_n}^\dagger$.} The set $\wh{\mathcal{M}}_c\oplus\epsilon_n$ is tube of radius $\epsilon_n$ around the submanifold $\wh{\mathcal{M}}_c$. Then $\wh\zeta(x,\epsilon):=\wh\zeta_x(\epsilon)$ defines a diffeomorphism between $\wh{\mathcal{M}}_c\times(-\epsilon_n,\epsilon_n)$ and $\wh{\mathcal{M}}_c\oplus \epsilon_n$.
The Weyl's volume element in the tube given in Wely (1939) (also see Gray, 2004) is $\text{det}(\mathbf{I}_d+\epsilon \wh S(x))  d\epsilon d\mathscr{H}(x),$ where $\wh S(x)$ is the shape operator on $\wh{\mathcal{M}}$ as defined in (\ref{shapeoperator}). Therefore
\begin{align*}
 \frac{1}{2\epsilon_n} \int_{\wh{\mathcal{M}}_c\oplus \epsilon_n} g(x) dx = \frac{1}{2\epsilon_n}  \int_{\wh{\mathcal{M}}_c} \int_{-\epsilon_n}^{\epsilon_n} g(x+\epsilon\wh\n(x)) \text{det}(\mathbf{I}_d+\epsilon \wh S(x)) d\epsilon d\mathscr{H}(x).
\end{align*}

Note that
\begin{align*}
g(x+\epsilon\wh\n(x)) =g(x) + \epsilon\wh\n(x)^T\nabla g(x) + \frac{1}{2}\epsilon^2 \wh\n(x)^T \nabla^2g(x+s\epsilon\wh\n(x))\wh\n(x),
\end{align*}
for some $0<s<1$, and $\text{det}(\mathbf{I}_d+\epsilon\wh S(x)) = 1+\epsilon \mathbf{tr}[\wh S(x)] + \epsilon^2 r(x,\epsilon),$ 
where $|r(x,\epsilon)|=O_p(1)$ using arguments similar to (\ref{Delta2ratepart1}) and (\ref{Delta2ratepart2}).\\

Let $R(x,\epsilon) =  \frac{1}{2} \wh\n(x)^T \nabla^2g(x+s\epsilon\wh\n(x))\wh\n(x) + g(x) r(x,\epsilon) + [\wh\n(x)^T\nabla g(x) + \frac{1}{2}\epsilon \wh\n(x)^T \nabla^2g(x+s\epsilon\wh\n(x))\wh\n(x)] \{\mathbf{tr}[\wh S(x)]+ \epsilon r(x,\epsilon)\}. $ 
Then it follows that $|R(x,\epsilon)|=O_p(1)$ and hence 
\begin{align*}
&|\lambda_{\epsilon_n}^\dagger(\wh f,g) - \lambda(\wh f,g)|\\
= & \left| \frac{1}{2\epsilon_n} \int_{\wh{\mathcal{M}}_c\oplus \epsilon_n} g(x) dx - \int_{\wh{\mathcal{M}}_c} g(x) d\mathscr{H}(x) \right| \\
 = &  \left|  \frac{1}{2\epsilon_n}  \int_{-\epsilon_n}^{\epsilon_n} \epsilon d\epsilon  \int_{\wh{\mathcal{M}}_c} \{ \wh\n(x)^T\nabla g(x) + g(x)\mathbf{tr}[\wh S(x)] \} d\mathscr{H}(x) +  \frac{1}{2\epsilon_n}  \int_{-\epsilon_n}^{\epsilon_n} \int_{\wh{\mathcal{M}}_c} \epsilon^2 R(x,\epsilon) d\mathscr{H}(x) d\epsilon \right| \\
 =& \left|  \frac{1}{2\epsilon_n} \int_{\wh{\mathcal{M}}_c}  \int_{-\epsilon_n}^{\epsilon_n} \epsilon^2 R(x,\epsilon) d\mathscr{H}(x)  d\epsilon \right| = O_p(\epsilon_n^2). 
 %
\end{align*}
$\hfill\square$
\end{proof}
%
%

{\bf Proof of Corollary \ref{Hausdorff}} 
\begin{proof}
First we consider the case of $\wt\lambda_{\epsilon_n}=\lambda$. Following the proof of Theorem~\ref{HausdorffIntegration}, 
we continue to use $\text{I}_n$, $\text{II}_n$ and $\text{III}_n$ as in (\ref{HausInexpr}), (\ref{HausIInexpr}) and (\ref{HausIIInexpr}), where we replace $g$ by $b_n$. Then $\int_{\wh U} b_n(x) d\mathscr{H}(x) - \int_{U} b_n(x) d\mathscr{H}(x) =\text{I}_n + \text{II}_n + \text{III}_n$, where 
\begin{align*}
\text{I}_n & =  \int_{U}  [b_n(P_n(x))-b_n(x)] d\mathscr{H}(x), \\
\text{II}_n & = \int_{U} b_n(x)\left\{  \left\{ \det[\mathbf{I}_{d-1} - Q_n(x)] \right\}^{1/2} -  1 \right\}d\mathscr{H}(x),\\
%
\text{III}_n  & = \int_{U}  [b_n(P_n(x))-b_n(x)]  \left\{  \left\{ \det[\mathbf{I}_{d-1} - Q_n(x)] \right\}^{1/2} -  1 \right\}d\mathscr{H}(x),
\end{align*}

where $Q_n$ is given (\ref{qndef}). By using Taylor expansion and Lemma~\ref{pointlevelset} we have
\begin{align}\label{gdiffsup}
&\sup_{x\in U}|b_n(P_n(x)) - b_n(x)|
\leq  \sup_{x\in U} |t_n(x)| \sup_{x\in \mathcal{I}(\delta_0)} \|\nabla b_n(x)\| = O_p\left\{\beta_n \left(\gamma_{n,h}^{(0)} + h^\nu\right)\right\}.
\end{align}
For $\text{I}_n$, we obtain from (\ref{gdiffsup}) that
\begin{align}\label{Iexpression}
\text{I}_n \leq   \int_{U} \left| b_n(P_n(x)) - b_n(x) \right| d\mathscr{H}(x) 
\leq  \sup_{x\in U}|b_n(P_n(x)) - b_n(x)| \mathscr{H}(U) 
=  O_p\left\{\beta_n \left(\gamma_{n,h}^{(0)} + h^\nu\right)\right\}.
\end{align}

For $\text{II}_n$, (\ref{detdiffrate}) leads to
\begin{align}\label{IIexp}
\text{II}_n \leq & \int_U |b_n(x)|   \left|  \{\det\left[\mathbf{I}_{d-1} - Q_n(x)\right]\}^{1/2} - 1\right| d\mathscr{H}(x) \nonumber\\
\leq & \sup_{x\in U} |b_n(x)|  \sup_{x\in U} \left|  \{\det\left[\mathbf{I}_{d-1} - Q_n(x)\right]\}^{1/2} - 1\right| \mathscr{H}(U) \nonumber\\
%
%
= & O_p\left\{\alpha_n \left(\gamma_{n,h}^{(0)} + h^\nu\right)\right\} .
\end{align}

For $\text{III}_n$, by using (\ref{gdiffsup}) and (\ref{detdiffrate}) we have
\begin{align}\label{IIIexp}
\text{III}_n = O_p\left\{\beta_n \left(\gamma_{n,h}^{(0)} + h^\nu\right)^2 \right\}.
\end{align}
By having (\ref{Iexpression}), (\ref{IIexp}) and (\ref{IIIexp}) we get (\ref{Dgn0rate2}) for $\wt\lambda_{\epsilon_n}=\lambda$. 
Next we consider the cases of $\wt\lambda_{\epsilon_n}=\lambda_{\epsilon_n}^*$  or $\lambda_{\epsilon_n}^\dagger$. Notice that 
\begin{align}\label{lambdadaggerdecomp}
\wt\lambda_{\epsilon_n} (\wh f, b_n) - \lambda (f, b_n) = [\lambda (\wh f, b_n) - \lambda (f, b_n)] + [\wt\lambda_{\epsilon_n} (\wh f, b_n) - \lambda (\wh f, b_n)], 
\end{align}
where $\lambda (\wh f, b_n) - \lambda (f, b_n)$ has been studied above in the case of $\wt\lambda_{\epsilon_n}=\lambda$. The proof is completed once we otain $ |\wt\lambda_{\epsilon_n} (f, b_n) - \lambda (f, b_n)|  = O_p\left(\epsilon_n^2 \left( \alpha_n + \beta_n + \eta_n\right)  \right)$,
which can be shown by replacing $g$ by $b_n$ in the proof of Theorem~\ref{lambdastardaggertheorem}.$\hfill\square$

%
%

\end{proof}
%
%
%
%
%

{\bf Proof of Theorem \ref{asympototicnormalitymknown} } 
\begin{proof}
%
In the proof we use $C>0$ as a generic constant that can vary as we need. Let $Q_\nu=\nu\vee(Q+1).$ Using Taylor expansion, we have 
$\wh g(x) - g(x) = \sum_{j=1}^{Q_\nu} T_{n,j}(x),$
where 
\begin{align}\label{Tnjdef}
T_{n,j}(x) = 
\begin{cases}
\frac{1}{j!} [d_{\wh f}(x)^T -d_f(x)^T ]^{\otimes j} \; \nabla^{\otimes j} \phi(d_f(x)) ,&  j=1,\cdots,Q_\nu-1,\\
\frac{1}{j!} [d_{\wh f,1}(x)^T -d_{f,1}(x)^T ]^{\otimes j} \; \nabla^{\otimes j} \phi(\wt{d_{f,1}}(x)), & j=Q_\nu,
\end{cases}
\end{align}
where $\wt{d_{f,1}}$ is between $d_{f,1}$ and $d_{\wh f,1}$. Using (\ref{supvar1}), (\ref{supbias1}) and assumptions (H) and (P), we have 
\begin{align}\label{qnv}
\lambda(f,T_{n,Q_\nu}) = O_p((\gamma_{n,h}^{(1)} + h^{\nu-1})^\nu) = O_p((h\gamma_{n,h}^{(2)} + h^{\nu-1})^\nu) = o_p(h^\nu).    
\end{align}
%
%
We can write $\lambda(f,T_{n,1}) - \mathbb{E} [\lambda(f,T_{n,1})] = \frac{1}{n}\sum_{i=1}^n (\xi_i - \mathbb{E} \xi_i)$, where $\xi_i= \xi_{i,1} + \xi_{i,2}$ with
\begin{align*}
& \xi_{i,l} = \frac{1}{h^{d+l}} \int_{\cal M}    \nabla_l \phi(d_f(x))^T d_{K,l}\left( \frac{x-X_i}{h} \right) d\mathscr{H}(x),\; l=1,2.
\end{align*}
Following the proof of Theorem~\ref{asympnorm}, one can show that for $l=1,2$,
$$\sqrt{nh^{1+2l}}\left[\frac{1}{n}\sum_{i=1}^n (\xi_{i,l} - \mathbb{E} \xi_{i,l}) \right] \rightarrow_D N(0,\sigma_l^2),$$
where $\sigma_l^2$ is given in (\ref{asyptoticvarexp}). Note that $\xi_{i,2}\equiv 0$ when $Q=0$. Therefore
\begin{align}\label{tn1asymp}
 \sqrt{nh^{1+2l_Q}} [\lambda(f,T_{n,1}) - \mathbb{E} [\lambda(f,T_{n,1})]] \rightarrow_D N(0,\sigma_{l_Q}^2).   
\end{align}

Following from standard calculation, one can show that 
\begin{align}
\mathbb{E} \lambda(f,T_{n,1}) =  h^{\nu} (\nu!)^{-1} \int_\mathcal{M} \left[\int_{\mathbb{R}^d} y^{\otimes\nu} K(y)dy\right]^T \nabla^{\otimes\nu} d_f(x)^T \nabla \phi(d_f(x))  \mathscr{H}(x)(1+o(1)).\label{asymptoticmu1exp}
%
\end{align}
For $j=2$, if $\nu=2$ and $Q\leq 1$, then $Q_\nu=2$ and (\ref{qnv}) applies, which implies $\mathbb{E}\lambda(f,T_{n,2})$ can be asymptotically absorbed into $\mathbb{E} \lambda(f,T_{n,1})$; otherwise by standard calculation we have
\begin{align}
\mathbb{E}\lambda(f,T_{n,2}) = \frac{1}{nh^{d+2l_{Q-1}}}\frac{c}{2}\int_{\mathcal{M}} \int_{\mathbb{R}^d} [d_{K,l_{Q-1}}(y)^T ]^{\otimes2}\nabla_{l_{Q-1}}^{\otimes2}\phi(d_f(x)) dyd\mathscr{H}(x)(1+o(1)) + O(h^{2\nu}),\label{asympototicmu2exp}  
\end{align}
which depends on $Q$ because the second partial derivatives of $f$ and $\wh f$ are not involved in $T_{n,2}$ if $Q\leq1$. Next we will show that $\lambda(f,T_{n,1}) + \mathbb{E}[\lambda(f,T_{n,2})]$ is the leading term in $\lambda(f,\sum_{j=1}^{Q_\nu} T_{n,j}) = \sum_{j=1}^{Q_\nu} \lambda(f,T_{n,j})$, i.e., their difference is relatively negligible. Note that $\lambda(f,T_{n,Q_\nu})$ has been discussed in (\ref{qnv}).\\

To simplify our discussion, in the following we assume that $Q\geq 3\vee (\nu-1)$ so that all $T_{n,j}$'s involve the second partial derivatives of $f$ and $\wh f$, for $j=1,\cdots,Q_\nu-1$. When $Q < 3\vee (\nu-1)$, some $T_{n,j}$'s may only include the first partial derivatives of $f$ and $\wh f$ and the calculation is very similar to (and simpler than) what we give below. For simplicity we write $\psi_j(x)$ as a generic element of the vector $\nabla^{\otimes j} \phi(d_f(x))$. In view of assumption in (\ref{nablajbound}), we have $\sup_{x\in\mathcal{M}} |\psi_j(x)| \leq C$.\\

For non-negative integers $p$ and $q$ with $1\leq p+q\leq Q_{\nu}-1$, we consider 
\begin{align}\label{gexpress}
H_{n,h}^{p,q} (x) := \prod_{i=1}^p \left[\wh f^{(\alpha_i)}(x) - f^{(\alpha_i)}(x) \right]\prod_{j=1}^q \left[\wh f^{(\beta_j)}(x) -  f^{(\beta_j)}(x) \right] \psi_{p+q}(x),
\end{align}
where $\alpha=\{\alpha_i\}_{i=1}^p$ and $\beta=\{\beta_j\}_{j=1}^q$, $\alpha_i$'s and $\beta_j$'s are $d$-dimensional indices with $|\alpha_i|=1$, $i=1,\dots,p$ and $|\beta_j|=2$, $j=1,\dots,q$. Here we use the convention that if $p=0$ (or $q=0$), the product of the first (or second) derivatives is not included. Note that $\lambda(f,T_{n,j})$ is a sum of $(a_d)^j$ functions in the generic form of $\frac{1}{j!}\lambda(f,H_{n,h}^{p,q})$ with $p+q=j$, for $j=1,\cdots,Q_\nu-1$. Recall that $a_d=d+d(d+1)/2$.\\

Next we study $\mathbb{E}\lambda(f,H_{n,h}^{p,q})$. Without loss of generality we assume $n\geq 4Q_\nu$. Let 
\begin{align*}
D_{h,1}(\alpha_i,X_{k_i},x) = \frac{1}{h^{d+1}}K^{(\alpha_i)}\left(\frac{x-X_{k_i}}{h}\right)  - f^{(\alpha_i)}(x),\\
D_{h,2}(\beta_j,X_{l_j},x) = \frac{1}{h^{d+2}}K^{(\beta_j)}\left(\frac{x-X_{l_j}}{h}\right)  - f^{(\beta_j)}(x).
\end{align*}

We write $\mathbf{k}=(k_1,\cdots,k_p)$, $\mathbf{l}=(l_1,\cdots,l_q)$, $\mathbf{k}^\prime=(k_1^\prime,\cdots,k_p^\prime)$, and $\mathbf{l}^\prime=(l_1^\prime,\cdots,l_q^\prime)$. We put braces around the notation when we want to emphasize it is a set instead of a sequence. For example, $\{\mathbf{k}\}=\{k_1,\cdots,k_p\}$.
Notice that we can write
\begin{align}\label{hnhexpress}
%
H_{n,h}^{p,q}(x) =\psi_{p+q}(x)\frac{1}{n^{p+q}} \sum_{\mathbf{k}\otimes\mathbf{l}\in\mathscr{V}_{p,q}} \prod_{i=1}^p  D_{h,1}(\alpha_i,X_{k_i},x)  \prod_{j=1}^q   D_{h,2}(\beta_j,X_{l_j},x) , 
\end{align}
where the index set $\mathscr{V}_{p,q}=\{\mathbf{k}\otimes\mathbf{l}: k_1,\cdots,k_p,l_1,\cdots,l_q \in\{1,\cdots,n\}\}.$ It has a partition$$\mathscr{V}_{p,q} = V(p,q,0) \bigcup \left( \mathop{\bigcup_{s=0}^p\bigcup_{t=0}^q}_{s+t\leq p+q-2}\bigcup_{u=1}^{\lfloor (p+q-s-t)/2\rfloor}V(s,t,u) \right),$$ where $V(s,t,u)$ is a subset of $\mathscr{V}_{p,q}$ such that there are $s$ of $\{\mathbf{k}\}$ and $t$ of $\{\mathbf{l}\}$ that have no duplicates in $\{\mathbf{k}\}\cup \{\mathbf{l}\}$, among which $u$ indices have duplicates. In other words, for any $\mathbf{k}\otimes\mathbf{l}\in V(s,t,u)$, $\{\mathbf{k}\}\cup \{\mathbf{l}\}$ has $s+t+u$ unique indices. Let $|A|$ be the cardinality of a set $A$. In general, we have $|V(s,t,u)| =O(n^{s+t+u})$ as $n\rightarrow\infty$. 
%
Then from (\ref{hnhexpress}) we can write
\begin{align}\label{Hnhpqexpress}
%
H_{n,h}^{p,q}(x) =\frac{1}{n^{p+q}}  \left[ \mathop{\sum_{s=0}^p \sum_{t=0}^q}_{s+t\leq p+q-2} \sum_{u=1}^{\lfloor (p+q-s-t)/2 \rfloor}\text{I}_{h}^{s,t,u}(x) + I_{h}^{p,q,0}(x)\right],
\end{align}
where $\text{I}_{h}^{s,t,u}(x) = \psi_{p+q}(x)\sum_{\mathbf{k}\otimes\mathbf{l}\in V(s,t,u)}  \prod_{i=1}^p D_{h,1}(\alpha_i,X_{k_i},x)  \prod_{j=1}^q D_{h,2}(\beta_j,X_{l_j},x) .$
%
%
We have
\begin{align}\label{generalform}
&\mathbb{E} \lambda(f,I_h^{s,t,u})\nonumber\\
=&\sum_{\mathbf{k}\otimes\mathbf{l}\in V(s,t,u)}\mathbb{E} \left[\int_{\mathcal{M}} \psi_{p+q}(x)\prod_{i=1}^p D_{h,1}(\alpha_i,X_{k_i},x)  \prod_{j=1}^q D_{h,2}(\beta_j,X_{l_j},x)  d\mathscr{H}(x) \right]\nonumber\\
=&\sum_{\mathbf{k}\otimes\mathbf{l}\in V(s,t,u)} \int_{\mathcal{M}}  \psi_{p+q}(x)\int_{\mathbb{R}^{d(s+t+u)}}  \prod_{i=1}^p D_{h,1}(\alpha_i,y_{k_i},x)  \prod_{j=1}^q D_{h,2}(\beta_j,y_{l_j},x) \nonumber \\
&\hspace{7cm}\times\prod_{\iota\in\text{uniq}(\{\mathbf{k}\}\cup\{\mathbf{l}\})} f(y_\iota) d\mathbf{y}d\mathscr{H}(x),
%
\end{align}
where ``uniq'' means the unique indices in the set, and $\mathbf{y}=(y_\iota)_{ \iota\in\text{uniq}(\{\mathbf{k}\}\cup\{\mathbf{l}\})}$.\\

If $u=0$ (i.e., $s=p$ and $t=q$), we have $|V(p,q,0)|= \frac{n!}{(n-p-q)!}$ and by applying change of variables for $p+q$ times and integration by parts in (\ref{generalform}) we have 
\begin{align*}
|\mathbb{E} \lambda(f, I_h^{p,q,0}) | \leq C n^{p+q}h^{(p+q)\nu},
%
%
%
\end{align*}
for some constant $C>0$. 
If $u\geq 1$, recall that $|V(s,t,u)|=O(n^{s+t+u})$. We apply change of variable $z_\iota=(x-y_\iota)/h$ for the unique $\iota$'s ($s+t+u$ of them) in (\ref{generalform}) and use integration by parts for the indices with no duplicates ($s+t$ of them) and get
\begin{align*}
|\mathbb{E} \lambda(f,I_h^{p,q,u})| \leq  C n^{s+t+u}h^{(s+t)\nu} h^{(s+t+u)d} \left(\frac{1}{h^{d+1}} \right)^{p-s} \left(\frac{1}{h^{d+2}} \right)^{q-t}.
%
\end{align*}

So using (\ref{Hnhpqexpress}) we get
\begin{align}
& |\mathbb{E}\lambda(f,H_{n,h}^{p,q})| \nonumber \\
 \leq &C \left( h^{(p+q)\nu}  + \frac{1}{n^{p+q}}  \mathop{\sum_{s=0}^p \sum_{t=0}^q}_{s+t\leq p+q-2} \sum_{u=1}^{\lfloor (p+q-s-t)/2 \rfloor} n^{s+t+u} h^{(s+t)\nu + ud- (p-s)(d+1) - (q-t)(d+2)}\right)  \nonumber \\
 = & C \left( h^{(p+q)\nu}  + \mathop{\sum_{s=0}^p \sum_{t=0}^q}_{s+t\leq p+q-2} \sum_{u=1}^{\lfloor (p+q-s-t)/2 \rfloor}  \frac{(nh^d)^u \label{sumofrate} h^{(s+t)\nu}}{(nh^{d+1})^{p-s}(nh^{d+2})^{q-t}} \right) : = C \rho_{n,h,1}^{p,q}.
\end{align}

%
%
%

We would like to simplify the upper bound for $|\mathbb{E}\lambda(f,H_{n,h}^{p,q})|$ under three possible scenarios a)-c) below.\\

a) When $nh^{d+4+2\nu} \rightarrow\gamma_1$ for $1<\gamma_1\leq +\infty$. Let $r_{n,h,1} = h^{(p+q)\nu}$. Then
\begin{align*}
&\frac{\rho_{n,h,1}^{p,q}}{r_{n,h,1} }  \\
= & 1+ \mathop{\sum_{s=0}^p \sum_{t=0}^q}_{s+t\leq p+q-2} \sum_{u=1}^{\lfloor (p+q-s-t)/2 \rfloor}  \frac{(nh^d)^u h^{(s+t-p-q)\nu}}{(nh^{d+1})^{p-s}(nh^{d+2})^{q-t}} \\  %
=& 1 + \mathop{\sum_{s=0}^p \sum_{t=0}^q}_{s+t\leq p+q-2} \sum_{u=1}^{\lfloor (p+q-s-t)/2 \rfloor}  \frac{1}{(nh^{d+2+2\nu})^{(p-s)/2}(nh^{d+4+2\nu})^{(q-t)/2} (nh^d)^{(p-s)/2+(q-t)/2-u}} \\
%
\leq & 1+ \left[\sum_{s=0}^p \frac{1}{(nh^{d+2+2\nu})^{(p-s)/2}}\right] \left[\sum_{t=0}^q  \frac{1}{(nh^{d+4+2\nu})^{(q-t)/2}} \right] \left[\sum_{u=1}^{\lfloor (p+q)/2 \rfloor} \frac{1}{(nh^d)^{p/2+q/2-u}}\right],
\end{align*}
which is bounded for $1\leq p+q \leq Q_\nu-1$.\\

b) When $nh^{d+4+2\nu} \rightarrow\gamma_1$ and $nh^{d+2+2\nu} \rightarrow \gamma_2$ for $0\leq\gamma_1\leq 1$ and $1<\gamma_2\leq +\infty$. Let $r_{n,h,2} = \frac{h^{p\nu}}{(nh^{d+4})^{q/2}}$. 
Then
\begin{align*}
&\frac{\rho_{n,h,1}^{p,q}}{2^{p+q} r_{n,h,2} }  \\
\leq & (nh^{d+4+2\nu})^{q/2} h^{q\nu} + \mathop{\sum_{s=0}^p \sum_{t=0}^q}_{s+t\leq p+q-2} \sum_{u=1}^{\lfloor (p+q-s-t)/2 \rfloor}  \frac{ (0.5nh^{d+4+2\nu})^{t/2}}{(nh^d)^{(p-s)/2+ (q-t)/2-u}(nh^{d+2+2\nu})^{(p-s)/2}} \\  %
\leq & (nh^{d+4+4\nu})^{q/2} + \left[\sum_{s=0}^p \frac{1}{(nh^{d+2+2\nu})^{(p-s)/2}}\right] \left[\sum_{t=0}^q  (0.5nh^{d+4+2\nu})^{t/2} \right] \left[\sum_{u=1}^{\lfloor (p+q)/2 \rfloor} \frac{1}{(nh^d)^{p/2+q/2-u}}\right],
%
\end{align*}
which is bounded for $1\leq p+q \leq Q_\nu-1$.\\

c) When $nh^{d+2+2\nu} \rightarrow \gamma_2$ for $0\leq\gamma_2\leq 1$. Let $r_{n,h,3} = \frac{ 1}{(nh^{d+2})^{p/2}(nh^{d+4})^{q/2}}$. Then
\begin{align*}
&\frac{\rho_{n,h,1}^{p,q}}{2^{p+q} r_{n,h,3} }  \\
\leq & (nh^{d+2+2\nu})^{p/2} (nh^{d+4+2\nu})^{q/2} + \mathop{\sum_{s=0}^p \sum_{t=0}^q}_{s+t\leq p+q-2} \sum_{u=1}^{\lfloor (p+q-s-t)/2 \rfloor} \frac{ (0.5nh^{d+2+2\nu})^{s/2} (nh^{d+4+2\nu})^{t/2}}{(nh^d)^{ (p-s)/2+(q-t)/2-u}}. \\
\leq & (nh^{d+2+2\nu})^{p/2} (nh^{d+4+2\nu})^{q/2} + \left[\sum_{s=0}^p (0.5nh^{d+2+2\nu})^{s/2}\right] \left[\sum_{t=0}^q  (nh^{d+4+2\nu})^{t/2} \right] \left[\sum_{u=1}^{\lfloor (p+q)/2 \rfloor} \frac{1}{(nh^d)^{p/2+q/2-u}}\right],
\end{align*}
which is bounded for $1\leq p+q \leq Q_\nu-1$.\\

Let $\beta_{n,h} = \max\{h^\nu, (nh^{d+4})^{-1/2}\}.$ In view of the discussion in a), b) and c), overall we can write $|\mathbb{E}\lambda(f,H_{n,h}^{p,q})| \leq C \beta_{n,h}^{p+q} $ for some constant $C$. \\

%
%
Therefore
\begin{align}\label{sumelambda}
\left|\sum_{j=3}^{Q_\nu-1}\mathbb{E}[\lambda(f,T_{n,j})]\right| \leq     \sum_{j=3}^{Q_\nu-1}\left|\mathbb{E}[\lambda(f,T_{n,j})]\right| \leq \sum_{j=3}^{Q_\nu-1} \frac{1}{j!} C(a_d)^j \beta_{n,h}^j = O(\beta_{n,h}^3).
\end{align}

Combining the above result with (\ref{asymptoticmu1exp}) and (\ref{asympototicmu2exp}), it follows that 
\begin{align}\label{biasexpression}
\mathbb{E} [\lambda (f,\wh g-g)] = \left[h^{\nu}\mu_1 + \frac{\mu_2}{nh^{d+4}}\right](1+o(1)).    
\end{align}

Next we consider the variance of $\lambda(f,\wh g)$. Notice that 
\begin{align*}
\text{Var}\left[\sum_{j=1}^{Q_\nu-1}\lambda(f,T_{n,j})\right]  =  \text{Var}\left[\lambda(f,T_{n,1})\right]  +\text{Var}\left[\sum_{j=2}^{Q_\nu-1}\lambda(f,T_{n,j})\right] + 2 \text{Cov} \left[\lambda(f,T_{n,1}), \sum_{j=2}^{Q_\nu-1}\lambda(f,T_{n,j})\right].
\end{align*}

In what follows we will show that 
\begin{align}\label{corremainder}
\text{Var}\left[\sum_{j=2}^{Q_\nu-1}\lambda(f,T_{n,j})\right] + 2 \text{Cov} \left[\lambda(f,T_{n,1}), \sum_{j=2}^{Q_\nu-1}\lambda(f,T_{n,j})\right]= o((nh^5)^{-1}).     
\end{align}
We consider the covariance between $\lambda(f,H_{n,h}^{p,q})$ and $\lambda(f,H_{n,h}^{p^\prime,q^\prime})$ for $1\leq p+q\leq Q_\nu-1$ and $2\leq p^\prime+q^\prime\leq Q_\nu-1$, which can be written as 
\begin{align}\label{varhexp}
\text{Cov}[\lambda(f,H_{n,h}^{p,q}),\lambda(f,H_{n,h}^{p^\prime,q^\prime})] &= \mathbb{E}\int_{\mathcal{M}}\int_{\mathcal{M}} H_{n,h}^{p,q}(x) H_{n,h}^{p^\prime,q^\prime}(x^\prime)d\mathscr{H}(x)d\mathscr{H}(x^\prime) \nonumber \\
&\hspace{2cm} - \mathbb{E}[\lambda(f,H_{n,h}^{p,q})]\mathbb{E}[\lambda(f,H_{n,h}^{p^\prime,q^\prime})].
\end{align}
The index set $\mathscr{V}_{p,q}\otimes \mathscr{V}_{p^\prime,q^\prime}$ has a partition$$\mathscr{V}_{p,q}\otimes \mathscr{V}_{p^\prime,q^\prime} = \wt V(*,*,*,*,0) \bigcup \left[ \mathop{\bigcup_{s=0}^{p}\bigcup_{t=0}^{q}}_{s+t\leq p+q-1} \mathop{\bigcup_{s^\prime=0}^{p^\prime}\bigcup_{t^\prime=0}^{q^\prime}}_{s^\prime+t^\prime\leq p^\prime+q^\prime-1} \bigcup_{u=1}^{\lfloor (p+p^\prime+q+q^\prime-s-t-s^\prime-t^\prime)/2 \rfloor}\wt V(s,t,s^\prime,t^\prime,u) \right],$$ 
where $\wt V(s,t,s^\prime,t^\prime,u)$ is a subset of $\mathscr{V}_{p,q}\otimes \mathscr{V}_{p^\prime,q^\prime}$ such that there are $s$ of $\{\mathbf{k}\}$, $t$ of $\{\mathbf{l}\}$, $s^\prime$ of $\{\mathbf{k}^\prime\}$ and $t^\prime$ of $\{\mathbf{l}^\prime\}$ that have no duplicates in $\{\mathbf{k}\}\cup \{\mathbf{l}\} \cup \{\mathbf{k}^\prime\}\cup \{\mathbf{l}^\prime\}$, among  which $u$ indices have duplicates. We also require that $\{\mathbf{k}\}\cup \{\mathbf{l}\}$ and $\{\mathbf{k}^\prime\}\cup \{\mathbf{l}^\prime\}$ share at least one duplicate. We use $V(*,*,*,*,0)$ to represent the set that no duplicate is shared between $\{\mathbf{k}\}\cup \{\mathbf{l}\}$ and $\{\mathbf{k}^\prime\}\cup \{\mathbf{l}^\prime\}$. 
Let 
\begin{align*}
\text{I}_h^{s,t,s^\prime,t^\prime,u}(x,x^\prime) = &\psi_{p+q}(x)\psi_{p+q}(x^\prime)\times\sum_{\mathbf{k}\otimes\mathbf{l}\otimes\mathbf{k}^\prime\otimes\mathbf{l}^\prime\in \wt V(s,t,s^\prime,t^\prime,u)}  \prod_{i=1}^p D_{h,1}(\alpha_i,X_{k_i},x)  \\
&\hspace{1cm} \times \prod_{j=1}^q D_{h,2}(\beta_j,X_{l_j},x)  \prod_{i=1}^{p^\prime}  D_{h,1}(\alpha_i,X_{k_i^\prime},x^\prime)   \prod_{j=1}^{q^\prime} D_{h,2}(\beta_j,X_{l_j^\prime},x^\prime) .
\end{align*}

Then we can write
\begin{align}\label{hhdecomp}
%
H_{n,h}^{p,q}(x) H_{n,h}^{p^\prime,q^\prime}(x^\prime) = \frac{1}{n^{p+p^\prime+q+q^\prime}} [\text{I}_h^{*,*,*,*,0}(x,x^\prime) + \text{II}_h(x,x^\prime)] ,
\end{align}
where
\begin{align*}
\text{II}_h(x,x^\prime)=  \mathop{\sum_{s=0}^p \sum_{t=0}^q}_{s+t\leq p+q-1}   \mathop{\sum_{s^\prime=0}^{p^\prime} \sum_{t^\prime=0}^{q^\prime}}_{s^\prime+t^\prime\leq p^\prime+q^\prime-1} \sum_{u=1}^{\lfloor (p+p^\prime+q+q^\prime-s-t-s^\prime-t^\prime)/2 \rfloor} \text{I}_h^{s,t,s^\prime,t^\prime,u}(x,x^\prime).
\end{align*}

Note that  
\begin{align}
& \frac{1}{n^{p+p^\prime+q+q^\prime}}\mathbb{E}\int_{\mathcal{M}}\int_{\mathcal{M}}  \text{I}_h^{*,*,*,*,0}(x,x^\prime)d\mathscr{H}(x)d\mathscr{H}(x^\prime) - \mathbb{E}[\lambda(f,H_{n,h}^{p,q})]\mathbb{E}[\lambda(f,H_{n,h}^{p^\prime,p^\prime})]\nonumber\\
 = &\frac{1}{n^{p+p^\prime+q+q^\prime}} \sum_{\substack{\mathbf{k}\otimes\mathbf{l}\otimes\mathbf{k}^\prime\otimes\mathbf{l}^\prime\\\in (\mathscr{V}_{p,q}\otimes\mathscr{V}_{p^\prime,q^\prime})\backslash \wt V(*,*,*,*,0)}}
\left[\mathbb{E} \int_{\mathcal{M}} \psi_{p+q}(x) \prod_{i=1}^p D_{h,1}(\alpha_i,X_{k_i},x) \prod_{j=1}^q D_{h,2}(\beta_j,X_{l_j},x)  d\mathscr{H}(x)\right.\nonumber\\
&\hspace{1cm} \left.\times  \mathbb{E} \int_{\mathcal{M}}  \psi_{p^\prime+q^\prime}(x^\prime)  \prod_{i=1}^{p^\prime} D_{h,1}(\alpha_i,X_{k_i^\prime},x^\prime)   \prod_{j=1}^{q^\prime} D_{h,2}(\beta_j,X_{l_j^\prime},x^\prime)d\mathscr{H}(x^\prime) \right]\nonumber\\
\leq &  \frac{1}{n}  C \rho_{n,h,1}^{p,q} \rho_{n,h,1}^{p^\prime,q^\prime}.\label{rateidiff}
\end{align}

The result above is derived by using (\ref{sumofrate}) and for $0\leq u\leq\lfloor(p+q)/2\rfloor$ and $0\leq u^\prime\leq\lfloor(p^\prime+q^\prime)/2\rfloor$,
\begin{align*}
&\frac{|\{\mathbf{k}\otimes\mathbf{l}\otimes\mathbf{k}^\prime\otimes\mathbf{l}^\prime\in (\mathscr{V}_{p,q}\otimes\mathscr{V}_{p^\prime,q^\prime})\backslash \wt V(*,*,*,*,0):\; \mathbf{k}\otimes\mathbf{l}\in V(s,t,u),\; \mathbf{k}^\prime\otimes\mathbf{l}^\prime \in V(s^\prime,t^\prime,u^\prime) \}|}{|V(s,t,u)|\times|V(s^\prime,t^\prime,u^\prime)|} \\
= &O\left(\frac{1}{n} \right).
\end{align*}

%

Next we consider $\frac{1}{n^{p+p^\prime}n^{q+q^\prime}} \int_{\mathcal{M}}\int_{\mathcal{M}} \text{II}_h(x,x^\prime)d\mathscr{H}(x)d\mathscr{H}(x^\prime)$. For $\mathbf{k}\otimes\mathbf{l}\otimes\mathbf{k}^\prime\otimes\mathbf{l}^\prime\in \wt V(s,t,s^\prime,t^\prime,u)$ with $u\geq1$, we have
\begin{align}
&\mathbb{E} \left\{\left[\int_{\mathcal{M}} \psi_{p+q}(x) \prod_{i=1}^p D_{h,1}(\alpha_i,X_{k_i},x)  \prod_{j=1}^q D_{h,2}(\beta_j,X_{l_j},x) d\mathscr{H}(x)\right] \right. \nonumber\\
&\hspace{1cm}\times\left.\left[\int_{\mathcal{M}} \psi_{p^\prime+q^\prime}(x^\prime) \prod_{i=1}^{p^\prime} D_{h,1}(\alpha_i,X_{k_i^\prime},x^\prime)  \prod_{j=1}^{q^\prime} D_{h,2}(\beta_j,X_{l_j^\prime},x^\prime) d\mathscr{H}(x^\prime)\right]\right\}\nonumber\\
=&\int_{\mathcal{M}} \psi_{p^\prime+q^\prime}(x^\prime) \int_{\mathcal{M}}\psi_{p+q}(x) \int_{\mathbb{R}^{d(s+t+s^\prime+t^\prime+u)}}\prod_{i=1}^p D_{h,1}(\alpha_i,y_{k_i},x)  \prod_{j=1}^q D_{h,2}(\beta_j,y_{l_j},x)    \nonumber  \\
&\hspace{1cm}\times  \prod_{i=1}^{p^\prime} D_{h,1}(\alpha_i,y_{k_i^\prime},x^\prime)   \prod_{j=1}^{q^\prime} D_{h,2}(\beta_j,y_{l_j^\prime},x^\prime) \prod_{\iota\in\text{uniq}(\{\mathbf{k}\}\cup\{\mathbf{l}\}\cup\{\mathbf{k}^\prime\}\cup\{\mathbf{l}^\prime\})} f(y_\iota)d\mathbf{y}  d\mathscr{H}(x)  d\mathscr{H}(x^\prime)\label{intermediatestep} \\
= &  h^{ud} h^{(s+t+s^\prime+t^\prime)\nu}h^{d-1} \int_{\cal M}   \psi_{p+q}(x) \int_{T_x(\mathcal{M})} \psi_{p^\prime+q^\prime}(x+hx^\prime) W(x,x^\prime)d\mathscr{H}(x^\prime)d\mathscr{H}(x) (1+o(1)),\nonumber
%
%
\end{align}

where $W$ is described below. The last expression is a result of the following steps: (i) $h^{ud}$ is obtained by applying change variables $z_\iota = (x-y_\iota)/h$ for index $\iota$ that has duplicates; (ii) $h^{(s+t+s^\prime+t^\prime)\nu}$ is obtained by applying change variables $z_\iota = (x-y_\iota)/h$ and using integration by part to exchange the derivatives between $K$ and $f$ for index $\iota$ that has no duplicates; and (iii) since duplicates in $\iota$ exists (because $u\geq 1$), $K(z_\iota + (x-x^\prime)/h)$ will appear after step (ii). Then $h^{d-1}$ is obtained using argument given in the proof of Theorem~\ref{asympnorm}, which approximates one of the surface integrals on $\mathcal{M}$ by the integral on the tangent space $T_x(\mathcal{M})$, and the error in the approximation is represented by the $o(1)$. After the above steps, the integration over $\mathbb{R}^{d(s+t+s^\prime+t^\prime+u)}$ in (\ref{intermediatestep}) is changed to $W(x,x^\prime)$ which satisfies $\sup_{x,x^\prime\in\mathbb{R}^d}|W(x,x^\prime)|\leq C \left(\frac{1}{h^{d+1}} \right)^{p+p^\prime-s-s^\prime} \left(\frac{1}{h^{d+2}} \right)^{q+q^\prime-t-t^\prime}$ for some constant $C$, where we have used the assumption that $K$ has bounded support. Therefore there exist a finite constant $C$ such that the above expectation is bounded by 
\begin{align*}
C  h^{ud} h^{(s+t+s^\prime+t^\prime)\nu}  \left(\frac{1}{h^{d+1}} \right)^{p+p^\prime-s-s^\prime} \left(\frac{1}{h^{d+2}} \right)^{q+q^\prime-t-t^\prime} h^{d-1}.
\end{align*}


Let $\eta_{n,h}(s,t,s^\prime,t^\prime) = n^{s+t+s^\prime+t^\prime+u} h^{(s+t+s^\prime+t^\prime)\nu +ud -(p+p^\prime-s-s^\prime)(d+1) - (q+q^\prime-t-t^\prime)(d+2) +(d-1)}.$ Note that $|\wt V(s,t,s^\prime,t^\prime,u)| =O(n^{s+t+s^\prime+t^\prime+u})$. 
Then 
\begin{align}\label{rateii}
&\frac{1}{n^{p+p^\prime}n^{q+q^\prime}} \int_{\mathcal{M}}\int_{\mathcal{M}} \text{II}_h(x,x^\prime)d\mathscr{H}(x)d\mathscr{H}(x^\prime) \nonumber\\
\leq & \frac{C}{n^{p+p^\prime}n^{q+q^\prime}}  \mathop{\sum_{s=0}^p \sum_{t=0}^q}_{s+t\leq p+q-1}   \mathop{\sum_{s^\prime=0}^{p^\prime} \sum_{t^\prime=0}^{q^\prime}}_{s^\prime+t^\prime\leq p^\prime+q^\prime-1} \sum_{u=1}^{\lfloor (p+p^\prime+q+q^\prime-s-t-s^\prime-t^\prime)/2 \rfloor}  \eta_{n,h}(s,t,s^\prime,t^\prime) \nonumber\\
= & C \rho_{n,h,2}^{p,q,p^\prime,q^\prime}, 
\end{align}
where $$\rho_{n,h,2}^{p,q,p^\prime,q^\prime} = \mathop{\sum_{s=0}^p \sum_{t=0}^q}_{s+t\leq p+q-1}   \mathop{\sum_{s^\prime=0}^{p^\prime} \sum_{t^\prime=0}^{q^\prime}}_{s^\prime+t^\prime\leq p^\prime+q^\prime-1} \sum_{u=1}^{\lfloor (p+p^\prime+q+q^\prime-s-t-s^\prime-t^\prime)/2 \rfloor} \frac{(nh^d)^uh^{(s+t+s^\prime+t^\prime)\nu + d-1}}{(nh^{d+1})^{p+p^\prime-s-s^\prime}(nh^{d+2})^{q+q^\prime-t-t^\prime}}.$$

%
%
%

Let $r_{n,h,4} = n^{-1}h^{(p+p^\prime+q+q^\prime)\nu-5-2\nu}$, $r_{n,h,5} = \frac{h^{(p+p^\prime)\nu+d-1}}{(nh^{d+4})^{(q+q^\prime)/2}}$, and $r_{n,h,6} =\frac{ h^{ d-1}}{(nh^{d+2})^{(p+p^\prime)/2}(nh^{d+4})^{(q+q^\prime)/2}}$. Then similar to the discussion in a), b), and c) above (\ref{sumelambda}), we can show that for some constant $C$
\begin{align*}
& \left| \frac{1}{n^{p+p^\prime}n^{q+q^\prime}} \int_{\mathcal{M}}\int_{\mathcal{M}} \text{II}_h(x,x^\prime)d\mathscr{H}(x)d\mathscr{H}(x^\prime) \right| \\
\leq &
\begin{cases}
C r_{n,h,4} & \text {if } nh^{d+4+2\nu} \rightarrow\gamma_1 \text { for } 1<\gamma_1\leq+\infty,\\
C r_{n,h,5} & \text {if } nh^{d+4+2\nu} \rightarrow\gamma_1 \; \& \; nh^{d+2+2\nu} \rightarrow \gamma_2  \text { for } 0\leq \gamma_1<1<\gamma_2\leq+\infty,\\
C r_{n,h,6} & \text {if } nh^{d+2+2\nu} \rightarrow \gamma_2 \text { for } 0\leq \gamma_2<1.
\end{cases}
\end{align*}

Therefore from (\ref{varhexp}), (\ref{hhdecomp}), (\ref{rateidiff}) and (\ref{rateii}), and noticing that $n^{-1}\rho_{n,h,1}^{p,q} \rho_{n,h,1}^{p^\prime,q^\prime} = O(\rho_{n,h,2}^{p,q,p^\prime,q^\prime}) $ we have $|\text{Cov}[\lambda(f,H_{n,h}^{p,q}),\; \lambda(f,H_{n,h}^{p^\prime,q^\prime})]| \leq C \rho_{n,h,2}^{p,q,p^\prime,q^\prime},$ and therefore 

$$|\text{Cov}[\lambda(f,H_{n,h}^{p,q}),\; \lambda(f,H_{n,h}^{p^\prime,q^\prime})] | \leq C ((\beta_{n,h,1}^{p+p^\prime+q+q^\prime}h^{d-1})\vee (\beta_{n,h,2}^{p+p^\prime+q+q^\prime}(nh^5)^{-1})),$$
where $\beta_{n,h,1} = 1/\sqrt{nh^d}$ and $\beta_{n,h,2}=h^\nu$.\\

Recall $\lambda(f,T_{n,j})$ is a sum of $(a_d)^j$ functions represented by $\frac{1}{j!}\lambda(f,H_{n,h}^{p,q})$ with $p+q=j$. So for $j,j^\prime\geq 1$, $$\left|\text{Cov} \left[\lambda(f,T_{n,j}),\lambda(f,T_{n,j^\prime})\right] \right| \leq \frac{1}{j!j^\prime !} [d+d(d+1)/2]^{j+j^\prime} C ((\beta_{n,h,1}^{j+j^\prime}h^{d-1}))\vee (\beta_{n,h,2}^{j+j^\prime}(nh^5)^{-1})).$$
Now we are ready to show (\ref{corremainder}). Notice that if we fix $j+j^\prime=k\geq 3$, then there are $k-2$ distinct ordered pairs of $[d_{\wh f}(x)^T - d_f(x)^T]^{\otimes j}$ and $[d_{\wh f}(x)^T - d_f(x)^T]^{\otimes j^\prime}$ in the double sum. So

\begin{align*}
%
&\left|\text{Var}\left[\sum_{j=2}^{Q_\nu-1}\lambda(f,T_{n,j})\right] + 2 \text{Cov} \left[\lambda(f,T_{n,1}), \sum_{j=2}^{Q_\nu-1}\lambda(f,T_{n,j})\right] \right|\\
\leq & 2\sum_{j^\prime=2}^{Q_\nu-1} \sum_{j=1}^{Q_\nu-1} \left|\text{Cov} \left[\lambda(f,T_{n,j}), \lambda(f,T_{n,j^\prime})\right] \right| \\
\leq & 2\sum_{j^\prime=2}^{Q_\nu-1} \sum_{j=1}^{Q_\nu-1}\left[ \frac{1}{j!} \frac{1}{j^\prime!} (a_d)^{j+j^\prime}C \Large((\beta_{n,h,1}^{j+j^\prime}h^{d-1})\vee (\beta_{n,h,2}^{j+j^\prime}(nh^5)^{-1})\Large)  \right] \\ 
\leq&  2\sum_{k=3}^{Q_\nu-1} \left[ (k-2)C (a_d)^{k} ((\beta_{n,h,1}^{k}h^{d-1})\vee (\beta_{n,h,2}^{k}(nh^5)^{-1})) \frac{1}{k!}\sum_{j+j^\prime = k}   \frac{k!}{j!j^\prime!}  \right] \\ 
= &  2\sum_{k=3}^{Q_\nu-1} \left[ (k-2)C(a_d)^{k} ((\beta_{n,h,1}^{k}h^{d-1})\vee (\beta_{n,h,2}^{k}(nh^5)^{-1})) \frac{2^k}{k!}  \right] \\
%
%
& = O(((\beta_{n,h,1}^3h^{d-1})\vee (\beta_{n,h,2}^3(nh^5)^{-1}))) = o((nh^5)^{-1}).
\end{align*}
Hence we have shown (\ref{corremainder}) and the proof for the case of $Q\geq 3\vee(\nu-1)$ is completed by combining this with (\ref{qnv}), (\ref{tn1asymp}), (\ref{asymptoticmu1exp}), (\ref{asympototicmu2exp}) and (\ref{biasexpression}). The proof for the case of $Q > 3\vee(\nu-1)$ follows from similar calculation above. $\hfill\square$
\end{proof}
%
%

{\bf Proof of Corollary~\ref{asympototicnormalityspec}}
\begin{proof}
Write the Taylor expansion $\wt\phi(d_{\wh f,1}^*(x))-\wt\phi(d_{f,1}(x)) = \sum_{p=1}^\nu T_{n,p}^*(x),$ where 
\begin{align*}
T_{n,p}^*(x) = 
\begin{cases}
\frac{1}{p!} \left[d_{\wh f,1}^*(x)^T -d_{f,1}(x)^T \right]^{\otimes p} \; \nabla^{\otimes p} \wt\phi(d_{f,1}(x))  ,&  j=1,\cdots,\nu-1,\\
\frac{1}{p!} \left[d_{\wh f,1}^*(x)^T -d_{f,1}(x)^T \right]^{\otimes p} \; \nabla^{\otimes p} \wt\phi(\wt{d_{f,1}}(x)), & j=\nu,
\end{cases}
\end{align*}
where $\wt{d_{f,1}}$ is between $d_{f,1}$ and $d_{\wh f,1}^*$. Note that
\begin{align*}
\wh g(x) =  \frac{(n-Q)!}{n!} \sum_{(l_1,\cdots,l_Q)\in\pi_n^Q} \left\{ \left[\wt\phi(d_{f,1}(x)) + \sum_{p=1}^\nu T_{n,p}^*(x)\right]\prod_{j=1}^Q \left[D_{h,2}(\beta_j,X_{l_j},x) + f^{(\beta_j)}(x)\right] \right\} 
\end{align*}
and we can write $\wh g(x) -g(x) = g_{n,1}(x) + g_{n,2}(x)$ where
\begin{align*}
g_{n,1}(x) = & \sum_{p=1}^\nu a_{n,p}(x):= \sum_{p=1}^\nu \left\{ \frac{(n-Q)!}{n!} \sum_{(l_1,\cdots,l_Q)\in\pi_n^q}  T_{n,p}^*(x)\prod_{j=1}^Q \left[D_{h,2}(\beta_j,X_{l_j},x) + f^{(\beta_j)}(x)\right] \right\}, \\
g_{n,2}(x)= & \frac{(n-Q)!}{n!} \wt\phi(d_{f,1}(x)) \sum_{(l_1,\cdots,l_q)\in\pi_n^Q} \left\{\prod_{j=1}^Q \left[D_{h,2}(\beta_j,X_{l_j},x) + f^{(\beta_j)}(x)\right] -  \prod_{j=1}^Q f^{(\beta_j)}(x) \right\}.
\end{align*}

Note that $a_{n,p}(x)$, $p=1,2,\cdots,\nu-1$ above 
is a sum of $d^p$ functions in the generic form of
\begin{align*}
%
\frac{1}{p!}\frac{(n-Q)!}{n!} \frac{1}{(n-Q)^p} \psi_{p}(x) \sum_{\mathbf{k}\otimes\mathbf{l}\in\mathscr{V}_{p,Q}} \left\{ \prod_{i=1}^p D_{h,1}(\alpha_i,X_{k_i},x) \prod_{j=1}^Q \left[D_{h,2}(\beta_j,X_{l_j},x) + f^{(\beta_j)}(x)\right] \right\}, 
%
\end{align*}
where $\sup_{x\in\mathcal{M}} |\psi_p(x)| \leq C$ and $$\mathscr{V}_{p,Q}=\{\mathbf{k}\otimes\mathbf{l}: (l_1,\cdots,l_Q) \in\pi_n^Q, k_1,\cdots,k_p\in\{1,\cdots,n\}\backslash \{l_1,\cdots,l_Q\} \}.$$
The structure in the above expressions of $g_{n,1}(x)$ and $g_{n,2}(x)$ is similar to $\phi(d_f(x))$ studied in Theorem~\ref{asympototicnormalitymknown}. Using arguments similar to that proof, we can show that $\mathbb{E}[\lambda(f,g_{n,1})] = O(h^\nu + \frac{1}{nh^{d+2}})$ and $\text{Var}[\lambda(f,g_{n,1})] = O((nh^3)^{-1})$. Also $\mathbb{E}[\lambda(f,g_{n,2})] = O(h^\nu)$ and $\sqrt{nh^5}\lambda(f, g_{n,2}-\mathbb{E}g_{n,2}) \rightarrow_D N(0,\sigma^2)$ for some $\sigma^2$. Here $\sigma^2$ is determined by the linear part of $\lambda(f, g_{n,2}-\mathbb{E}g_{n,2})$, which is
\begin{align*}
& \int_{\mathcal{M}}\frac{(n-Q)!}{n!} \wt\phi(d_{f,1}(x)) \sum_{(l_1,\cdots,l_Q)\in\pi_n^Q} \left\{\sum_{j=1}^Q \left[D_{h,2}(\beta_j,X_{l_j},x) \times \prod_{j^\prime=1,j^\prime\neq j}^Q f^{(\beta_{j^\prime})}(x)\right]  \right\}  d\mathscr{H}(x)\\
= &  \int_{\mathcal{M}}\wt\phi(d_{f,1}(x))\sum_{j=1}^Q  \left\{ \frac{1}{n}\sum_{i=1}^n \left[ D_{h,2}(\beta_j,X_{i},x) \times \prod_{j^\prime=1,j^\prime\neq j}^Q f^{(\beta_{j^\prime})}(x)\right]\right\}d\mathscr{H}(x).
\end{align*}
So calculation following (\ref{asyptoticvarexp}) leads to
\begin{align}\label{sigma2express}
&\sigma^2 = c \int_{\mathcal{M}} \wt\phi(d_{f,1}(x))^2 \int_{\mathbb{R}} \left( \int_{T_{x,t}(\mathcal{M})}  \sum_{j=1}^Q\prod_{j^\prime=1,j^\prime\neq j}^Q \left(f^{(\beta_{j^\prime})}(x) K^{(\beta_{j})}(u)  d\mathscr{H}(u)\right) \right)^2dt d\mathscr{H}(x).
%
%
\end{align}
$\hfill\square$
\end{proof}

{\bf Proof of Theorem~\ref{finaltheorem}}
\begin{proof}
Due to the decomposition of $\wt\lambda_{\epsilon_n}(\wh f,\wh g) - \lambda(f,g)$ in (\ref{lambdadiffdecomp}), (\ref{asymptotunrat}) and (\ref{asymptotnormunkn}) follow from Theorem~\ref{HausdorffIntegration}, Corollary~\ref{Hausdorff} and Theorem~\ref{asympototicnormalitymknown}. When Corollary~\ref{Hausdorff} is applied, notice that if $Q\geq1$ then $\alpha_n=\gamma_{n,h}^{(2)} + h^\nu$, $\beta_n=\gamma_{n,h}^{(3)} + h^{\nu-1}$ and $\eta_n=\gamma_{n,h}^{(4)} + h^{\nu-2}$; while if $Q=0$ then $\alpha_n=\gamma_{n,h}^{(1)} + h^\nu$, $\beta_n=\gamma_{n,h}^{(2)} + h^{\nu-1}$ and $\eta_n=\gamma_{n,h}^{(3)} + h^{\nu-2}$.\\

To get (\ref{asymptotnormunknpl}), it suffices to show $\wh\sigma_l^2-\sigma_l^2=o_p(1)$ for $l=1,2$. Let $q_l(x,u)=\nabla_l\phi(d_f(x))^T d_{K,l}(u)$ and $\wh q_l(x,u)=\nabla_l\phi(d_{\wh f}(x))^T d_{K,l}(u)$. Then $\sigma_l^2 = c\int_{\mathcal{M}} \int_{\mathbb{R}} \left(\int_{T_{x,t}(\mathcal{M})} q_l(x,u)d\mathscr{H}(u) \right)^2dtdd\mathscr{H}(x)$ and $\wh\sigma_l^2 = c\int_{\wh{\mathcal{M}}} \int_{\mathbb{R}} \left(\int_{T_{x,t}(\wh{\mathcal{M}})} \wh q_l(x,u)d\mathscr{H}(u) \right)^2dt d\mathscr{H}(x)$.
Using (\ref{supvar0}) - (\ref{supbias2}), we have $$\sup_{x\in\mathcal{I}(\delta_0)}\sup_{v\in\mathbb{R}^{d}} |\wh q_l(x,v)- q_l(x,v)|=o_p(1),$$ which yields $\wh\sigma_l^2 - \wt\sigma_l^2=o_p(1)$, where $\wt\sigma_l^2 = c\int_{\wh{\mathcal{M}}} \int_{\mathbb{R}} \left(\int_{T_{x,t}(\wh{\mathcal{M}})} q_l(x,u)d\mathscr{H}(u) \right)^2dt d\mathscr{H}(x).$ We only need to show $\wt\sigma_l^2-\sigma_l^2=o_p(1)$. 
%
%
%
We use the finite atlas $\{(U_\alpha,\psi_\alpha):\alpha\in\mathscr{A}\}$ defined in the proof of Theorem~\ref{HausdorffIntegration} and again we suppress the subscript $\alpha$. Recall $P_n$ defined in (\ref{map}) and $\wh U =\{P_n(x):x\in U\}$. Following the proof of of Theorem~\ref{HausdorffIntegration}, we can write 
\begin{align*}
&\int_{\wh{U}} \int_{\mathbb{R}} \left(\int_{T_{x,t}(\wh{\mathcal{M}})} q_l(x,u)d\mathscr{H}(u) \right)^2dt d\mathscr{H}(x) - \int_{U} \int_{\mathbb{R}} \left(\int_{T_{x,t}(\mathcal{M})} q_l(x,u)d\mathscr{H}(u) \right)^2dtdd\mathscr{H}(x) \\
= & \text{I}_n + \text{II}_n+\text{III}_n,    
\end{align*}
where
\begin{align*}
&\text{I}_n = \int_{U} \int_{\mathbb{R}} \left[\left(\int_{T_{P_n(x),t}(\wh{\mathcal{M}})} q_l(x,u)d\mathscr{H}(u) \right)^2 - \left(\int_{T_{x,t}(\mathcal{M})} q_l(x,u)d\mathscr{H}(u) \right)^2 \right]dt d\mathscr{H}(x),\\
%
&\text{II}_n = \int_{U} \int_{\mathbb{R}} \left[\left(\int_{T_{P_n(x),t}(\wh{\mathcal{M}})} q_l(P_n(x),u)d\mathscr{H}(u) \right)^2 - \left(\int_{T_{P_n(x),t}(\wh{\mathcal{M}})} q_l(x,u)d\mathscr{H}(u) \right)^2 \right]dt d\mathscr{H}(x),\\
%
&\text{III}_n = \int_{U} \left\{\{\text{det}[\mathbf{I}_d-Q_n(\psi^{-1}(x))]\}^{1/2}-1\right\} \int_{\mathbb{R}}\left(\int_{T_{P_n(x),t}(\wh{\mathcal{M}})} q_l(P_n(x),u)d\mathscr{H}(u) \right)^2 dt d\mathscr{H}(x),
\end{align*}
where $Q_n$ is given in (\ref{qndef}). Following similar arguments as given in the proof of Theorem~\ref{HausdorffIntegration}, it can be seen that $\text{II}_n=o_p(1)$ and $\text{III}_n=o_p(1)$. Notice that for any $x\in \mathcal{M}$, $\wh \n(P_n(x))-\n(x)=o_p(1)$ and so $\text{I}_n=o_p(1)$. Hence $\wt\sigma_l^2-\sigma_l^2=o_p(1)$ and we conclude the proof. $\hfill\square$
\end{proof}


{\bf Proof of Theorem~\ref{eulercharactheorem}}

\begin{proof}
Lemma~\ref{compatible} implies that there exists $\epsilon_0>0$ such that for all $|\epsilon|\leq\epsilon_0$, ${\cal M}_c$, $\wh{\cal M}_c$, $\wh{\cal M}_{c+\epsilon}$, and $\wh{\cal M}_c\uplus\epsilon$ are all homeomorphic for $n$ large enough with probability one. The result in the theorem is a consequence of (\ref{eulerestimator1}), (\ref{eulerestimator2}), and the fact the Euler characteristic is invariant under homeomorphisms. $\hfill\square$
\end{proof}


\section*{References}
\begin{description}
\itemsep0em 
\item Ambrosio, L., Colesanti, A. and Villa, E. (2008). Outer Minkowski content for some classes of closed sets. {\em Math. Ann.} {\bf 342}, 727-748.
\item Arias-Castro, E., Mason, D., and Pelletier, B. (2016). On the estimation of the gradient lines of a density and the consistency of the mean-shift algorithm. {\em Journal of Machine Learning Research} {\bf 17} 1-28.
\item Arias-Castro, E., and Rodr\'{i}guez-Casal, A. (2017). On estimating the perimeter using the alpha-shape. {\em Ann. Inst. H. Poincar\'{e} Probab. Statist.} {\bf 53} 1051-1068.
\item Armend\'{a}riz, I., Cuevas, A., and Fraiman, R. (2009). Nonparametric estimation of boundary measures and related functionals: asymptotic results. {\em Advances in Applied Probability} {\bf 41} 311-322.
\item Baddeley, A. J. and Gundersen, H. J. G. and Cruz-Orive, L. M. (1986). Estimation of surface area from vertical sections. {\em Journal of Microscopy} {\bf 142} 259-276.
\item Baddeley, A.J. and Jensen, E.B.V. (2005). {\em Stereology for Statisticians}, CRC Press, Boca Raton, FL.
\item Ba\'{i}llo, A. (2003). Total error in a plug-in estimator of level sets. {\em Statistics \& Probability Letters} {\bf 65} 411-417.
\item Ba\'{i}llo, A., Cuevas, A. and Justel, A. (2000). Set estimation and nonparametric detection. {\em The Canadian Journal of Statistics} {\bf 28} 765-782.
%
%
\item Bickel, P.J. and Ritov, Y. (1988). Estimating integrated squared density derivatives: sharp best order of convergence estimates. {\em Sankhy\={a} Ser. A} {\bf 50} 381-393.
\item Bobrowski, O., Mukherjee, S. and Taylor, J.E. (2017). Topological consistency via kernel estimation. {\em Bernoulli}, {\bf 23}, 288 - 328.
\item Cadre, B. (2006). Kernel estimation of density level sets. {\em J. Multivariate Anal.} {\bf 97} 999-1023.
\item Cadre, B., Pelletier, B. and Pudlo, P. (2009). Clustering by estimation of density level sets at a fixed probability. Available at \url{http://hal.archives-ouvertes.fr/docs/00/39/74/37/PDF/tlevel.pdf}.
\item Cannings, T.I., Berrett, T.B., Samworth, R.J. (2017). Local nearest neighbour classification with applications to semi-supervised learning. {\em arXiv: 1704.00642}.
\item Canzonieri, V. and Carbone, A. (1998). Clinical and biological applications of image analysis in non-Hodgkin's lymphomas. {\em Hematological Oncology} {\bf 16} 15-28.
\item Caselles, V., Haro, G., Sapiro, G., Verdera, J. (2008). On geometric variational models for inpainting surface holes. {\em Computer Vision and Image Understanding}. {\bf 111} 351-373.
\item Chac\'{o}n J.E. and Duong T. (2018). {\em Multivariate Kernel Smoothing and Its Applications}. Chapman \& Hall/CRC, Boca Raton (USA).
\item Chazal, F., Lieutier, A., and Rossignac, J. (2007). Normal-map between normal-compatible manifolds. {\em International Journal of Computational Geometry \& Applications} {\bf 17} 403-421.
\item Chen, Y.-C. (2017). Nonparametric inference via bootstrapping the debiased estimator. {\em arXiv: 1702.07027}.
\item Chen, Y.-C., Genovese, C.R., and Wasserman, L. (2017). Density Level Sets: Asymptotics, Inference, and Visualization. {\em J. Amer. Statist. Assoc.}, {\bf 112} 1684-1696.
%
%
\item Cuevas, A. and Fraiman, R., and Rodr\'{i}guez-Casal, A. (2007). A nonparametric approach to the estimation of lengths and surface areas. {\em Ann. Statist.} {\bf 35} 1031-1051.
\item Cuevas, A. and Fraiman, R., and Pateiro-L\'{o}pez, B. (2012). On statistical properties of sets fulfilling rolling-type conditions. {\em Advances in Applied Probability} {\bf 44} 311-329.
\item Cuevas, A., Gonz\'{a}lez-Manteiga, W., and  Rodr\'{i}guez-Casal, A. (2006). Plug-in estimation of general level sets. {\em Australian \& New Zealand Journal of Statistics} {\bf 48} 7-19.
%
%
\item Dony J. and Mason D.M. (2008). Uniform in bandwidth consistency of conditional U-statistics. {\em Bernoulli} {\bf14} 1108-1133.
\item Evans, L.C. and Gariepy, R.F. (1992). {\em Measure Theory and Fine Properties of Functions}. CRC Press, Boca Raton, FL.
\item Federer, H. (1959). Curvature measures. {\em Trans. Amer. Math. Soc.} {\bf 93} 418-491.
\item Garcia-Dorado, D., Th{\'e}roux, P., Duran, J. M., Solares, J., Alonso, J., Sanz, E., Munoz, R., Elizaga, J., Botas, J. and Fernandez-Avil{\'e}s, F. (1992). Selective inhibition of the contractile apparatus. A new approach to modification of infarct size, infarct composition, and infarct geometry during coronary artery occlusion and reperfusion. {\em Circulation} {\bf 85} 1160-1174.
\item Gin\'{e}, E. and Nickl, R. (2008). A simple adaptive estimator of the integrated square of a density. {\em Bernoulli} {\bf 14} 47-61.
\item Gin\'{e}, E. and Mason, D. (2008). Uniform in bandwidth estimation of integral functionals of the density function. {\em Scandinavian Journal of Statistics} {\bf 35} 739-761.
\item Gokhale, A. M. (1990). Unbiased estimation of curve length in 3-D using vertical slices. {\em Journal of Microscopy} {\bf 159}, 133-141.
\item Goldman R. (2005). Curvature formulas for implicit curves and surfaces. {\em Computer Aided Geometric Design} {\bf 22}, 632-658.
\item Gray, A. (2004). {\em Tubes}. 2nd ed. {\em Progress in Mathematics} {\bf 221}. Birkh\"{a}user, Basel.
%
%
\item Hall, P. and Kang, K.-H. (2005). Bandwidth choice for nonparametric classification. {\em Ann. Statist.} {\bf 33} 284-306.
\item Hall, P. and Marron, J.S. (1987). Estimation for integrated squared density derivatives. {\em Statist. Probab. Lett.} {\bf6} 109-155.
\item Hall, P. and Murison, R.D. (1993). Correcting the negativity of high-order kernel density estimators. {\em J. Multivariate Anal.} {\bf 47} 103-122.
\item Hartigan, J.A. (1987). Estimation of a convex density contour in two dimensions. {\em J. Amer. Statist. Assoc.} {\bf 82} 267-270.
\item Ipsen, I.C.F. and Rehman, R. (2008). Perturbation bounds for determinants and characteristic polynomials. {\em SIAM J. Matrix Anal. Appl.} {\bf30}(2), 762-776.
%
%
\item Jim\'{e}nez, R. and Yukich, J. E. (2011). Nonparametric estimation of surface integrals. {\em Ann. Statist.} {\bf 39} 232-260.
\item Kerscher, M. (2000). Statistical analysis of large-scale structure in the Universe. In {\em Statistical Physics and Spatial Statistics: The Art of Analyzing and Modeling Spatial Structures and Pattern Formation} (ed. Mecke, K.R. \& Stoyan, D.). Lecture Notes in Physics, vol. 554. Springer.
\item Kindlmann, G., Whitaker, R., Tasdizen, T. and M\"{o}ller, T. (2003). Curvature-based transfer functions for direct volume rendering: methods and applications. In {\em Proc. IEEE Visualization 2003}, pages 513-520.
%
%
\item Lee, J.M. (1997). {\em Riemannian Manifolds An Introduction to Curvature}. Springer-Verlag New York.
\item Levit, B. Ya. (1978). Asymptotically efficient estimation of nonlinear functionals [in Russian]. {\em Problems Inform. Transmission} {\bf14}, 65-72.
\item Lim, P.H., Bagci, U, and Bai, L. (2013) Introducing Willmore flow Into level set segmentation of spinal vertebrae. {\em IEEE Transactions on Biomedical Engineering} {\bf 60}(1), 115-122.
\item  Mammen, E. and Polonik, W. (2013). Confidence sets for level sets. {\em Journal of Multivariate Analysis} {\bf 122} 202-214.
\item Mammen, E., and Tsybakov, A. B. (1999). Smooth Discrimination Analysis. {\em Ann. Statist.} {\bf27} 1808-1829.
\item Mart\'{i}nez, V. J., and Saar, E. (2001). {\em Statistics of the Galaxy Distribution}. Chapman \& Hall/CRC, Roca Raton, Florida.
\item Mason, D.M. and Polonik, W. (2009). Asymptotic normality of plug-in level set estimates. {\em The Annals of Applied Probability} {\bf 19} 1108-1142.
%
%
\item Mecke, K. R., Buchert, T., and Wagner, H. (1994): Robust morphological measures for large-scale structure in the Universe. {\em Astron. Astrophys.}, {\bf 288}, 697-704.
%
%
\item Mirsky, L. (1955). {\em An introduction to Linear Algebra}. Clarendon Press, Oxford.
%
%
\item Pateiro-L\'{o}pez, B. and Rodr\'{i}guez-Casal, A. (2008). Length and surface area estimation under smoothness restrictions. {\em Advances in Applied Probability} {\bf 40} 348-358.
\item Pateiro-L\'{o}pez, B. and Rodr\'{i}guez-Casal, A. (2009). Surface area estimation under convexity type assumptions. {\em Journal of Nonparametric Statistics} {\bf 21} 729-741.
\item Polonik, W. (1995). Measuring mass concentrations and estimating density contour clusters - an excess mass approach. {\em Ann. Statist.} {\bf 23} 855-881.
\item Pranav, P., Edelsbrunner, H., van de Weygaert, R., Vegter, G., Kerber, M., Jones, B. J. T.,  Wintraecken, M. (2017). The topology of the cosmic web in terms of persistent Betti numbers. {\em Mon. Not. R. Astron. Soc.}, {\bf 465}(4), 4281-4310. 
\item Qiao, W. (2018). Asymptotics and optimal bandwidth selection for nonparametric estimation of density level sets. {\em arXiv: 1707.09697}.
%
%
%
\item Qiao, W. and Polonik, W. (2019). Nonparametric confidence regions for level sets: statistical properties and geometry.  {\em Electronic Journal of Statistics} 13(1), 985-1030.
\item Rigollet, P. and Vert, R. (2009). Optimal rates for plug-in estimators of density level sets. {\em Bernoulli} {\bf 15} 1154-1178.
\item Rinaldo, A. and Wasserman, L. (2010). Generalized density clustering. {\em Ann. Statist.} {\bf 38} 2678-2722.
\item Salas, W.A., Boles, S. H., Frolking, S., Xiao, X. and Li, C. (2003). The perimeter/area ratio as an index of misregistration bias in land cover change estimates. {\em International Journal of Remote Sensing} {\bf 24} 1165-1170.
\item Samworth, R.J. (2012). Optimal weighted nearest neighbour classifiers. {\em Ann. Statist.} {\bf 40} 2733-2763.
\item Schmalzing, J. and G\'{o}rski, K. M. (1998). Minkowski functionals used in the morphological analysis of cosmic microwave background anisotropy maps. {\em Mon. Not. Roy. Astron. Soc.}, {\bf297} 355-365.
\item Seifert, U. (1997). Configurations of fluid membranes and vesicles. {\em Advances in Physics}, {\bf 46}(1),13-137.
\item Silverman, B.W. (1986). {\em Density Estimation for Statistics and Data Analysis}. Chapman and Hall, London.
\item Singh, A., Scott C. and Nowak, R. (2009). Adaptive Hausdorff estimation of density level sets. {\em Ann. Statist.} {\bf 37} 2760-2782.
\item Steinwart, I. (2015). Fully adaptive density-based clustering. {\em Ann. Statist.} {\bf 43} 2132-2167.
\item Steinwart, I., Hush, D. and Scovel, C. (2005). A classification framework for anomaly detection. {\em J. Machine Learning Research} {\bf 6} 211-232.
\item Tsybakov, A.B. (1997): Nonparametric estimation of density level sets. {\em Ann. Statist.} {\bf 25} 948-969.
\item Trillo, N.G., Slep\v{c}ev, D. and Von Brecht, J. (2017). Estimating perimeter using graph cuts. {\em Advances in Applied Probability} {\bf 49} 1067-1090.
\item Troutt, M.D., Pang, W.K., and Hou, S.H. (2004). {\em Vertical Density Representation and Its Applications}. World Scientific, River Edge, NJ.
\item Turner, K., Mukherjee, S., and Boyer, D.M. (2014). Persistent homology transform for modeling shapes and surfaces. {\em Information and Inference: A Journal of the IMA}. {\bf3}(4) 310-344.
\item Walther, G. (1997). Granulometric smoothing.  {\em Ann. Statist.} {\bf } 2273-2299.
%
%
\item Weyl, H. (1939). On the volume of tubes. {\em Amer. J. Math.} {\bf 61} 461-472.
\item Willmore, T.J. (1965). Note on embedded surfaces. {\em An. \c Sti. Univ. ``Al. I. Cuza'' Ia\c{s}i Sec\c{t}. I a Mat. (N.S.).} {\bf 11B} 493-496
\end{description}
\end{document}